\documentclass[final,onefignum,onetab]{environment/siamart251216}

\usepackage{environment/obara_siam_basic}
\usepackage{environment/global_convergence_specific}

\ifpdf
\hypersetup{
  pdftitle={Riemannian Interior Point Trust Region Method},
  pdfauthor={M. Obara, T. Okuno, and A. Takeda}
}
\fi

\makeatletter
\newcommand{\customfootnote}[2]{%
  \xdef\@thefnmark{}\@footnotetext{$\!\!\!$\textsuperscript{#1}#2}
}
\makeatother

\begin{document}
    \maketitle

    \customfootnote{$\dagger$}{Graduate School of Information Science and Technology, The University of Tokyo, Tokyo, Japan (\email{polyelo2plumo@gmail.com,
    takeda@mist.i.u-tokyo.ac.jp})}
    \customfootnote{$\ddagger$}{Center for Advanced Intelligence Project, RIKEN, Tokyo, Japan}
    \customfootnote{$\S$}{Faculty of Science and Technology, Seikei University, Tokyo, Japan (\email{takayuki-okuno@st.seikei.ac.jp})}

    \begin{abstract}
        We consider Riemannian inequality-constrained optimization problems.
        Such problems inherit the benefits of Riemannian approach developed in the unconstrained setting and naturally arise from applications in control, machine learning, and other fields.
        We propose a Riemannian primal-dual interior point trust region method (\RIPTRM{}) for solving them.
        We prove its global convergence to an approximate Karush-Kuhn-Tucker point and a weak second-order stationary point.
        Under the strict complementarity condition, this result reduces to global convergence to a second-order stationary point.
        To the best of our knowledge, this is the first algorithm that incorporates the trust region strategy for constrained optimization on Riemannian manifolds, and has the second-order convergence property for optimization problems on Riemannian manifolds with nonlinear inequality constraints.
        We conduct numerical experiments in which we introduce a truncated conjugate gradient method and an eigenvalue-based subsolver for \RIPTRM{} to approximately and exactly solve the trust region subproblems, respectively.
        Empirical results show that 
        \RIPTRM{}s consistently find solutions with high accuracy.
        \isextendedversion{Additionally, we observe that \RIPTRM{} with the exact search direction shows promising performance in an instance where the Hessian of the Lagrangian has a large negative eigenvalue.}{}
    \end{abstract}
    
    \begin{keywords}
        Riemannian optimization, Inequality-constrained optimization, Interior point trust region method, Eigenvalue-based solver.
    \end{keywords}
    
    \begin{AMS}
      65K05, 90C30
    \end{AMS}    

    \section{Introduction}
        In this paper, we consider the following optimization problem:
        \isextendedversion{
        \begin{mini}
            {\pt \in \mani}{\objfun\paren*{\pt}}
            {\label{prob:RICO}}{}
            \addConstraint{\ineqfun[\ineqidx]\paren*{\pt}}{\geq 0,}{\quad \ineqidx \in \ineqset\coloneqq\brc*{1,\ldots,\ineqdime},}
        \end{mini}
        }{
        \begin{equation}\label{prob:RICO}
            \underset{\pt \in \mani} {\minimize} ~ \objfun\paren*{\pt} \quad \subjectto ~ \ineqfun[\ineqidx]\paren*{\pt}\geq 0, \quad \ineqidx \in \ineqset\coloneqq\brc*{1,\ldots,\ineqdime},
        \end{equation}}
        where $\mani$ is a $\dime$-dimensional, connected, complete Riemannian manifold, and $\objfun\colon\mani\to\setR$ and $\brc*{\ineqfun[\ineqidx]}_{\ineqidx\in\ineqset}\colon\mani\to\setR[]$ are twice continuously differentiable functions.
        If $\mani$ is not connected, we can consider our result on the connected components of $\mani$ separately.
        We call problem \cref{prob:RICO} the Riemannian inequality-constrained optimization problem and abbreviate it as the \RICO{} problem. 
        \RICO{} is a natural extension of the standard inequality-constrained optimization problem from a Euclidean space to a Riemannian manifold. 
        Indeed, when $\mani=\setR[\dime]$, \RICO{}~\eqref{prob:RICO} reduces to the standard problem on $\setR[\dime]$.
        Due to its versatility, \RICO{}~\cref{prob:RICO} has applications across various fields~\cite{Obaraetal2024StabLinSysIdentifwithPriorKnwlbyRiemSQO,ZassShashua2006NNPCA,BrossettEscandeKheddar2018MCPostureComputonMani}.

        The Riemannian approach offers several advantages over its Euclidean approach.
        First, the Riemannian modeling ensures the feasibility for every iterate being in $\mani$, while
        their feasibility of iterates and the output in Euclidean optimization depends on the behavior of numerical algorithms.
        This may lead to an infeasible output~ \cite[Section 5.1]{Andreanietal2024GlobConvALMforNLPviaRiemOptim} and slow convergence to accurate solutions~\cite[Section~5]{LiuBoumal2020SimpleAlgoforOptimonRiemManiwithCstr} in the constrained settings.
        Second, the Riemannian approach serves as a means to handle specific structures in applications, such as the set of positive-definite matrices and the geometry of quotient manifolds~\cite{Obaraetal2024StabLinSysIdentifwithPriorKnwlbyRiemSQO,MisawaSato2022H2OptimReductofPosNetusingRiemALM,Satoetal20RiemIdentifSymPosSys,SatoSato17RiemSysIdentifLinearMIMO,SatoSato17StructurePreservByRTR}, the special orthogonal group~\cite{BrossettEscandeKheddar2018MCPostureComputonMani,Boumal23IntroOptimSmthMani}, the Grassmann manifold~\cite{Bendokatetal2024GrassmannManiHandBookBasicGeomandComputAsp}, the fixed-rank manifold~\cite{Vandereycken2013LowRankMatComplbyRiemOptim} and the flag manifold~\cite{Yeetal2022OptimonFlagMani}.
        Furthermore, certain possibly nonsmooth sets, such as bounded-rank matrices~\cite{KhrulkovOseledets2018DesingofBddRankMatSets,Levinetal2023FindStatPtsonBddrankMatsaGeomHurSmthRemedy,RebjockBoumal2024OptimOverBddRankMatsThruADesingEnabJointGlobandLocalGuar}, can be addressed via a smooth parameterization~\cite{Levinetal2024EffofSmthParamonNonconvOptimLandsc}.
        These special structures often come together with additional constraints; for example, the references~\cite{Obaraetal2024StabLinSysIdentifwithPriorKnwlbyRiemSQO,MisawaSato2022H2OptimReductofPosNetusingRiemALM,BrossettEscandeKheddar2018MCPostureComputonMani} study constrained optimization on Riemannian manifolds.
        Third, leveraging the inherent geometry of optimization problems provides better understanding and improves computational performance.
        For example, the formulations of geodesically convex optimization problems
        guarantee that any local optimum is also a global optimum~\cite{Boumal23IntroOptimSmthMani}.
        The formulation covers the geodesically convex constraints, which arise in optimistic likelihood~\cite{Nguyenetal2019CalcOptimLUsingGeodsicConvOptim} and the Riemannian barycenter~\cite[Section 6.4]{Hsiehetal2023RiemStochOptimMethAvoidStrSaddlePt}, for example.
        
        Riemannian optimization has been extensively investigated over the last two decades especially for the unconstrained case.
        Absil~et al.~\cite{Absiletal08OptimAlgoonMatMani} established the modern theory of Riemannian optimization, where they proposed the geometric Newton method and the Riemannian trust region methods.
        Based on their work, classical algorithms for the Euclidean optimization have been extended to the Riemannian settings with the guarantees of the global convergence to a first-order stationary point~\cite{HuangAbsilGallivan2015RiemSymRankOneTRM,Huangetal18RiemBFGSMethwithoutDiffRetrforNonconvOptim,Narushimaetal2023MemlessQuasiNewtonMethBsdonSpectrScalingBroydenFamforRiemOptim,SakaiIiduka2024ModifMemlessSpectrScalingBroydenFamonRiemMani,Sato2022RiemConjGradMethGenFwandSpecAlgowithConvAnal}.
        Several algorithms are designed to achieve the global convergence to a second-order stationary point (\SOSP{}) and efficiently escape the saddle points~\cite{Agarwaletal2021AdaptRegwithCubiconMani,ZhangZhang2018CubicRegNewtonMethoverRiemMani,BoumalAbsilCartis2018GlobRateofConvforNonconvOptimonMani,GoyensRoyer2024RiemTRMethforStrictSaddleFunwithComplexGuar,KasaiMishra2018InexactTRAlgoonRiemMani,CriscitielloBoumal2019EffEscSaddlePtonMani,SunFlammarionFazel2019EscfromSaddlePtonRiemMani}.

        Riemannian optimization with nonlinear constraints is also under development.
        Yang~et~al.~\cite{YangZhangSong14OptimCondforNLPonRiemMani} derived optimality conditions for Riemannian optimization with inequality and equality constraints.
        Constraint qualifications have been examined for the Riemannian cases~\cite{BergmannHerzog2019IntrinsicFormulationKKTcondsConstrQualifSmthMani,YamakawaSato2022SeqOptimCondforNLOonRiemManiandGlobConvALM,Andreanietal2024ConstrQualifStrongGlobConvPropALMonRiemMani,Andreanietal2024GlobConvALMforNLPviaRiemOptim}.
        Liu and Boumal~\cite{LiuBoumal2020SimpleAlgoforOptimonRiemManiwithCstr} developed an exact penalty method (\EPM{}) combined with smoothing and an augmented Lagrangian method (\ALM{}).
        They proved the global convergence of \EPM{} to a Karush-Kuhn-Tucker (\KKT{}) point and that of \ALM{} to a \KKT{} point, an \SOSP{} in the equality-constrained case, or the global optimum depending on the subsolver's quality.
        Note that, however, the concrete analyses of the sequence by the subsolver are beyond the scope of their work, and computing the global optimum is generally hard due to nonconvexity.
        \ALM{} and \EPM{} have been refined to achieve the global convergence to an approximate-\KKT{} (\AKKT{}) point~\cite{YamakawaSato2022SeqOptimCondforNLOonRiemManiandGlobConvALM,SmthL1ExactPenaMethforIntrinsicConstRiemOptimProb} or a positive-\AKKT{} (\PAKKT{}) point~\cite{Andreanietal2024ConstrQualifStrongGlobConvPropALMonRiemMani,Andreanietal2024GlobConvALMforNLPviaRiemOptim}.
        Obara~et~al.~\cite{ObaraOkunoTakeda2021SQOforNLOonRiemMani} proposed Riemannian sequential quadratic optimization (\SQO{}) for optimization problems with inequality and equality constraints, proving its global convergence to a \KKT{} point and the local quadratic convergence to a local optimum under assumptions. 
        Schiela and Ortiz~\cite{SchielaOrtiz2021SQPMethforEqCstrOptimonHilbertMani} proposed \SQO{} for optimization problems on Hilbert manifolds with equality constraints and proved the local quadratic convergence. 
        
        Interior point method, abbreviated as \IPM{}, is known as one of the most prominent algorithms for optimization problems with nonlinear constraints in the Euclidean setting; we refer readers to \cite{NocedalWright2006NumerOptim,Forsgrenetal2002IntMethforNonlinOptim,GouldOrbanToint2005NumerMethforLgeScaleNLO} for comprehensive reviews of \IPM{}s in the Euclidean setting.
        Interior point trust region method, or \IPTRM{}, is a variant of \IPM{} that uses the trust region approach as the globalization mechanism; that is, it employs the trust region scheme when a trial iterate fails to make sufficient progress toward the solution set and is rejected~\cite[Chapter 13]{ConnGouldToint2000TRMeth}.
        The strengths of \IPTRM{} include strong convergence properties and the practicality.
        For instance, W\"{a}chter and Biegler~\cite{WumlautachterBiegler2000FailGlobalConvforClofIntPtMethforNLP} provided an example where several \IPM{}s using line search fail to achieve feasibility when starting at reasonable infeasible points, while the \IPTRM{} by Byrd et al.~\cite{ByrdGilbertNocedal2000TRMethBsdIntPtTechniqueforNLP} avoids this issue.
        Additionally, the trust region framework allows for the use of the exact Hessian of the Lagrangian when computing the search directions by solving subproblems, whereas the line search methods cannot since the quadratic term in the subproblems need to be positive semidefinite.
        In the Euclidean setting, primal-dual \IPTRM{}s are designed to have global convergence to \SOSP{}s~\cite{Connetal2000PrimalDualTRAlgoforNonconvexNLP,Tseng2002ConvInfeasiIntptTRMethforConstMin}.
        Various \IPTRM{}s are also proposed and proved to possess global convergence to \KKT{} points and local convergence~\cite{ByrdGilbertNocedal2000TRMethBsdIntPtTechniqueforNLP,ByrdLiuNocedal1997LocalBehavIntPtMethforNLP,YamashitaYabeTanabe2005GlobalSuperlinConvPrimalDualIntPtTRNethforLgeScaleCstrOptim,UlbrichUlbrichVicente2004GlobalConvPrimalDualIntPtFltrMethforNLP,Silvaetal2008GlobalConvPrimalDualIntPtFltrMethforNLPNewFltrOptimMeasComputRslt,Gouldetal2001SuperlinConvofPDIntptAlgoforNonlinProgram}.

        In the Riemannian setting, Lai and Yoshise~\cite{LaiYoshise2024RiemIntPtMethforCstrOptimonMani} proposed two types of \IPM{}s that have global convergence to a \KKT{} point and the local quadratic convergence properties, respectively.
        Hirai~et~al.~\cite{HiraiNieuwboerWalter2023IntptMethonManiTheoandAppl} analyzed \IPM{} for convex optimization problems on Riemannian manifolds.
        Note that these \IPM{}s are based on the line search for the globalization mechanism.
        To the best of our knowledge, no \IPTRM{}s have been developed in the Riemannian setting. 
        Building on prior studies in the Euclidean setting~\cite[Chapters 4 and 19]{NocedalWright2006NumerOptim}, the trust region approach is expected to exhibit strong convergence properties and high practicality for constrained optimization on Riemannian manifolds.

    \subsection{Our contribution}
        In this paper, we propose a primal-dual Riemannian \IPTRM{}, abbreviated as \RIPTRM{}, for \RICO{}~\cref{prob:RICO}, inspired by a Euclidean \IPTRM{}~\cite{Connetal2000PrimalDualTRAlgoforNonconvexNLP}.
        Compared with \cite{Connetal2000PrimalDualTRAlgoforNonconvexNLP}, we employ the exact Hessian rather than its approximation and omit the treatment of linear equality constraints, focusing on extending the \IPTRM{} to the Riemannian setting.
        \RIPTRM{} consists of outer and inner iterations.
        In outer iteration, starting from an interior point of \RICO{}~\cref{prob:RICO}, \RIPTRM{} generates a sequence by adjusting the barrier parameter and the tolerance level for residuals used in the inner iterations.
        In the inner iteration, the algorithm computes a search direction by approximately or exactly solving a trust region subproblem defined on a tangent space of $\mani$.
        The subproblem uses the full Hessian of $\objfun$ and $\brc*{\ineqfun[\ineqidx]}_{\ineqidx\in\ineqset}$, which may be indefinite.
        Then, the search direction is evaluated to decide if the algorithm updates the iterate and adjusts the trust region radius.
        Here, we make use of the retraction, a smooth map for computing the next iterate in $\mani$ according to the search direction, and the log barrier function as a merit function.
        
        Our contributions are summarized as follows:
        \begin{enumerate}
            \item 
            We propose \RIPTRM{} for solving \RICO{}~\cref{prob:RICO}.
            To the best of our knowledge, this is the first algorithm incorporating the trust region strategy for optimization problems on Riemannian manifolds with nonlinear constraints.
            \item We prove that our \RIPTRM{} achieves the global convergence to an \AKKT{} point and a weak \SOSP{} (\wSOSP{}), 
            depending on the chosen quality of the search direction.
            Under the strict complementarity condition, this result reduces to global convergence to an \SOSP{}.
            To the best of our knowledge, this is the first algorithm achieving the second-order global convergence property for Riemannian optimization with inequality constraints.
            
            The \SOSP{} is a stronger solution concept than a \KKT{} point, in that it excludes strict saddle points.
            In Euclidean constrained optimization, the vanilla projected gradient method, 
            which is globally convergent only to \KKT{} points,
            may converge to a strict saddle point with positive probability~\cite{Nouiehedetal2018Convto2ndOrdStnrtyforConstNonConvOptim}.
            Since \RICO{} generalizes the Euclidean setting, such problematic instances fall within our framework.
            From this perspective, guaranteeing global convergence to an \SOSP{} is beneficial.
            \item
            We conduct numerical experiments in which we introduce a truncated conjugate gradient (\tCG{}) method and eigenvalue-based subsolver~\cite{Adachietal2017SolvingTRSbyGenEigenProb} for \RIPTRM{} to approximately and exactly solve the trust region subproblems, respectively.
            Empirical results show that 
            \RIPTRM{}s consistently find solutions with high accuracy.
            \isextendedversion{Additionally, we observe that \RIPTRM{} with the exact search direction shows promising performance in an instance where the Hessian of the Lagrangian has a large negative eigenvalue.}{}
        \end{enumerate}
        
        Note also that, combined with a companion work~\cite{Obaraetal2025LocalConvofRiemIPMs} on a local convergence of Riemannian \IPM{}s, our method is the first Riemannian \IPM{} that possesses both global and local convergence.
        We summarize the comparison of our algorithm with the existing ones for Riemannian optimization with nonlinear constraints in \cref{tabl:RIPTRMcomparison}.

        \begin{table}[]
        \caption{Comparison of algorithms for Riemannian optimization with nonlinear constraints. 
        The symbols $\ineqset$ and $\eqset$ denote the capability of handling inequality and equality constraints, respectively.
        The symbols $\doublecheck$ and $\doublecheck[ptred]$ represent the global convergence to an \AKKT{} and a \PAKKT{} point, respectively. 
        }
        
        \begin{center}    
        \begin{tabular}{cccccc}
        \hline
        Reference                                                                                                                      & Method & Constraints        & \KKT{} & \wSOSP{}$^{\dagger}$ 
        & Local    \\ \hline
        \cite[Algorithm 1]{LiuBoumal2020SimpleAlgoforOptimonRiemManiwithCstr}                                                                      & ALM    & $\ineqset, \eqset$ & $\checkmark$                & $\checkmark^{\ddagger}$  &          \\
        \cite[Algorithm 2]{LiuBoumal2020SimpleAlgoforOptimonRiemManiwithCstr}                                                                      & EPM    & $\ineqset, \eqset$ & $\checkmark$                &                                           &          \\
        \cite{SchielaOrtiz2021SQPMethforEqCstrOptimonHilbertMani}                                                                     & SQO    & $\eqset$           &                    &                                           & quadratic      \\
        \cite{ObaraOkunoTakeda2021SQOforNLOonRiemMani}                                                                                & SQO    & $\ineqset, \eqset$ & $\checkmark$                &                                           & quadratic      \\
        \cite{YamakawaSato2022SeqOptimCondforNLOonRiemManiandGlobConvALM} & ALM    & $\ineqset, \eqset$ & $\doublecheck $        &                                           &          \\
        \cite[Algorithm 2]{LaiYoshise2024RiemIntPtMethforCstrOptimonMani}                                                                          & IPM    & $\ineqset, \eqset$ &                    &                                           & quadratic      \\
        \cite[Algorithm 5]{LaiYoshise2024RiemIntPtMethforCstrOptimonMani}                                                                          & IPM    & $\ineqset, \eqset$ & $\checkmark$                &                                           &          \\
        \cite{Andreanietal2024ConstrQualifStrongGlobConvPropALMonRiemMani,Andreanietal2024GlobConvALMforNLPviaRiemOptim}                                                                          & ALM    & $\ineqset, \eqset$ & $\doublecheck[ptred]$                &                                           &          \\
        \cite{SmthL1ExactPenaMethforIntrinsicConstRiemOptimProb}                                                                          & EPM    & $\ineqset, \eqset$ & $\doublecheck $                &                                           &          \\
       \textbf{Ours}
        & \textbf{\IPTRM}  & \textbf{$\ineqset$}         & \textbf{$\doublecheck$}                & \textbf{$\checkmark$}                                    & $\text{near-quadratic}^{\S}$ \\ \hline
        \end{tabular}
        \end{center}
        
        \footnotesize{
        $\dagger$: Under the strict complementarity condition, the results reduce to global convergence to \SOSP{}s.\\
        $\ddagger$: \ALM{} in \cite{LiuBoumal2020SimpleAlgoforOptimonRiemManiwithCstr} has the second-order convergence property under only equality-constrained setting.\\
        $\S$: The local convergence analysis can be found in \cite{Obaraetal2025LocalConvofRiemIPMs}.
        }\label{tabl:RIPTRMcomparison}
        \end{table}

    \subsection{Paper organization}
        In \cref{sec:preliminaries}, we review fundamental concepts from Riemannian geometry and Riemannian optimization.
        In \cref{sec:propmeth}, we describe \RIPTRM{} consisting of outer and inner iterations.
        We prove its global convergence in \cref{sec:globconv}.
        In \cref{sec:experiment}, we provide numerical experiments on the stable linear system identification and compare our algorithm with the existing methods.
        In \cref{sec:conclusion}, we summarize our research and state future work.


    \section{Preliminaries}\label{sec:preliminaries}
        Define $\setNz\coloneq\brc*{0, 1, 2, \ldots}$ and 
        let $\setR[\dime], \setRp[\dime]$, and $\setRpp[\dime]$ be $\dime$-dimensional Euclidean space, its nonnegative orthant, and its positive orthant, respectively. 
        We denote by $\onevec[\ineqdime]$ the $\ineqdime$-dimensional vector of ones. 
        We omit the subscript $\ineqdime$ when the context is clear. 
        A continuous function $\forcingfun\colon\setRp[]\to\setRp[]$ is said to be a forcing function if $\forcingfun\paren*{\barrparam[]}=0$ holds if and only if $\barrparam[]=0$.
        \isextendedversion{
        For related positive quantities $\orderconstone$ and $\orderconsttwo$, we write $\orderconstone=\bigO[\orderconsttwo]$ if there exists a constant $\ordercoeff > 0$ such that $\orderconstone\leq\ordercoeff\orderconsttwo$ for all $\orderconsttwo$ sufficiently small.
        }{}
        Given two normed vector spaces $\vecspc, \vecspctwo$ and a linear operator $\opr\colon\vecspc\to\vecspctwo$, we define the operator norm as $\opnorm{\opr}\coloneqq\sup\brc*{\norm{\opr\vecone}_{\vecspctwo} \colon \vecone \in \vecspc \text{ and } \norm{\vecone}_{\vecspc} \leq 1}$, where $\norm{\plchold}_{\vecspc}$ and $\norm{\plchold}_{\vecspctwo}$ are the norms on $\vecspc$ and $\vecspctwo$, respectively.
    
    \subsection{Notation and terminology from Riemannian geometry}\label{subsec:notationRiemggeo}
        We briefly review some concepts from Riemannian geometry by following the notation of \cite{Absiletal08OptimAlgoonMatMani,Boumal23IntroOptimSmthMani}.
        Let $\ptone[] \in \mani$, and let $\tanspc[{\ptone[]}]\mani$ be the tangent space to $\mani$ at $\ptone[]$.
        A vector field on $\mani$ is a map $\funthr\colon\mani\to\tanspc[]\mani$ with $\funthr\paren*{\ptone}\in\tanspc[\ptone]\mani$, where $\tanspc[]\mani$ is the tangent bundle. 
        Let $\intvl\subseteq\setR[]$ be an open interval, and let $\curve\colon\intvl\to\mani$ be a smooth curve.
        A Riemannian metric on $\mani$ is a choice of inner product $\metr[{\ptone[]}]{\plchold}{\plchold}\colon\tanspc[{\ptone[]}]\mani\times\tanspc[{\ptone[]}]\mani\to\setR[]$ for every $\ptone\in\mani$ satisfying that, for all smooth vector fields $\funthr, \funfou$ on $\mani$, the function $\ptone\mapsto\metr[{\ptone[]}]{\funthr\paren*{\ptone}}{\funfou\paren*{\ptone}}$ is smooth from $\mani$ to $\setR[]$.
        A Riemannian manifold is a smooth manifold endowed with a Riemannian metric.
        The Riemannian metric induces the norm $\Riemnorm[\ptone]{\tanvecone[\ptone]} \coloneqq \sqrt{\metr[{\ptone[]}]{\tanvecone[\ptone]}{\tanvecone[\ptone]}}$ for $\tanvecone[\ptone]\in\tanspc[\ptone]\mani$ and the Riemannian distance $\Riemdist{\plchold}{\plchold}\colon\mani\times\mani\to\setR$.
        It follows from \cite[Theorem 13.29]{Lee12IntrotoSmthManibook2ndedn} that $\mani$ is a metric space under the Riemannian distance.
        According to the Hopf-Rinow theorem, every closed bounded subset of $\mani$ is compact for a finite-dimensional, connected, complete Riemannian manifold by regarding $\mani$ as a metric space~\cite[Chapter 7, Theorem 2.8]{doCarmo92RiemGeo}. 

        For two smooth manifolds $\mani_{1}, \mani_{2}$ and a differentiable map $\funone\colon\mani_{1}\to\mani_{2}$, we denote by $\D\funone\paren*{\ptone[]}\colon\tanspc[\ptone]\mani_{1}\to\tanspc[\funone\paren*{\ptone}]\mani_{2}$ the differential of $\funone$ at $\ptone\in\mani_{1}$.
        We use the canonical identification $\tanspc[\ptone]\vecspc\cid\vecspc$ for a vector space $\vecspc$ and $\ptone\in\vecspc$.
        Let $\funset\paren*{\mani}$ denote the set of sufficiently differentiable real-valued functions.
        Given $\funtwo\in\funset\paren*{\mani}$, $\D\funtwo\paren*{\ptone}\sbra*{\tanvecone[\ptone]} \in \tanspc[\funtwo\paren*{\ptone}]\setR[]\cid\setR$ is the differential of $\funtwo$ at $\ptone\in\mani$ along $\tanvecone[\ptone]\in\tanspc[\ptone]\mani$.
        The Riemannian gradient of $\funtwo$ at $\ptone$, denoted by $\gradstr\funtwo\paren*{\ptone}$, is defined as a unique element of $\tanspc[\ptone]\mani$ that satisfies
        \isextendedversion{\begin{equation}\label{eq:riemgrad}
                \metr[\ptone]{\gradstr\funtwo\paren*{\ptone}}{\tanvecone[\ptone]} = \D\funtwo\paren*{\ptone}\sbra*{\tanvecone[\ptone]}    
            \end{equation}}{
            $\metr[\ptone]{\gradstr\funtwo\paren*{\ptone}}{\tanvecone[\ptone]} = \D\funtwo\paren*{\ptone}\sbra*{\tanvecone[\ptone]}$}
        for any $\tanvecone[\ptone]\in\tanspc[\ptone]\mani$.
        Here, $\gradstr\funtwo\colon\mani\to\tanspc[]\mani\colon\ptone\mapsto\gradstr\funtwo\paren*{\ptone}$ is the gradient vector field.
        Let $\Riemcxt{}{}$ be the Levi-Civita connection and $\vecfldset\paren*{\mani}$ be the set of sufficiently differentiable vector fields. 
        For any $\funthr\in\vecfldset\paren*{\mani}$, we define the Jacobian of $\funthr$ at $\ptone$ as $\Jacobian[\funthr]\paren*{\ptone}\colon\tanspc[\ptone]\mani\to\tanspc[\ptone]\mani\colon\tanvecone[\ptone]\mapsto\Riemcxt{\tanvecone[\ptone]}{\funthr}$.
        In particular, for the case $\funthr=\gradstr\funtwo$, the operator $\Hess\funtwo\paren*{\ptone}\colon\tanspc[\ptone]\mani\to\tanspc[\ptone]\mani$ denotes the Riemannian Hessian of $\funtwo$ at $\ptone$; that is, $\Hess\funtwo\paren*{\ptone}\sbra*{\tanvecone[\ptone]}\coloneqq\Riemcxt{\tanvecone[\ptone]}{\gradstr\funtwo}$ for all $\tanvecone[\ptone]\in\tanspc[\ptone]\mani$.
        When $\mani$ is a Euclidean space, we have $\Hess\funtwo\paren*{\ptone}\sbra*{\tanvecone[\ptone]} = \D\paren*{\gradstr\funtwo}\paren*{\ptone}\sbra*{\tanvecone[\ptone]}$
        and $\metr[\ptone]{\Hess\funtwo\paren*{\ptone}\sbra*{\tanvecone[\ptone]^{1}}}{\tanvecone[\ptone]^{2}} = \D[2]\funtwo\paren*{\ptone}\sbra*{\tanvecone[\ptone]^{1}, \tanvecone[\ptone]^{2}}$ for all $\tanvecone[\ptone], \tanvecone[\ptone]^{1}, \tanvecone[\ptone]^{2} \in \tanspc[\ptone]\mani$.
        For each $\pt[1]\in\mani$ and any $\pt[2]\in\mani$ sufficiently close to $\pt[1]$, we denote by $\partxp[]{\pt[2]}{\pt[1]}\colon\tanspc[{\pt[1]}]\mani\to\tanspc[{\pt[2]}]\mani$ the parallel transport of $\tanspc[{\pt[1]}]\mani$ to $\tanspc[{\pt[2]}]\mani$ along the unique minimizing geodesic that connects $\pt[1]$ and $\pt[2]$.
        \isextendedversion{
        Note that the parallel transport is isometric, i.e., $\Riemnorm[{\pt[2]}]{\partxp[]{\pt[2]}{\pt[1]}\sbra*{\tanvecone[{\pt[1]}]}} = \Riemnorm[{\pt[1]}]{\tanvecone[{\pt[1]}]}$ for any $\tanvecone[{\pt[1]}]\in\tanspc[{\pt[1]}]\mani$ and that $\partxp[]{\pt[1]}{\pt[1]}\colon\tanspc[{\pt[1]}]\mani\to\tanspc[{\pt[1]}]\mani$ is the identity map.
        The adjoint of the parallel transport corresponds with its inverse, that is, $\metr[{\pt[2]}]{\partxp[]{\pt[2]}{\pt[1]}\sbra*{\tanvecone[{\pt[1]}]}}{\tanvectwo[{\pt[2]}]} = \metr[{\pt[1]}]{\tanvecone[{\pt[1]}]}{\partxp[]{\pt[1]}{\pt[2]}\sbra*{\tanvectwo[{\pt[2]}]}}$ for any $\tanvecone[{\pt[1]}]\in\tanspc[{\pt[1]}]\mani$ and any $\tanvectwo[{\pt[2]}]\in\tanspc[{\pt[2]}]\mani$.}{}
        
        A retraction $\retr[]\colon\tanspc[]\mani\to\mani$ is a smooth map with the following properties: by letting $\retr[\ptone]\colon\tanspc[\ptone]\mani\to\mani$ be the restriction of $\retr[]$ to $\tanspc[\ptone]\mani$, it satisfies 
        \isextendedversion{
            \begin{subequations}\label{eq:retrdef}
            \begin{align}
                \retr[\ptone]\paren*{\zerovec[\ptone]} = \ptone,\label{eq:retrzero}\\
                \D\retr[\ptone]\paren*{\zerovec[\ptone]} = \id[{\tanspc[\ptone]\mani}]\label{eq:retrdiffzero}
            \end{align}
            \end{subequations}}{
            $\retr[\ptone]\paren*{\zerovec[\ptone]} = \ptone \text{ and }
            \D\retr[\ptone]\paren*{\zerovec[\ptone]} = \id[{\tanspc[\ptone]\mani}]$}
        under $\tanspc[{\zerovec[\ptone]}]\paren*{\tanspc[\ptone]\mani}\cid\tanspc[\ptone]\mani$, where $\zerovec[\ptone]$ is the zero vector of $\tanspc[\ptone]\mani$ and $\id[{\tanspc[\ptone]\mani}]$ denotes the identity map on
        \isextendedversion{
        \begin{equation}\label{def:pullbackfun}
            \pullback[]{\funtwo}\coloneqq\funtwo\circ\retr[] \text{ and } \pullback[\ptone]{\funtwo}\coloneqq\funtwo\circ\retr[\ptone]
        \end{equation}
        }{$\tanspc[\ptone]\mani$.
        Let $\pullback[]{\funtwo}\coloneqq\funtwo\circ\retr[]$ and $\pullback[\ptone]{\funtwo}\coloneqq\funtwo\circ\retr[\ptone]$}
        denote the pullback of the function $\funtwo\colon\mani\to\setR[]$ and the restriction of $\pullback{\funtwo}$ to $\tanspc[\ptone]\mani$, respectively.
        \isextendedversion{Note that, it follows from \cref{eq:retrdef} that
        }{Note that, it follows from the definition of the retraction that}
        \begin{equation}\label{eq:pullbackgradzero}
            \gradstr\pullback[]{\funtwo}\paren*{\zerovec[\ptone]} = \gradstr\funtwo\paren*{\ptone}.
        \end{equation}
        A second-order retraction is a retraction satisfying that, for all $\ptone[]\in\mani$ and all $\tanvecone[\ptone]\in\tanspc[\ptone]\mani$, the curve $\curve\paren*{\tmefiv}=\retr[\ptone]\paren*{\tmefiv\tanvecone[\ptone]}$ has zero acceleration at $\tmefiv=0$.
        Second-order retractions are not a restrictive class.
        For example, metric projection retractions and the exponential maps meet the condition; see \cite[Section~5.12]{Boumal23IntroOptimSmthMani} and \cite[Section~5.5]{Absiletal08OptimAlgoonMatMani} for details.
        Second-order retractions have the following property:
        \begin{proposition}[{\cite[Proposition~5.5.5]{Absiletal08OptimAlgoonMatMani}}]\label{prop:secodrretrHess}
            If the retraction $\retr[]$ is second order, then, for any $\funtwo\in\funset\paren*{\mani}$,
            \isextendedversion{
            \begin{align}\label{eq:secondretrHesszero}
                \Hess\pullback[\ptone]{\funtwo}\paren*{\zerovec[\ptone]} = \Hess\funtwo\paren*{\ptone},
            \end{align}
            }{$\Hess\pullback[\ptone]{\funtwo}\paren*{\zerovec[\ptone]} = \Hess\funtwo\paren*{\ptone}$,}
            where the left-hand side is the Hessian of $\pullback[\ptone]{\funtwo}\colon\tanspc[\ptone]\mani\to\setR$ at $\zerovec[\ptone]\in\tanspc[\ptone]\mani$.
        \end{proposition}
      
        We say that $\pullback[]{\funtwo}$ is radially Lipschitz continuously differentiable (\rLCone{}) on $\subsetmani\subseteq\mani$ if there exist positive scalars $\betaRL[\funtwo], \deltaRL[\funtwo] \in \setRpp$ such that, for any $\ptone\in\subsetmani$ and all $\tmetwo \geq 0, \tanvecone[\ptone]\in\tanspc[\ptone]\mani$ with $\tmetwo \Riemnorm[\ptone]{\tanvecone[\ptone]} \leq \deltaRL[\funtwo]$, it holds that
        \begin{equation}\label{def:radialLCone}
            \abs*{\metr[\ptone]{\gradstr\pullback[\ptone]{\funtwo}\paren*{\tmetwo\tanvecone[\ptone]} - \gradstr\pullback[\ptone]{\funtwo}\paren*{\zerovec[\ptone]}}{\tanvecone[\ptone]}} \leq \betaRL[\funtwo]\tmetwo \Riemnorm[\ptone]{\tanvecone[\ptone]}^{2}.
        \end{equation}
        Note that the definition above is equivalent to that in \cite[Definition~7.4.1]{Absiletal08OptimAlgoonMatMani} with the trivial case of $\tanvecone[\ptone] = \zerovec[\ptone]$.
        Similarly, we say that $\pullback[]{\funtwo}$ is radially Lipschitz twice continuously differentiable (\rLCtwo{}) on $\subsetmani\subseteq\mani$ if there exist positive scalars $\betaRLtwo[\funtwo], \deltaRLtwo[\funtwo] \in \setRpp[]$ such that, for any $\ptone\in\subsetmani$ and all $\tmetwo \geq 0, \tanvecone[\ptone]\in\tanspc[\ptone]\mani$ with $\tmetwo \Riemnorm[\ptone]{\tanvecone[\ptone]} \leq \deltaRLtwo[\funtwo]$, it holds that
        \begin{align}\label{def:radialLCtwo}
            \abs*{\metr[\ptone]{\paren*{\Hess\pullback[\ptone]{\funtwo}\paren*{\tmetwo\tanvecone[\ptone]} - \Hess\pullback[\ptone]{\funtwo}\paren*{\zerovec[\ptone]}}\sbra*{\tanvecone[\ptone]}}{\tanvecone[\ptone]}} \leq \betaRLtwo[\funtwo]\tmetwo \Riemnorm[\ptone]{\tanvecone[\ptone]}^{3}
        \end{align}
        under $\tanspc[{\tmetwo\tanvecone[\ptone]}]\paren*{\tanspc[\pt]\mani}\simeq\tanspc[\pt]\mani$ and $\tanspc[{\zerovec[\pt]}]\paren*{\tanspc[\pt]\mani}\simeq\tanspc[\pt]\mani$.

    \subsection{Optimality conditions for \RICO{}}\label{subsec:optimcond}
        Define the Lagrangian of \RICO{}~\cref{prob:RICO} as
        \isextendedversion{
        \begin{align}
            \Lagfun\paren*{\allvar[]} \coloneqq \objfun\paren*{\pt[]} - \sum_{\ineqidx\in\ineqset} \ineqLagmult[\ineqidx]\ineqfun[\ineqidx]\paren*{\pt[]},
        \end{align}}{
        $\Lagfun\paren*{\allvar[]} \coloneqq \objfun\paren*{\pt[]} - \sum_{\ineqidx\in\ineqset} \ineqLagmult[\ineqidx]\ineqfun[\ineqidx]\paren*{\pt[]}$,}
        where $\allvar[] \coloneqq \paren*{\pt[], \ineqLagmult[]} \in \mani \times \setR[\ineqdime]$ and $\ineqLagmult[]\in\setR[\ineqdime]$ is the vector of Lagrange multipliers for the inequality constraints.
        The Riemannian gradient and the Riemannian Hessian of the Lagrangian with respect to $\pt[]\in\mani$ are represented as 
        \isextendedversion{
        \begin{align}
            &\gradstr[\pt]\Lagfun\paren*{\allvar[]} = \gradstr\objfun\paren*{\pt[]} - \sum_{\ineqidx\in\ineqset} \ineqLagmult[\ineqidx]\gradstr\ineqfun[\ineqidx]\paren*{\pt[]}, \label{eq:gradLagfundef}\\
            &\Hess[\pt]\Lagfun\paren*{\allvar[]} = \Hess\objfun\paren*{\pt[]} - \sum_{\ineqidx\in\ineqset}\ineqLagmult[\ineqidx]\Hess\ineqfun[\ineqidx]\paren*{\pt[]},
        \end{align}
        }{
        $\gradstr[\pt]\Lagfun\paren*{\allvar[]} = \gradstr\objfun\paren*{\pt[]} - \sum_{\ineqidx\in\ineqset} \ineqLagmult[\ineqidx]\gradstr\ineqfun[\ineqidx]\paren*{\pt[]}, \text{ and } \Hess[\pt]\Lagfun\paren*{\allvar[]} = \Hess\objfun\paren*{\pt[]} - \sum_{\ineqidx\in\ineqset}\ineqLagmult[\ineqidx]\Hess\ineqfun[\ineqidx]\paren*{\pt[]}$,
        }
        respectively.
        For each $\paren*{\pt,\ineqLagmult[]}\in\mani\times\setR[\ineqdime]$, any $\vecone\in\setR[\ineqdime]$, and all $\tanvecone[\ptone]\in\tanspc[\ptone]\mani$, we define 
        \isextendedversion{
        \begin{align}
            &\Ineqfunmat[]\paren*{\pt} \coloneqq \diag\paren*{\ineqfun[]\paren*{\pt}}\in\setR[\ineqdime\times\ineqdime], \quad \IneqLagmultmat[] \coloneqq \diag\paren*{\ineqLagmult[]}\in\setR[\ineqdime\times\ineqdime],\\
            &\ineqgradopr[\pt]\sbra*{\vecone} \coloneqq \sum_{\ineqidx\in\ineqset} \vecone[\ineqidx] \gradstr\ineqfun[\ineqidx]\paren*{\pt}\in\tanspc[\pt]\mani,\\ &\coineqgradopr[\pt]\sbra*{\tanvecone[\ptone]} \coloneqq \trsp{\paren*{\metr[\pt]{\gradstr\ineqfun[1]\paren*{\pt}}{\tanvecone[\ptone]}, \ldots,\metr[\pt]{\gradstr\ineqfun[\ineqdime]\paren*{\pt}}{\tanvecone[\ptone]}}}\in\setR[\ineqdime],
        \end{align}}{
        \begin{align}
            &\Ineqfunmat[]\paren*{\pt} \coloneqq \diag\paren*{\ineqfun[]\paren*{\pt}}, \;
            \IneqLagmultmat[] \coloneqq \diag\paren*{\ineqLagmult[]}, \;
            \ineqgradopr[\pt]\sbra*{\vecone} \coloneqq \sum_{\ineqidx\in\ineqset} \vecone[\ineqidx] \gradstr\ineqfun[\ineqidx]\paren*{\pt}\in\tanspc[\pt]\mani,\\ &\coineqgradopr[\pt]\sbra*{\tanvecone[\ptone]} \coloneqq \trsp{\paren*{\metr[\pt]{\gradstr\ineqfun[1]\paren*{\pt}}{\tanvecone[\ptone]}, \ldots,\metr[\pt]{\gradstr\ineqfun[\ineqdime]\paren*{\pt}}{\tanvecone[\ptone]}}}\in\setR[\ineqdime],
        \end{align}
        }
        where $\diag\colon\setR[\ineqdime]\to\setR[\ineqdime\times\ineqdime]$ is the diagonal operator.
        Let \isextendedversion{\begin{align}
            \feasirgn \coloneqq \brc*{\pt[] \in \mani \relmiddle{|} \ineqfun[\ineqidx]\paren*{\pt[]} \geq 0 \text{ for all } \ineqidx\in\ineqset} \text{ and } \strictfeasirgn \coloneqq \brc*{\pt[] \in \mani \relmiddle{|} \ineqfun[\ineqidx]\paren*{\pt[]} > 0 \text{ for all } \ineqidx\in\ineqset}
        \end{align}}{
        $\feasirgn \coloneqq \brc*{\pt[] \in \mani \relmiddle{|} \ineqfun[\ineqidx]\paren*{\pt[]} \geq 0 \text{ for all } \ineqidx\in\ineqset} \text{ and } \strictfeasirgn \coloneqq \brc*{\pt[] \in \mani \relmiddle{|} \ineqfun[\ineqidx]\paren*{\pt[]} > 0 \text{ for all } \ineqidx\in\ineqset}$
        }
        be the feasible region and the strictly feasible region of \RICO{}~\cref{prob:RICO}, respectively.
        We define the index set of active inequalities at $\pt[] \in \mani$ as 
        \isextendedversion{
        \begin{equation}\label{def:activeineqset}
            \activeineqset\paren*{\pt[]} \coloneqq \brc*{\ineqidx \in \ineqset \relmiddle{|} \ineqfun[\ineqidx]\paren*{\pt[]} = 0}.
        \end{equation}}{$\activeineqset\paren*{\pt[]} \coloneqq \brc*{\ineqidx \in \ineqset \relmiddle{|} \ineqfun[\ineqidx]\paren*{\pt[]} = 0}$.}
        
        We introduce the optimality conditions and related concepts.       
        \begin{definition}[{\cite[Equation~(4.3)]{YangZhangSong14OptimCondforNLPonRiemMani}}]
            The linear independence constraint qualification (\LICQ{}) holds at $\ptaccum\in\mani$ if $\brc*{\gradstr\ineqfun[\ineqidx]\paren*{\ptaccum}}_{\ineqidx\in\activeineqset\paren*{\ptaccum}}$ is linearly independent on $\tanspc[\ptaccum]\mani$.
        \end{definition}

        \begin{theorem}[{\cite[Theorem~4.1]{YangZhangSong14OptimCondforNLPonRiemMani}}]
            Suppose that $\ptaccum \in \mani$ is a local minimum of \RICO{}~\cref{prob:RICO} and the \LICQ{} holds at $\ptaccum$.
            Then, there exists a vector of Lagrange multipliers for the inequality constraints $\ineqLagmultaccum[]\in\setR[\ineqdime]$ such that the following hold:
            \isextendedversion{
            \begin{equation}\label{eq:KKTconditions}
                \gradstr[\pt]\Lagfun\paren*{\allvaraccum} = \zerovec[\ptaccum], \ineqLagmultaccum[\ineqidx] \geq 0, \, \ineqfun[\ineqidx]\paren*{\ptaccum} \geq 0, \text{ and } \ineqLagmultaccum\ineqfun[\ineqidx]\paren*{\ptaccum} = 0 \text{ for all } \ineqidx \in \ineqset.
            \end{equation}
            We call \cref{eq:KKTconditions}}{$\gradstr[\pt]\Lagfun\paren*{\allvaraccum}$ $= \zerovec[\ptaccum]$, $\ineqLagmultaccum[\ineqidx] \geq 0$, $\ineqfun[\ineqidx]\paren*{\ptaccum} \geq 0$, and $\ineqLagmultaccum\ineqfun[\ineqidx]\paren*{\ptaccum} = 0$ for all $\ineqidx \in \ineqset$.
            We call them}
            the \KKT{} conditions of \RICO{}~\cref{prob:RICO} and $\ptaccum$ a \KKT{} point of \RICO{}~\cref{prob:RICO}.
        \end{theorem}
        We also introduce a sequential optimality condition that works without any constraint qualification:
        \begin{theorem}[{\cite[Definition~5, Theorem~1]{YamakawaSato2022SeqOptimCondforNLOonRiemManiandGlobConvALM}}]
            Suppose that $\ptaccum \in \mani$ is a local minimum of \RICO{}~\cref{prob:RICO}.
            Then, there exist sequences $\brc*{\ptotriter}_{\otriteridx}\subseteq\mani$ and $\brc*{\ineqLagmultotriter[]}_{\otriteridx}\subseteq\setRp[\ineqdime]$ such that
            \isextendedversion{
            \begin{equation}\label{def:AKKTconditions}
                \lim_{\otriteridx\to\infty} \ptotriter = \ptaccum, \, \lim_{\otriteridx\to\infty} \gradstr[\pt]\Lagfun\paren*{\allvarotriter} = 0, \, \lim_{\otriteridx\to\infty} \sum_{\ineqidx=1}^{\ineqdime} \max\brc*{ \ineqLagmultotriter[\ineqidx]\ineqfun[\ineqidx]\paren*{\ptotriter}, 0} = 0,
            \end{equation}
            where $\allvarotriter\coloneqq\paren*{\ptotriter, \ineqLagmultotriter[]}$.
            We call \cref{def:AKKTconditions}}{$\lim_{\otriteridx\to\infty} \ptotriter = \ptaccum$, $\lim_{\otriteridx\to\infty} \gradstr[\pt]\Lagfun\paren*{\allvarotriter}$ $= 0$, $\lim_{\otriteridx\to\infty} \sum_{\ineqidx=1}^{\ineqdime}\hspace{-1mm}\max\brc*{ \ineqLagmultotriter[\ineqidx]\ineqfun[\ineqidx]\paren*{\ptotriter}, 0}$ $= 0$,
            where $\allvarotriter\coloneqq\paren*{\ptotriter, \ineqLagmultotriter[]}$.
            We call them}
            the \AKKT{} conditions of \RICO{}~\cref{prob:RICO} and $\ptaccum$ satisfying the \AKKT{} conditions an \AKKT{} point of \RICO{}~\cref{prob:RICO}.
        \end{theorem}

        \begin{proposition}[{\cite[Theorem~2]{YamakawaSato2022SeqOptimCondforNLOonRiemManiandGlobConvALM}}]\label{prop:eqivKKTAKKT}
            Suppose that $\ptaccum\in\mani$ is an arbitrary point satisfying the \LICQ.
            Then, the following two statements are equivalent:
            \begin{enumerate}
                \item $\ptaccum$ is a \KKT{} point of \RICO{}~\cref{prob:RICO}.
                \item $\ptaccum$ is an \AKKT{} point of \RICO{}~\cref{prob:RICO}.
            \end{enumerate}
        \end{proposition}

        \begin{definition}[{\cite[Definition 2.5]{ObaraOkunoTakeda2021SQOforNLOonRiemMani}}]
            Given $\ptaccum \in \feasirgn$ satisfying the \KKT{} conditions with an associated Lagrange multiplier vector $\ineqLagmultaccum[] \in \setRp[\ineqdime]$, we say that the strict complementarity condition (\SC{}) holds if exactly one of $\ineqLagmultaccum$ and $\ineqfun[\ineqidx]\paren*{\ptaccum}$ is zero for each index $\ineqidx \in \ineqset$. 
            Hence, we have $\ineqLagmultaccum > 0$ for every $\ineqidx \in \activeineqset\paren*{\ptaccum}$.
        \end{definition}
        We define the critical cone associated with $\allvaraccum\coloneqq\paren*{\ptaccum, \ineqLagmultaccum[]} \in \mani \times \setRp[\ineqdime]$ as
        \begin{equation}\label{def:criticalcone}
            \criticalcone\paren*{\allvaraccum} \coloneqq \left\{ \tanvecone[\ptaccum] \in \tanspc[\ptaccum]\mani \middle|
            \begin{aligned}
                &\metr[\ptaccum]{\gradstr\ineqfun[\ineqidx]\paren*{\ptaccum}}{\tanvecone[\ptaccum]} = 0 \text{ for all } \ineqidx \in \activeineqset\paren*{\ptaccum} \text{ with } \ineqLagmultaccum > 0, \\
                &\metr[\ptaccum]{\gradstr\ineqfun[\ineqidx]\paren*{\ptaccum}}{\tanvecone[\ptaccum]} \geq 0 \text{ for all } \ineqidx \in \activeineqset\paren*{\ptaccum} \text{ with } \ineqLagmultaccum = 0 \\
            \end{aligned}
            \right\}.
        \end{equation}
        
        \begin{theorem}[{\cite[Theorem~4.2]{YangZhangSong14OptimCondforNLPonRiemMani}}]
            Suppose that $\ptaccum \in \mani$ is a local minimum of \RICO{}~\cref{prob:RICO} and the \LICQ{} holds at $\ptaccum$.
            Let $\ineqLagmultaccum[] \in \setRp[\ineqdime]$ be the vector of Lagrange multipliers for the inequality constraints.
            Then,
            \isextendedversion{
            \begin{equation}\label{eq:secondordernecessarycond}
                \metr[\ptaccum]{\Hess[\pt]\Lagfun\paren*{\allvaraccum}\sbra*{\tanvecone[\ptaccum]}}{\tanvecone[\ptaccum]} \geq 0 \text{ for all } \tanvecone[\ptaccum] \in \criticalcone\paren*{\allvaraccum}.
            \end{equation}
            We call \cref{eq:secondordernecessarycond}
            }{
            $\metr[\ptaccum]{\Hess[\pt]\Lagfun\paren*{\allvaraccum}\sbra*{\tanvecone[\ptaccum]}}{\tanvecone[\ptaccum]} \geq 0$ for all $\tanvecone[\ptaccum] \in \criticalcone\paren*{\allvaraccum}$.
            We call it}
            the second-order necessary condition and such $\ptaccum$ an \SOSP{}.
        \end{theorem}
        We next define the weak second-order necessary condition as 
        \isextendedversion{
        \begin{align}\label{eq:weaksecondordernecessarycond}
            \metr[\ptaccum]{\Hess[\pt]\Lagfun\paren*{\allvaraccum}\sbra*{\tanvecone[\ptaccum]}}{\tanvecone[\ptaccum]} \geq 0 \text{ for all } \tanvecone[\ptaccum] \in \weakcriticalcone\paren*{\ptaccum},
        \end{align}
        }{
        $\metr[\ptaccum]{\Hess[\pt]\Lagfun\paren*{\allvaraccum}\sbra*{\tanvecone[\ptaccum]}}{\tanvecone[\ptaccum]}$ $\geq 0$ for all $\tanvecone[\ptaccum] \in \weakcriticalcone\paren*{\ptaccum}$,
        }
        where 
        \isextendedversion{
        \begin{equation}\label{def:weakcriticalcone}
            \weakcriticalcone\paren*{\ptaccum} \coloneqq \left\{ \tanvecone[\ptaccum] \in \tanspc[\ptaccum]\mani \relmiddle{|}
            \metr[\ptaccum]{\gradstr\ineqfun[\ineqidx]\paren*{\ptaccum}}{\tanvecone[\ptaccum]} = 0 \text{ for all } \ineqidx \in \activeineqset\paren*{\ptaccum}
            \right\}.
        \end{equation}}{\begin{equation}
            \weakcriticalcone\paren*{\ptaccum} \coloneqq \left\{ \tanvecone[\ptaccum] \in \tanspc[\ptaccum]\mani \relmiddle{|}
            \metr[\ptaccum]{\gradstr\ineqfun[\ineqidx]\paren*{\ptaccum}}{\tanvecone[\ptaccum]} = 0 \text{ for all } \ineqidx \in \activeineqset\paren*{\ptaccum}
            \right\}.
        \end{equation}}
        We call such a point $\ptaccum\in\feasirgn$ a \wSOSP{}~\cite[Definition~3.3]{LiuBoumal2020SimpleAlgoforOptimonRiemManiwithCstr}.
        By definition, an \SOSP{} is also a \wSOSP{}.
        The converse is not necessarily true.
        Under the strict complementarity condition, however, these two points are identical since $\criticalcone = \weakcriticalcone$ holds.
        \begin{remark}\label{prop:wSOSPequivSOSPudrSC}
            Let $\ptaccum \in \mani$ be a \wSOSP{} of \RICO{}~\cref{prob:RICO}.
            Suppose that the \SC{} holds at $\ptaccum$.
            Then, $\ptaccum$ is an \SOSP{} of \RICO{}~\cref{prob:RICO}.
        \end{remark}
 

    \section{Proposed method}\label{sec:propmeth}
        In this section, we propose \RIPTRM{} for solving \RICO{} \cref{prob:RICO}.
        First, we introduce the \KKT{} vector field, which was originally proposed by \cite[Section~3.1]{LaiYoshise2024RiemIntPtMethforCstrOptimonMani}.
        The vector field will be the basis of the subproblems to be solved in our algorithm.
        Then, we propose \RIPTRM{}, which consists of outer and inner iterations.

    \subsection{Overview of \RIPTRM{}}\label{subsec:overviewRIPTRM}
        Let $\barrparam[] > 0$ be a barrier parameter.
        We define the barrier \KKT{} vector field as 
        \begin{equation}\label{eq:barrKKTvecfld}
            \barrKKTvecfld\paren*{\allvar[];\barrparam[]}\coloneqq 
            \begin{bmatrix}
                \gradstr[\pt]\Lagfun\paren*{\allvar[]}\\
                \IneqLagmultmat[]\ineqfun\paren*{\pt} - \barrparam[]\onevec
            \end{bmatrix}.
        \end{equation}
        This vector field originates from the \isextendedversion{\KKT{} conditions~\cref{eq:KKTconditions}}{\KKT{} conditions}; it consists of the Riemannian gradient of the Lagrangian and the complementarity condition relaxed by the barrier parameter.
        The condition $\barrKKTvecfld\paren*{\allvar[];0} = 0$ is equivalent to the \isextendedversion{\KKT{} conditions~\cref{eq:KKTconditions}}{\KKT{} conditions} for $\allvar[] \in \feasirgn \times \setRp[\ineqdime]$.
        For brevity, we often write $\barrKKTvecfld\paren*{\allvar[]}$ for $\barrKKTvecfld\paren*{\allvar[];0}$.
        
        \RIPTRM{} aims to find the solution of $\barrKKTvecfld\paren*{\allvar[]}=0$ with $\allvar[]\in\feasirgn\times\setRp[\ineqdime]$ by approximately solving
        \begin{equation}\label{eq:barrKKTequation}
            \barrKKTvecfld\paren*{\allvar[];\barrparamotriter}=0, \quad \allvar[]\in\strictfeasirgn\times\setRpp[\ineqdime]
        \end{equation}
        for a sequence of barrier parameters $\brc*{\barrparamotriter}_{\otriteridx}\subset \setRpp[]$ that converges to zero.
        Hereafter, for $\allvar[] = \paren*{\pt[], \ineqLagmult[]}$, we refer to $\pt[]$ as the primal variable and $\ineqLagmult[]$ as the dual variable, respectively.
        Let $\allvarotriter\coloneqq\paren*{\ptotriter, \ineqLagmultotriter[]} \in \strictfeasirgn \times \setRpp[\ineqdime]$ denote an approximate solution of \cref{eq:barrKKTequation}.
        We call the sequence $\brc*{\allvarotriter}_{\otriteridx}$ the outer iterates. 
        The associated adjustment of the barrier parameter and tolerances for residuals defines the outer iteration, whose index is the subscript $\otriteridx\in\setNz$.
        Each $\allvarotriter$ is an output of the inner iterations, where a corresponding sequence of inner iterates is indexed by the superscript $\inriteridx\in\setNz$.
        At each inner iteration, the algorithm approximately solves \cref{eq:barrKKTequation} with $\barrparamotriter > 0$ fixed while keeping $\allvarinriter\coloneqq\paren*{\ptinriter, \ineqLagmultinriter[]}\in\strictfeasirgn\times\setRpp[\ineqdime]$.
        In the following, we explain the details of the outer and inner iterations in turn.

    \subsection{Outer iteration}\label{subsec:outeriter}
        In this subsection, we present the specific outer process of the proposed algorithm.
        Let $\forcingfungradLag, \forcingfuncompl, \forcingfunsosp \colon \setRp[]\to\setRp[]$ be forcing functions.
        At the $\otriteridx$-th outer iteration, 
        we choose a barrier parameter $\barrparamotriter>0$ such that $\lim_{\otriteridx\to\infty}\barrparamotriter=0$.
        Here, the update can depend on computations performed by the algorithm up to outer iteration $\otriteridx$.
        We find a point $\allvarotriterp\in\strictfeasirgn\times\setRpp[\ineqdime]$ that satisfies the following stopping conditions with $\barrparam > 0$:
        \begin{subequations}\label{eq:stopcond}
            \begin{align}
                    &\Riemnorm[{\ptotriterp}]{\gradstr[{\pt[]}]\Lagfun\paren*{\allvarotriterp}} \leq \forcingfungradLag\paren*{\barrparam},\label{eq:stopcondKKT}\\
                    &\norm*{\IneqLagmultmatotriterp \ineqfun[]\paren*{\ptotriterp} - \barrparam \onevec} \leq \forcingfuncompl\paren*{\barrparam},\label{eq:stopcondbarrcompl}\\
                    &\mineigval\sbra*{\trsquadotriterp} \geq - \forcingfunsosp\paren*{\barrparam},\label{eq:stopcondsecondord} \\
                    & \ineqfun[]\paren*{\ptotriterp} > 0, \text{ and } \ineqLagmultotriterp[] > 0,\label{eq:stopcondstrictfeasi}
            \end{align}
        \end{subequations}
        where we write $\IneqLagmultmatotriterp$ and $\trsquadotriterp$ for $\diag\paren*{\ineqLagmultotriterp[]}$ and $\trsquad\paren*{\allvarotriterp}$, \isextendedversion{
        as defined later in \cref{eq:trsquaddef}, 
        }{which will be defined in the next subsection,}
        respectively.
        If our goal is to compute not an \SOSP{} but an \AKKT{} point, we can remove \cref{eq:stopcondsecondord} from the stopping conditions.
        To gain such $\allvarotriterp$, for example, we use \cref{algo:RIPTRMInner} with an initial point $\allvarotriter\in\strictfeasirgn\times\setRpp[\ineqdime]$, the $\otriteridx$-th barrier parameter $\barrparamotriter > 0$, and the $\otriteridx$-th initial trust region radius $\inittrradius_{\otriteridx} \in \left(0,\maxtrradius\right]$.
        We also store the final trust region radius $\finaltrradius_{\otriteridx}$ used in the inner iteration.
        Then, the algorithm defines $\inittrradius_{\otriteridx+1} \coloneqq \max\paren*{\finaltrradius_{\otriteridx}, \mininittrradius}$, where $\mininittrradius > 0$ is a predefined parameter.
        The outer iteration is formally presented as \cref{algo:RIPTRMOuter}.        

        \begin{algorithm}[t]
            \SetKwInOut{Require}{Require}
            \SetKwInOut{Input}{Input}
            \SetKwInOut{Output}{Output}
            \Require{Riemannian manifold $\mani$, twice continuously differentiable functions $\objfun:\mani\to\setR[]$ and $\brc*{\ineqfun[\ineqidx]}_{\ineqidx\in\ineqset}\colon\mani\to\setR[]$, maximal trust region radius $\maxtrradius > 0$, initial trust region radius $\inittrradius_{0} \in \left(0,\maxtrradius\right]$, 
            minimum initial trust region radius $\mininittrradius\in \left(0,\maxtrradius\right]$, 
            forcing functions $\forcingfungradLag,  \forcingfuncompl, \forcingfunsosp\colon\setRp[]\to\setRp[]$, ${\rm \seccondflag\in\brc*{\True, \False}}$.}
            \Input{Initial point $\allvar[0] = \paren*{\pt[0], \ineqLagmult[0]} \in \strictfeasirgn \times \setRpp[\ineqdime]$,
            initial barrier parameter $\barrparam[-1] > 0$.}
            \For{$\otriteridx = 0, 1, \ldots$}{
                Set $\barrparamotriter > 0$ such that 
                $\lim_{\otriteridx \to \infty} \barrparamotriter = 0$.
                Set \stopcond{} $\gets$ \cref{eq:stopcond} if \seccondflag{} is \True{}. Else, set \stopcond{} $\gets$ \cref{eq:stopcondKKT,eq:stopcondbarrcompl,eq:stopcondstrictfeasi}.
                
                Compute $\allvarotriterp\in \strictfeasirgn \times \setRpp[\ineqdime]$ satisfying \stopcond{} and $\finaltrradius_{\otriteridx} > 0$ using, e.g., \cref{algo:RIPTRMInner}~$\paren*{\allvarotriter, \barrparamotriter, \inittrradius_{\otriteridx}, {\rm \stopcond{}}}$.

                Set $\inittrradius_{\otriteridx+1}\gets\max\paren*{\finaltrradius_{\otriteridx}, \mininittrradius}$.
            }
            \caption{Outer iteration of \RIPTRM{} 
            }\label{algo:RIPTRMOuter}
        \end{algorithm}
        
    \subsection{Inner iteration}\label{subsec:inneriter}
        To approximately solve \cref{eq:barrKKTequation} with $\barrparam[] = \barrparamotriter > 0$ fixed, 
        we consider the Newton system of \cref{eq:barrKKTequation}:
        for 
        $\allvar=\paren*{\pt,\ineqLagmult[]}\in\strictfeasirgn\times\setR[\ineqdime]$,
        the Newton system is given by
        \isextendedversion{
        \begin{align}\label{eq:RiemNewtoneq}
            \Jacobian[\barrKKTvecfld]\paren*{\allvar[]}\sbra*{\dirallvar} = - \barrKKTvecfld\paren*{\allvar[];\barrparam[]},
        \end{align}
        where
        \begin{align}
            \begin{split}\label{def:JacobbarrKKTvecfld}
                \Jacobian[\barrKKTvecfld]\paren*{\allvar[]}\colon\tanspc[\pt]\mani\times\setR[\ineqdime]&\to\tanspc[\pt]\mani\times\setR[\ineqdime]\\
                \dirallvar\coloneqq\paren*{\dirpt, \dirineqLagmult[]} &\mapsto 
                \begin{bmatrix}
                    \Hess[\pt]\Lagfun\paren*{\allvar}\sbra*{\dirpt} - \ineqgradopr[\pt]\sbra*{\dirineqLagmult[]}\\
                    \IneqLagmultmat[]\coineqgradopr[\pt]\sbra*{\dirpt} + \Ineqfunmat[]\paren*{\pt}\dirineqLagmult[]
                \end{bmatrix}
            \end{split}
        \end{align}}{
        $\Jacobian[\barrKKTvecfld]\paren*{\allvar[]}\sbra*{\dirallvar} = - \barrKKTvecfld\paren*{\allvar[];\barrparam[]}$,
        where
        \begin{align}
            \begin{split}\label{def:JacobbarrKKTvecfld}
                \Jacobian[\barrKKTvecfld]\paren*{\allvar[]}\colon\tanspc[\pt]\mani\times\setR[\ineqdime]&\to\tanspc[\pt]\mani\times\setR[\ineqdime]\\
                \dirallvar\coloneqq\paren*{\dirpt, \dirineqLagmult[]} &\mapsto 
                \begin{bmatrix}
                    \Hess[\pt]\Lagfun\paren*{\allvar}\sbra*{\dirpt} - \ineqgradopr[\pt]\sbra*{\dirineqLagmult[]}\\
                    \IneqLagmultmat[]\coineqgradopr[\pt]\sbra*{\dirpt} + \Ineqfunmat[]\paren*{\pt}\dirineqLagmult[]
                \end{bmatrix}
            \end{split}
        \end{align}
        }
        is the Jacobian of $\barrKKTvecfld\paren*{\cdot;\barrparam[]}$ at $\allvar$ under $\tanspc[{\ineqLagmult[]}]\setR[\ineqdime]\cid\setR[\ineqdime]$.
        Note that the Jacobian is independent of the value of the barrier parameter $\barrparam[]$.
        Since $\Ineqfunmat[]\paren*{\pt}\in\setR[\ineqdime\times\ineqdime]$ is nonsingular due to $\pt\in\strictfeasirgn$, rearranging the Newton system with respect to $\dirineqLagmult[]$ gives the following equation:
        \begin{align}
            &\trsquad\paren*{\allvar}\sbra*{\dirpt} = - \trslin[{\barrparam[]}]\paren*{\pt},\label{eq:primalNewtoneq}\\
            &\dirineqLagmult[] = - \ineqLagmult[] + \barrparam[]\inv{\Ineqfunmat[]\paren*{\pt}}\onevec - \IneqLagmultmat[]\inv{\Ineqfunmat[]\paren*{\pt}}\coineqgradopr[\pt]\sbra*{\dirpt},\label{eq:dualNewtoneq}
        \end{align}
        where
        \isextendedversion{
        \begin{align}
            &\trsquad\paren*{\allvar} \coloneqq \Hess[\pt]\Lagfun\paren*{\allvar} + \ineqgradopr[\pt]\IneqLagmultmat[]\inv{\Ineqfunmat[]\paren*{\pt}}\coineqgradopr[\pt],\label{eq:trsquaddef}\\
            &\trslin[{\barrparam[]}]\paren*{\pt}\coloneqq\gradstr[]\objfun\paren*{\pt} - \barrparam[]\ineqgradopr[\pt]\sbra*{\inv{\Ineqfunmat[]\paren*{\pt}}\onevec}.\label{eq:trslindef}
        \end{align}}{
        $\trsquad\paren*{\allvar} \coloneqq \Hess[\pt]\Lagfun\paren*{\allvar} + \ineqgradopr[\pt]\IneqLagmultmat[]\inv{\Ineqfunmat[]\paren*{\pt}}\coineqgradopr[\pt] \text{ and } \trslin[{\barrparam[]}]\paren*{\pt}\coloneqq\gradstr[]\objfun\paren*{\pt} - \barrparam[]\ineqgradopr[\pt]\sbra*{\inv{\Ineqfunmat[]\paren*{\pt}}\onevec}$.
        }
         
        Based on \cref{eq:primalNewtoneq}, for each $\inriteridx$-th iteration, \RIPTRM{} computes a search direction for the primal variable using the trust region strategy; we introduce the trust region subproblem at $\pt\in\strictfeasirgn$ with the dual variable $\ineqLagmult[]\in\setRpp[\ineqdime],$ the trust region radius $\trradius[] > 0$, and the barrier parameter $\barrparam[] > 0$ as follows:
        \isextendedversion{
        \begin{mini}
            {\dir[] \in \tanspc[\pt]\mani}{\modelfun[{\allvar, \barrparam[]}]\paren*{\dir[]} \coloneqq \frac{1}{2}\metr[\pt]{\trsquad\paren*{\allvar}\sbra*{\dir[]}}{\dir[]} + \metr[\pt]{\trslin[{\barrparam[]}]\paren*{\pt}}{\dir[]}}
            {\label{prob:TRS}}{}
            \addConstraint{\Riemnorm[\pt]{\dir[]}}{\leq \trradius[].\label{cstr:trradius}}
        \end{mini}}{\begin{equation}\label{prob:TRS}
            \underset{\dir[] \in \tanspc[\pt]\mani}{\minimize}~\modelfun[{\allvar, \barrparam[]}]\paren*{\dir[]} \coloneqq \frac{1}{2}\metr[\pt]{\trsquad\paren*{\allvar}\sbra*{\dir[]}}{\dir[]} + \metr[\pt]{\trslin[{\barrparam[]}]\paren*{\pt}}{\dir[]}~\subjectto~\Riemnorm[\pt]{\dir[]} \leq \trradius[].
        \end{equation}}
        Despite the fact that the subproblem is generally a nonconvex quadratic optimization problem,
        it is possible to solve the subproblem exactly; see, e.g., \cite{MoracuteeSorensen1983ComputTrustRegionStep,Adachietal2017SolvingTRSbyGenEigenProb,BeckVaisbourd2018GlobSolvingTRSUsingSimpleFirstOrdMeths}. 
        \isextendedversion{
        }{We provide further discussion in \cref{subsec:experimentenv}.}
        \RIPTRM{} computes the search direction for the primal variable, denoted by $\dirptinriter \in \tanspc[\ptinriter]\mani$, by approximately or exactly solving the trust region subproblem \cref{prob:TRS} at $\ptinriter\in\strictfeasirgn$ with $\ineqLagmultinriter[]\in\setRpp[\ineqdime]$ and $\trradiusinriter, \barrparam[] \in \setRpp[]$.
        
        By using $\dirptinriter\in\tanspc[\ptinriter]\mani$ that satisfies $\Riemnorm[\ptinriter]{\dirptinriter}\leq \trradiusinriter$, \RIPTRM{} computes the search direction for the dual variable $\ineqLagmultinriter[]\in\setR[\ineqdime]$, denoted by $\dirineqLagmultinriter[]\in\tanspc[{\ineqLagmultinriter[]}]\setR[\ineqdime]\cid\setR[\ineqdime]$, according to \cref{eq:dualNewtoneq}. 
        If the point $\paren*{\retr[\ptinriter]\paren*{\dirptinriter}, \ineqLagmultinriter[] + \dirineqLagmultinriter[]}$ satisfies the stopping conditions \cref{eq:stopcond}
        with $\barrparam[] > 0$, the algorithm outputs this point.
        Otherwise, several tests are performed to determine whether to update the trust region radius and the variables.
        First, \RIPTRM{} tests whether $\retr[\ptinriter]\paren*{\dirptinriter}$ belongs to $\strictfeasirgn$.
        If not, it shrinks the trust region radius $\trradiusinriterp \gets \contcoeff\Riemnorm[\ptinriter]{\dirptinriter}$ with $0<\contcoeff<1$ and starts over from the beginning of the next iteration, computing another search direction due to the smaller trust region radius.
        Otherwise, the algorithm proceeds to the other tests: 
        for these tests, we introduce the log barrier function 
        \isextendedversion{
        \begin{equation}\label{eq:logbarrfun}
            \meritfun[{\barrparam[]}]\paren*{\pt[]} \coloneqq \objfun\paren*{\pt[]} - \barrparam[] \sum_{\ineqidx\in\ineqset} \log\ineqfun[\ineqidx]\paren*{\pt}
        \end{equation}}{
        \begin{align}\label{eq:logbarrfun}
            \meritfun[{\barrparam[]}]\paren*{\pt[]} \coloneqq \objfun\paren*{\pt[]} - \barrparam[] \sum_{\ineqidx\in\ineqset} \log\ineqfun[\ineqidx]\paren*{\pt}
        \end{align}
        }
        for $\pt\in\strictfeasirgn$, which serves as a merit function.
        We also define two reductions as \isextendedversion{
        \begin{subequations}
            \begin{align}
                &\ared[{\barrparam[]}]\paren*{\dirpt} \coloneqq \pullback[\pt]{\meritfun[{\barrparam[]}]}\paren*{\zerovec[\pt]} - \pullback[\pt]{\meritfun[{\barrparam[]}]}\paren*{\dirpt},\label{def:aredinriter}\\
                &\pred[{\allvar, \barrparam[]}]\paren*{\dirpt} \coloneqq \modelfun[{\allvar, \barrparam[]}]\paren*{\zerovec[\pt]} - \modelfun[{\allvar, \barrparam[]}]\paren*{\dirpt}, \label{def:predinriter}
            \end{align} 
        \end{subequations}
        }{
        \begin{align}
            \ared[{\barrparam[]}]\paren*{\dirpt} \coloneqq \pullback[\pt]{\meritfun[{\barrparam[]}]}\paren*{\zerovec[\pt]} - \pullback[\pt]{\meritfun[{\barrparam[]}]}\paren*{\dirpt} \text{ and } \pred[{\allvar, \barrparam[]}]\paren*{\dirpt} \coloneqq \modelfun[{\allvar, \barrparam[]}]\paren*{\zerovec[\pt]} - \modelfun[{\allvar, \barrparam[]}]\paren*{\dirpt},
        \end{align}
        }
        which we call the actual and predicted reductions, respectively.
        Here, we recall that $\pullback[\pt]{\meritfun[{\barrparam[]}]} = \meritfun[{\barrparam[]}] \circ \retr[\pt]$.
        The actual reduction is the change in \isextendedversion{\cref{eq:logbarrfun}}{the merit function}
        produced by the step, and the predicted reduction is that of the model
        in \cref{prob:TRS}.
        \RIPTRM{} sets the trust region radius at the next iteration according to the ratio of the two reductions:
        \begin{equation}\label{eq:trradiusupdate}
            \trradiusinriterp \gets
            \begin{cases}
                \frac{1}{4} \trradiusinriter & \aredinriter < \frac{1}{4}\predinriter, \\
                \min\paren*{2\trradiusinriter, \maxtrradius} & \aredinriter \geq \frac{3}{4} \predinriter \text{ and } \Riemnorm[\ptinriter]{\dirptinriter} = \trradiusinriter,\\
                \trradiusinriter & \text{otherwise},
            \end{cases}
        \end{equation}
        where $\maxtrradius > 0$ is a predefined parameter called the maximal trust region radius, and we write $\aredinriter$ and $\predinriter$ for $\ared[{\barrparam[]}]\paren*{\dirptinriter}$ and $\pred[{\allvarinriter, \barrparam[]}]\paren*{\dirptinriter}$, respectively.
        Here, we omit the subscript $\barrparam[]$ for brevity.
        Thereafter, the algorithm updates the iterate as
        \begin{equation}\label{eq:primalupdate}
            \ptinriterp\gets
            \begin{cases}
                \retr[\ptinriter]\paren*{\dirptinriter} & \aredinriter > \trratiothold \predinriter, \\
                \ptinriter & \text{otherwise,}
            \end{cases}
        \end{equation}
        where $\trratiothold\in\paren*{0,\frac{1}{4}}$ is a predefined parameter working as the threshold of the ratio.
        We say that the $\inriteridx$-th iteration is {\it successful} if $\retr[\ptinriter]\paren*{\dirptinriter}\in\strictfeasirgn$ and $\aredinriter > \trratiothold \predinriter$ 
        hold and the iterate is redefined.
        Otherwise, we say that the iteration is {\it unsuccessful}.
        \RIPTRM{} also updates $\ineqLagmultinriterp[] \in \setRpp[\ineqdime]$ using $\ineqLagmultinriter[]$ and $\dirineqLagmultinriter[]$ when the iteration is successful.
        From a theoretical perspective, any update rule for the dual variable $\ineqLagmult[]$ is admissible provided that it satisfies \cref{assu:ineqLagmultbounded,assu:ineqLagmultconverged}, which will be stated in \cref{subsec:globconvinner}.        \isextendedversion{
        In this paper, we introduce clipping for this update originally introduced in \cite[Section~4.3]{Connetal2000PrimalDualTRAlgoforNonconvexNLP}; see \cref{subsec:computdualvar} for the details.}{
        In this paper, we employ the clipping update for the dual variables, which was originally proposed in the Euclidean setting \cite[Section~4.3]{Connetal2000PrimalDualTRAlgoforNonconvexNLP}.
        Let $\constfou, \constfiv \in \setR$ be constants satisfying $0 < \constfou < 1 < \constfiv$.
        Then, define $\clipintvlinritermin[\ineqidx] \coloneqq \constfou\min\brc*{1, \ineqLagmultinriter[\ineqidx], \frac{\barrparam[]}{\pullback[\ptinriter]{\ineqfun[\ineqidx]}\paren*{\dirptinriter}}}$ and $\clipintvlinritermax[\ineqidx] \coloneqq \max\brc*{\constfiv, \ineqLagmultinriter[\ineqidx], \frac{\constfiv}{\barrparam[]}, \frac{\constfiv}{\pullback[\ptinriter]{\ineqfun[\ineqidx]}\paren*{\dirptinriter}}}$ for each $\ineqidx\in\ineqset$. 
        We update the dual variable as follows:
        \begin{align}\label{def:clippingupdate}
            \ineqLagmultinriterp[] \leftarrow 
            \begin{cases}
                \clip[{\clipintvl[\inriteridx]}]\paren*{\ineqLagmultinriter[] + \dirineqLagmultinriter[]} & \text{if } \ptinriterp=\retr[\ptinriter]\paren*{\dirptinriter}, \\
               \ineqLagmultinriter[] & \text{if } \ptinriterp=\ptinriter,\\
            \end{cases}
        \end{align}
        where $\clip[{\clipintvl[\inriteridx]}]\colon\setR[\ineqdime]\to\setR[\ineqdime]$ is the operator whose $\ineqidx$-th component is defined by $\sbra*{\clip[{\clipintvl[\inriteridx]}]\paren*{\vecfou}}_{\ineqidx}$ $\coloneqq \max\brc*{\clipintvlinritermin[\ineqidx], \min\brc*{\clipintvlinritermax[\ineqidx], \vecfou_{\ineqidx}}}$.}
        \cref{algo:RIPTRMInner} formally states the inner iteration.

        \begin{algorithm}[t]
            \SetKwInOut{Require}{Require}
            \SetKwInOut{Input}{Input}
            \SetKwInOut{Output}{Output}
            \Require{Riemannian manifold $\mani$, twice continuously differentiable functions $\objfun:\mani\to\setR[]$ and $\brc*{\ineqfun[\ineqidx]}_{\ineqidx\in\ineqset}\colon\mani\to\setR[]$, maximal trust region radius $\maxtrradius > 0$, threshold $\trratiothold \in \paren*{0, \frac{1}{4}}$, coefficient $\contcoeff \in \paren*{0,1}$.}
            \Input{Initial point ${\allvar[]}^{0} = \paren*{{\pt[]}^{0}, {\ineqLagmult[]}^{0}} \in \strictfeasirgn \times \setRpp[\ineqdime]$,
            barrier parameter $\barrparam[] > 0$,
            initial trust region radius $\trradius[]^{0} \in \left(0,\maxtrradius\right]$, stopping conditions \stopcond{}.}
            \Output{Final point $\allvar[]=\paren*{\pt[],\ineqLagmult[]} \in \strictfeasirgn \times \setRpp[\ineqdime]$ and $\finaltrradius > 0$.}
            \For{$\inriteridx = 0, 1, \ldots$\label{algo:forline}}{
                Compute $\dirptinriter \in \tanspc[\ptinriter]\mani$ by approximately or exactly solving the subproblem~\cref{prob:TRS}.\;
                
                Compute $\dirineqLagmultinriter\in\tanspc[{\ineqLagmultinriter[]}]\setR[\ineqdime]$ according to \cref{eq:dualNewtoneq}.\;
                
                \If{$\paren*{\retr[\ptinriter]\paren*{\dirptinriter}, \ineqLagmultinriter[] + \dirineqLagmultinriter[]}$ satisfies {\rm \stopcond{}}}
                {\Return $\allvar[]=\paren*{\retr[\ptinriter]\paren*{\dirptinriter}, \ineqLagmultinriter[] + \dirineqLagmultinriter[]}$ and $\finaltrradius=\trradiusinriter$.}                    
                \If{$\retr[\ptinriter]\paren*{\dirptinriter}\notin\strictfeasirgn$}
                {Set $\ptinriterp\leftarrow\ptinriter$, $\ineqLagmultinriterp[]\leftarrow\ineqLagmultinriter[]$, and $\trradiusinriterp \gets \contcoeff\Riemnorm[\ptinriter]{\dirptinriter}$ and go to line~\ref{algo:forline}.}
                \isextendedversion{
                Compute $\aredinriter$ and $\predinriter$ by \cref{def:aredinriter,def:predinriter}, respectively.\; 
                }{
                Compute $\aredinriter$ and $\predinriter$, respectively.\;}        

                Update $\trradiusinriterp$ according to \cref{eq:trradiusupdate}.\;

            \isextendedversion{
            Define $\ptinriterp$ by \cref{eq:primalupdate}
            and update $\ineqLagmultinriterp[] \in \setRpp[\ineqdime]$ by $\ineqLagmultinriter[]$ and $\dirineqLagmultinriter[]$; see \cref{def:clippingupdate} in \cref{subsec:computdualvar}, for example.
            }{
            Define $\ptinriterp$ by \cref{eq:primalupdate}
            and update $\ineqLagmultinriterp[] \in \setRpp[\ineqdime]$ by $\ineqLagmultinriter[]$ and $\dirineqLagmultinriter[]$ using, e.g., clipping~\cref{def:clippingupdate}.}
            }
            \caption{Inner iteration of \RIPTRM{}
            }\label{algo:RIPTRMInner}
        \end{algorithm}


    \section{Global convergence analysis}\label{sec:globconv}
        In this section, we establish the global convergence properties of \RIPTRM{}.
        We first analyze \cref{algo:RIPTRMInner}; we prove its 
        validity
        and the global convergence in \cref{subsec:wellposedness,subsec:globconvinner}, respectively.
        Then, we prove the global convergence of \cref{algo:RIPTRMOuter} in \cref{subsec:globconvouter}.
        
        \stepcounter{assumptionsection}
        
    \subsection{
    Validity
    of 
    \texorpdfstring{\cref{algo:RIPTRMInner}}{inner iterations}}\label{subsec:wellposedness}
        In the subsection, we prove 
        that \cref{algo:RIPTRMInner} is 
        valid
        in the sense that, if the current point is not a solution of $\barrKKTvecfld\paren*{\cdot;\barrparam[]} = 0$ and the trust region radius is sufficiently small,
        the iteration becomes successful; that is, $\retr[\ptinriter]\paren*{\dirptinriter}\in\strictfeasirgn$ and $\aredinriter > \trratiothold \predinriter$ hold.
        Hereafter, for brevity, we write $\modelfuninriter$, $\trsquadinriter$, and $\trslininriter$ for $\modelfun[{\allvarinriter, \barrparam[]}]$, $\trsquad\paren*{\allvarinriter}$, and $\trslin[{\barrparam[]}]\paren*{\ptinriter}$ by omitting the subscript $\barrparam[]$, respectively.
        \isextendedversion{}{See \cref{subsec:inneriter} for the definitions of $\ared$, $\pred$, $\modelfun[{\allvar, \barrparam[]}]$, $\trsquad$, and $\trslin[{\barrparam[]}]$.}
        For the 
        validity,
        we assume the following:        \begin{assumption}\label{assu:Cauchydecreasing}
            There exists $\constCauchy > 0$ such that, for all $\inriteridx\in\setNz$, the search direction $\dirptinriter\in\tanspc[\ptinriter]\mani$ satisfies
            \isextendedversion{
            \begin{align}
                \modelfuninriter\paren*{\zerovec[\ptinriter]} - \modelfuninriter\paren*{\dirptinriter} \geq \constCauchy\Riemnorm[\ptinriter]{\trslininriter}\min\paren*{\trradiusinriter, \frac{\Riemnorm[\ptinriter]{\trslininriter}}{\opnorm{\trsquadinriter}}}.
            \end{align}}{$\modelfuninriter\paren*{\zerovec[\ptinriter]} - \modelfuninriter\paren*{\dirptinriter} \geq \constCauchy\Riemnorm[\ptinriter]{\trslininriter}\min\paren*{\trradiusinriter, \frac{\Riemnorm[\ptinriter]{\trslininriter}}{\opnorm{\trsquadinriter}}}$.}
        \end{assumption}
        Note that, under \cref{assu:Cauchydecreasing}, $\dirptinriter \neq \zerovec[\ptinriter]$ holds if $\Riemnorm[\ptinriter]{\trslininriter} \neq 0$.
        \isextendedversion{From \cref{lemm:Cauchydecreasing}, the \Cauchypoint{}}{The \Cauchypoint{}~\cite[Lemma~6.15]{Boumal23IntroOptimSmthMani}} satisfies \cref{assu:Cauchydecreasing} with $\constCauchy=\frac{1}{2}$, and so does the \exactsolution{}, that is, the global optimum of \cref{prob:TRS}\isextendedversion{; see \cref{sec:searchdir} for the detail.}{.}
        \cref{assu:Cauchydecreasing} is standard in the literature; it is made in \cite[A6.3]{Boumal23IntroOptimSmthMani} and \cite[Section~7.4.1]{Absiletal08OptimAlgoonMatMani}, for example.
                
        We first derive a bound on the values of the inequality constraints around each feasible point: recall the definitions of the restriction of the pullback
        \isextendedversion{\cref{def:pullbackfun}. }{
        $\pullback[\pt]{\cdot}$ in \cref{subsec:notationRiemggeo}.}
        The lemma holds by the continuity of $\brc*{\ineqfun[\ineqidx]}_{\ineqidx\in\ineqset}$.
        \begin{lemma}\label{lemm:ineqfunfeasi}
            Choose $\pt\in\strictfeasirgn$ arbitrarily.
            Then, there exists $\deltaone{} > 0$ such that $\pullback[\pt]{\ineqfun[\ineqidx]}\paren*{\tanvecone[\pt]} \geq \frac{1}{2} \ineqfun[\ineqidx]\paren*{\pt} > 0$ holds for all $\ineqidx\in\ineqset$ and any $\tanvecone[\pt]\in\tanspc[\pt]\mani$ with $\Riemnorm[\pt]{\tanvecone[\pt]}\leq\deltaone{}$.
        \end{lemma}
        \isextendedversion{
        \begin{proof}
            See \cref{appx:proofconsistency}.
        \end{proof}
        }{
        \begin{proof}
            The statement follows from the continuity of $\retr[]$ and of $\{\ineqfun[\ineqidx]\}_{\ineqidx\in\ineqset}$.
            See \cite[Appendix~D.1]{Obaraetal2025APrimalDualIPTRMfor2ndOrdStnryPtofRiemIneqCstrOptimProbs} for a complete proof.
        \end{proof}
        }
        Using \cref{lemm:ineqfunfeasi}, we conclude that the search direction obtained from \cref{prob:TRS} passes the test in line~15 of \cref{algo:RIPTRMInner} when the trust region radius is sufficiently small.
        We formally state this at the end of this section.
        Next, we consider the tests on the ratio of the actual and predicted reductions.
        To this end, we derive the first and second directional derivatives of the log barrier function \cref{eq:logbarrfun} by direct calculation: \isextendedversion{}{we refer the reader to \cite[Appendix~D.1]{Obaraetal2025APrimalDualIPTRMfor2ndOrdStnryPtofRiemIneqCstrOptimProbs} for the proofs.}
        \begin{lemma}\label{lemm:firstdirderivmeritfun}
            Choose $\pt\in\strictfeasirgn$ and $\barrparam[] > 0$ arbitrarily.
            Then, 
            for all $\tanvecone[\pt]\in\tanspc[\pt]\mani$, 
            \isextendedversion{
            \begin{equation}
                \D\pullback[\pt]{\meritfun[{\barrparam[]}]}\paren*{\zerovec[\pt]}\sbra*{\tanvecone[\pt]} = \metr[\pt]{\trslin[{\barrparam[]}]\paren*{\pt}}{\tanvecone[\pt]}\label{eq:dirderivmeritfun}.
            \end{equation}
            }{$\D\pullback[\pt]{\meritfun[{\barrparam[]}]}\paren*{\zerovec[\pt]}\sbra*{\tanvecone[\pt]} = \metr[\pt]{\trslin[{\barrparam[]}]\paren*{\pt}}{\tanvecone[\pt]}$ holds.}
        \end{lemma}
        \isextendedversion{
        \begin{proof}
            See \cref{appx:proofconsistency}.
        \end{proof}
        }{}
        \begin{lemma}\label{lemm:twicedirderivmeritfun}
            Choose $\pt[]\in\strictfeasirgn$ and $\barrparam[] > 0$ arbitrarily.
            Let $\tanvecone[\pt] \in \tanspc[\pt]\mani$ be any tangent vector satisfying $\pullback[\pt]{\ineqfun[\ineqidx]}\paren*{\tanvecone[\pt]} \neq 0$ for all $\ineqidx \in \ineqset$.
            Then, for all $\tanvectwo[\pt], \tanvecthr[\pt] \in \tanspc[{\tanvecone[\pt]}]\paren*{\tanspc[\pt]\mani}\cid\tanspc[\pt]\mani$, 
            \isextendedversion{
            \begin{align}
                \begin{split}\label{eq:twicedirderivmeritfun}
                    &\D[2]\pullback[\pt]{\meritfun[{\barrparam[]}]}\paren*{\tanvecone[\pt]}\sbra*{\tanvectwo[\pt],\tanvecthr[\pt]} =\metr[\pt]{\paren*{\Hess\pullback[\pt]{\objfun}\paren*{\tanvecone[\pt]} -  \sum_{\ineqidx\in\ineqset}\frac{\barrparam[]}{\pullback[\pt]{\ineqfun[\ineqidx]}\paren*{\tanvecone[\pt]}}\Hess\pullback[\pt]{\ineqfun[\ineqidx]}\paren*{\tanvecone[\pt]}}\sbra*{\tanvectwo[\pt]}}{\tanvecthr[\pt]}\\
                    &\quad + \sum_{\ineqidx\in\ineqset}\frac{\barrparam[]}{{\pullback[\pt]{\ineqfun[\ineqidx]}\paren*{\tanvecone[\pt]}}^{2}}\metr[\pt]{\gradstr\pullback[\pt]{\ineqfun[\ineqidx]}\paren*{\tanvecone[\pt]}}{\tanvectwo[\pt]}\metr[\pt]{\gradstr\pullback[\pt]{\ineqfun[\ineqidx]}\paren*{\tanvecone[\pt]}}{\tanvecthr[\pt]}.
                \end{split}
            \end{align}  
            }{$\D[2]\pullback[\pt]{\meritfun[{\barrparam[]}]}\paren*{\tanvecone[\pt]}\sbra*{\tanvectwo[\pt],\tanvecthr[\pt]} =\metr[\pt]{\paren*{\Hess\pullback[\pt]{\objfun}\paren*{\tanvecone[\pt]} -  \sum_{\ineqidx\in\ineqset}\frac{\barrparam[]}{\pullback[\pt]{\ineqfun[\ineqidx]}\paren*{\tanvecone[\pt]}}\Hess\pullback[\pt]{\ineqfun[\ineqidx]}\paren*{\tanvecone[\pt]}}\sbra*{\tanvectwo[\pt]}}{\tanvecthr[\pt]} + \sum_{\ineqidx\in\ineqset}\frac{\barrparam[]}{{\pullback[\pt]{\ineqfun[\ineqidx]}\paren*{\tanvecone[\pt]}}^{2}}\metr[\pt]{\gradstr\pullback[\pt]{\ineqfun[\ineqidx]}\paren*{\tanvecone[\pt]}}{\tanvectwo[\pt]}\metr[\pt]{\gradstr\pullback[\pt]{\ineqfun[\ineqidx]}\paren*{\tanvecone[\pt]}}{\tanvecthr[\pt]}$ holds.}
        \end{lemma}
        \isextendedversion{
        \begin{proof}
            See \cref{appx:proofconsistency}.
        \end{proof}}{}
        Using \cref{lemm:firstdirderivmeritfun,lemm:twicedirderivmeritfun}, we derive a bound on the gap between the predicted and actual reductions.
        \begin{lemma}\label{lemm:diffpredaredfixedpt}
            Choose $\allvar=\paren*{\pt, \ineqLagmult[]}\in\strictfeasirgn\times\setRpp[\ineqdime]$ and $\barrparam[] > 0$ arbitrarily. 
            Then, there exists $\coeffone > 0$ such that 
            \isextendedversion{
            \begin{align}\label{ineq:diffpredaredatpt}
                \abs*{\pred[{\allvar, \barrparam[]}]\paren*{\dirpt}-\ared[{\barrparam[]}]\paren*{\dirpt}} \leq \coeffone\Riemnorm[\pt]{\dirpt}^{2}
            \end{align}
            }{$\abs*{\pred[{\allvar, \barrparam[]}]\paren*{\dirpt}-\ared[{\barrparam[]}]\paren*{\dirpt}} \leq \coeffone\Riemnorm[\pt]{\dirpt}^{2}$}
            for any $\dirpt\in\tanspc[\pt]\mani$ sufficiently small.
        \end{lemma}
        \begin{proof}
            See \cref{appx:proofdiffpredared}.
        \end{proof}
        Now, we prove that the iteration is successful if the current iterate is not a solution of $\barrKKTvecfld\paren*{\plchold;\barrparam[]} = 0$ for a given $\barrparam[] > 0$ and the trust region radius is sufficiently small.
        Let $\succset\subseteq\setNz$ be the set of successful iterations.
        \begin{proposition}\label{prop:existsuccidx}
            Suppose \cref{assu:Cauchydecreasing}.
            Let $\inriteridxunsucc\in\setNz\backslash\succset$ be any unsuccessful iteration of \cref{algo:RIPTRMInner} with $\barrKKTvecfld\paren*{\allvarinriterunsucc;\barrparam[]} \neq 0$ for a given $\barrparam[] > 0$.
            Then, there exists $\inriteridxsucc[]\in\succset$ such that $\inriteridxsucc[]\geq\inriteridxunsucc$. 
        \end{proposition}        
        \begin{proof}
            We first consider the test in line~15 of \cref{algo:RIPTRMInner} on the strict feasibility of the primal variable.
            Recall that, when the iterations are unsuccessful, the iterates remain the same and only the trust region radius continues to shrink, meaning that the norm of the search directions is made sufficiently small since $\Riemnorm[\ptinriter]{\dirptinriter} \leq \trradiusinriter$ holds.
            Therefore, it follows from \cref{lemm:ineqfunfeasi} that these sufficiently small directions pass the test in line~15 of \cref{algo:RIPTRMInner}.

            In the following, we focus on such iterations and additionally prove the existence of $\inriteridxsucc\in\setN$ satisfying that $\inriteridxsucc[]\geq\inriteridxunsucc$ and $\ared^{\inriteridxsucc}>\trratiothold\pred^{\inriteridxsucc}$ hold.
            We argue by contradiction.
            Suppose that $\aredinriter\leq\trratiothold\predinriter$ holds for every $\inriteridx\geq\inriteridxunsucc$.
            Since $\aredinriter\leq\trratiothold\predinriter<\frac{1}{4}\predinriter$ holds for every $\inriteridx\geq\inriteridxunsucc$, the sequence $\brc*{\trradiusinriter}_{\inriteridx}$ converges to zero as $\inriteridx$ tends to infinity according to \cref{eq:trradiusupdate}.  
            Combining this limit with $\Riemnorm[\ptinriter]{\trslininriter} = \Riemnorm[\ptinriterunsucc]{\trslininriterunsucc}$ and $\opnorm{\trsquadinriter}=\opnorm{\trsquadinriterunsucc}$ for all $\inriteridx\geq\inriteridxunsucc$ yields 
            \isextendedversion{
            \begin{equation}\label{eq:minCauchytrradius}
                \min\paren*{\trradiusinriter, \frac{\Riemnorm[\ptinriter]{\trslininriter}}{\opnorm{\trsquadinriter}}}=\min\paren*{\trradiusinriter, \frac{\Riemnorm[\ptinriterunsucc]{\trslininriterunsucc}}{\opnorm{\trsquadinriterunsucc}}} = \trradiusinriter
            \end{equation}}{
            \vspace{-1\baselineskip}
            \begin{equation}\label{eq:minCauchytrradius}
                \min\paren*{\trradiusinriter, \frac{\Riemnorm[\ptinriter]{\trslininriter}}{\opnorm{\trsquadinriter}}}=\min\paren*{\trradiusinriter, \frac{\Riemnorm[\ptinriterunsucc]{\trslininriterunsucc}}{\opnorm{\trsquadinriterunsucc}}} = \trradiusinriter
            \end{equation}
            }
            for any $\inriteridx$ sufficiently large.
            We have that, for every $\inriteridx\geq\inriteridxunsucc$ sufficiently large,
            \isextendedversion{
            \begin{align}
                \begin{split}\label{ineq:diffpredaredCauchyatfixpt}
                    &\predinriter - \aredinriter \geq \paren*{1 - \trratiothold} \predinriter\\
                    &\geq \paren*{1 - \trratiothold} \constCauchy\Riemnorm[\ptinriter]{\trslininriter}\min\paren*{\trradiusinriter, \frac{\Riemnorm[\ptinriter]{\trslininriter}}{\opnorm{\trsquadinriter}}} = \paren*{1 - \trratiothold} \constCauchy\Riemnorm[\ptinriter]{\trslininriter}\trradiusinriter,
                \end{split}
            \end{align}
            }{
            \begin{align}
                &\predinriter - \aredinriter \geq \paren*{1 - \trratiothold} \predinriter\\
                &\geq \paren*{1 - \trratiothold} \constCauchy\Riemnorm[\ptinriter]{\trslininriter}\min\paren*{\trradiusinriter, \frac{\Riemnorm[\ptinriter]{\trslininriter}}{\opnorm{\trsquadinriter}}} = \paren*{1 - \trratiothold} \constCauchy\Riemnorm[\ptinriter]{\trslininriter}\trradiusinriter,
            \end{align}
            }
            where the first inequality follows from $\aredinriter\leq\trratiothold\predinriter$, the second one from \cref{assu:Cauchydecreasing}, and the equality from \cref{eq:minCauchytrradius}.
            \isextendedversion{
            Hence, it follows from \cref{ineq:diffpredaredCauchyatfixpt,lemm:diffpredaredfixedpt}
            }{
            Using this, it follows from \cref{lemm:diffpredaredfixedpt}}
            and $\Riemnorm[\ptinriter]{\dirptinriter}\leq\trradiusinriter$ that 
            \isextendedversion{\begin{align}
                \paren*{1 - \trratiothold} \constCauchy\Riemnorm[\ptinriterunsucc]{\trslininriterunsucc}\trradiusinriter \leq \predinriter - \aredinriter \leq\abs*{\predinriter - \aredinriter} \leq \coeffone\Riemnorm[\ptinriter]{\dirptinriter}^{2} \leq \coeffone\paren*{\trradiusinriter}^{2}
            \end{align}}{
            \vspace{-1.5\baselineskip}
            \begin{align}
                \paren*{1 - \trratiothold} \constCauchy\Riemnorm[\ptinriterunsucc]{\trslininriterunsucc}\trradiusinriter \leq \predinriter - \aredinriter \leq\abs*{\predinriter - \aredinriter} \leq \coeffone\Riemnorm[\ptinriter]{\dirptinriter}^{2} \leq \coeffone\paren*{\trradiusinriter}^{2}
            \end{align}
            }
            for all $\inriteridx\geq\inriteridxunsucc$ sufficiently large.
            This is, however, a contradiction since the left-hand side is $\bigO[\trradiusinriter]$, whereas the right-hand side is $\bigO[\paren*{\trradiusinriter}^{2}]$.
            The proof is complete.
        \end{proof}


    \subsection{Global convergence of \texorpdfstring{\cref{algo:RIPTRMInner}}{inner iteration}}\label{subsec:globconvinner}
        In this subsection, we prove the global convergence of \cref{algo:RIPTRMInner} to a solution of \cref{eq:barrKKTequation} with the second-order condition, which ensures that \cref{algo:RIPTRMInner} terminates in a finite number of iterations to satisfy the stopping conditions~\cref{eq:stopcond}.
        For the global convergence, we additionally assume the following:
        \begin{assumption}\label{assu:retrsecondord}
            The retraction $\retr[]$ is second order.
        \end{assumption}
        \begin{assumption}\label{assu:radialLCone}
            The functions $\brc*{\ineqfun[\ineqidx]}_{\ineqidx\in\ineqset}$ are \rLCone{} on some $\opensubset\subseteq\mani$ with $\brc*{\ptinriter}_{\inriteridx}\subseteq\opensubset$.
        \end{assumption}
        \begin{assumption}\label{assu:radialLCtwo}
            The functions $\objfun$ and $\brc*{\ineqfun[\ineqidx]}_{\ineqidx\in\ineqset}$ are \rLCtwo{}  on some $\opensubset\subseteq\mani$ with $\brc*{\ptinriter}_{\inriteridx}\subseteq\opensubset$.
        \end{assumption}
        \begin{assumption}\label{assu:ineqfunbounded}
            The following hold:
            \begin{enumerate}
                \item The sequence $\brc*{\ineqfun[\ineqidx]\paren*{\ptinriter}}_{\inriteridx}$ is bounded above for every $\ineqidx\in\ineqset$.\label{assu:ineqfunboundedvalue}
                \item The sequences $\brc*{\Riemnorm[\ptinriter]{\gradstr\ineqfun[\ineqidx]\paren*{\pt}}}_{\inriteridx}$ and $\brc*{\opnorm{\Hess\ineqfun[\ineqidx]\paren*{\pt}}}_{\inriteridx}$ are bounded above for every $\ineqidx\in\ineqset$.\label{assu:ineqfunboundedgradhess}
            \end{enumerate}
        \end{assumption}
        \begin{assumption}\label{assu:objfunbounded}
        The following hold:
            \begin{enumerate}
                \item The sequence $\brc*{\objfun\paren*{\ptinriter}}_{\inriteridx}$ is bounded below. \label{assu:objfunboundedbelow}
                \item The sequence $\brc*{\opnorm{\Hess\objfun\paren*{\ptinriter}}}_{\inriteridx}$ is bounded above. \label{assu:objfunHessbounded}
            \end{enumerate}
        \end{assumption}
        \begin{assumption}\label{assu:pullbackineqfuncont}
            For each $\ineqidx\in\ineqset$ and any $\epsone > 0$, there exists $\thldtrradiusthr > 0$ such that, for any $\inriteridx\in\setNz$ and all $\tanvecone[\ptinriter]\in\tanspc[\ptinriter]\mani$ with $\Riemnorm[\ptinriter]{\tanvecone[\ptinriter]} \leq \thldtrradiusthr$,
            \isextendedversion{
            \begin{equation}\label{eq:pullbackineqfuncont}
                \abs*{\pullback[\ptinriter]{{\ineqfun[\ineqidx]}}\paren*{\tanvecone[\ptinriter]} - \pullback[\ptinriter]{{\ineqfun[\ineqidx]}}\paren*{\zerovec[\ptinriter]}} \leq \epsone.
            \end{equation}}{$\abs*{\pullback[\ptinriter]{{\ineqfun[\ineqidx]}}\paren*{\tanvecone[\ptinriter]} - \pullback[\ptinriter]{{\ineqfun[\ineqidx]}}\paren*{\zerovec[\ptinriter]}} \leq \epsone$.}
        \end{assumption}
        \begin{assumption}\label{assu:ineqLagmultbounded}
            The sequence $\brc*{\Riemnorm[]{\ineqLagmultinriter[]}}_{\inriteridx}$ is bounded above.
        \end{assumption}
        \begin{assumption}\label{assu:ineqLagmultconverged}
            $\lim_{\inriteridx\to\infty}\norm{\ineqLagmultinriter[]-\barrparam[]\inv{\Ineqfunmat[]\paren*{\ptinriter}}\onevec}=0$.
        \end{assumption}
        \begin{assumption}\label{assu:retrdistbound}
            There exist positive scalars $\epstwo, \deltatwo \in \setRpp[]$ such that, for all $\inriteridx\in\setNz$ and all $\tanvecone[\ptinriter]\in\tanspc[\ptinriter]\mani$ with $\Riemnorm[\ptinriter]{\tanvecone[\ptinriter]}\leq\epstwo$,
            \begin{equation}\label{ineq:retrdistbound}
                \Riemdist{\ptinriter}{\retr[\ptinriter]\paren*{\tanvecone[\ptinriter]}}\leq\deltatwo\Riemnorm[\ptinriter]{\tanvecone[\ptinriter]}.
            \end{equation}
        \end{assumption}
        \begin{assumption}\label{assu:Eigendecreasing}
            There exists $\constEigen > 0$ such that, 
            for any $\epseigen>0$ and any $\inriteridx\in\setNz$ with $\mineigval\sbra*{\trsquadinriter} < -\epseigen$, the search direction $\dirptinriter\in\tanspc[\ptinriter]\mani$ satisfies
            \isextendedversion{
            \begin{equation}\label{ineq:Eigendecreasing}
                \modelfuninriter\paren*{\zerovec[\ptinriter]} - \modelfuninriter\paren*{\dirptinriter} \geq \constEigen\epseigen\paren*{\trradiusinriter}^{2}.
            \end{equation}}{$\modelfuninriter\paren*{\zerovec[\ptinriter]} - \modelfuninriter\paren*{\dirptinriter} \geq \constEigen\epseigen\paren*{\trradiusinriter}^{2}$.}           
        \end{assumption}
        \isextendedversion{
        In \cref{subsec:computdualvar}, we will discuss \cref{assu:ineqLagmultbounded,assu:ineqLagmultconverged} when using clipping to update the dual variable.
        We defer the discussion for the other assumptions to \cref{appx:suffcondassus}.
        In short, 
        }{}
        \cref{assu:retrsecondord} is not restrictive as in \cref{subsec:notationRiemggeo}.
        \cref{assu:radialLCone,assu:radialLCtwo,assu:ineqfunbounded,assu:objfunbounded,assu:pullbackineqfuncont,assu:retrdistbound} are fulfilled if the generated sequence $\brc*{\ptinriter}_{\inriteridx}$ is bounded and all functions $\objfun, \brc*{\ineqfun[\ineqidx]}_{\ineqidx\in\ineqset}$ are of class $C^3$, for example.
        The \eigenstep{}~\isextendedversion{}{\cite[Lemma~6.16]{Boumal23IntroOptimSmthMani}}
        and the \exactsolution{}, the global optimum of \cref{prob:TRS}, satisfy \cref{assu:Eigendecreasing} with $\constEigen=\frac{1}{2}$\isextendedversion{; see \cref{sec:searchdir} for the detail.}{.} 
        \isextendedversion{}{
        The clipping \cref{def:clippingupdate} satisfies \cref{assu:ineqLagmultbounded,assu:ineqLagmultconverged} under other assumptions.
        Since we can prove it in a similar manner to the Euclidean version~\cite[Section~4.3]{Connetal2000PrimalDualTRAlgoforNonconvexNLP}, 
        we omit the detail.\footnote{
        The reader can find the detail in 
        the extended version of the paper~\cite[Appendix E]{Obaraetal2025APrimalDualIPTRMfor2ndOrdStnryPtofRiemIneqCstrOptimProbs}.}}
        We also note that all the assumptions are standard in the literature.
      
        For the global convergence analysis, we first investigate the properties of the sequence of the merit function's values.
        By enumerating the elements of $\succset\subseteq\setNz$ in increasing order, we obtain the sequence of successful iteration indices $\brc*{\inrsuccsubseqiter}_{\succsubseqidx\in\setNz}$.
        For example, $\inriteridx_{1}$ denotes the first successful iteration and $\inriteridx_{2}$ the second.
        Note that any iterate between $\inrsuccsubseqiter$ and $\inrsuccsubseqiterp$, if it exists, is unsuccessful by definition; that is, the $\paren*{1+\inrsuccsubseqiter}$-, $\paren*{2+\inrsuccsubseqiter}$-, $\ldots$, $\paren*{\inrsuccsubseqiterp-1}$-th iterations are all unsuccessful.
        The following holds by the definitions of the successful iterates and the lower-boundedness of $\brc*{\meritfun[{\barrparam[]}]\paren*{\ptinriter}}_{\inriteridx}$.
        \begin{lemma}\label{lemm:meritfunanal}
            Suppose \cref{assu:Cauchydecreasing}.
            The following hold:
            \begin{enumerate}
                \item The sequence $\brc*{\meritfun[{\barrparam[]}]\paren*{\ptinriter}}_{\inriteridx}$ is monotonically non-increasing for a given $\barrparam[] > 0$.\label{lemm:meritfunmonnonincr}
                \item Under \assuenumicref{assu:ineqfunbounded}{assu:ineqfunboundedvalue} and \assuenumicref{assu:objfunbounded}{assu:objfunboundedbelow}, $\brc*{\meritfun[{\barrparam[]}]\paren*{\ptinriter}}_{\inriteridx}$ is convergent for a given $\barrparam[] > 0$.\label{lemm:meritfunconv}
            \end{enumerate}
        \end{lemma}
        \isextendedversion{        \begin{proof}
            See \cref{appx:proofinnerglobal}.
        \end{proof}}{
        \begin{proof}
            We defer the complete proof to \cite[Appendix~D.2]{Obaraetal2025APrimalDualIPTRMfor2ndOrdStnryPtofRiemIneqCstrOptimProbs} and here provide an overview.
            By the update rule~\cref{eq:primalupdate}, we have $\meritfun[{\barrparam[]}]\paren*{\ptinrsuccsubseqiter}
            - \meritfun[{\barrparam[]}]\paren*{\ptinrsuccsubseqiterp}
            \geq \trratiothold\paren*{\modelfuninrsuccsubseqiter\paren*{\zerovec[\ptinrsuccsubseqiter]} - \modelfuninrsuccsubseqiter\paren*{\dirinrsuccsubseqiter}}$ for all $\succsubseqidx\in\setN$.
            Under \cref{assu:Cauchydecreasing}, the right-hand side is nonnegative, which yields \cref{lemm:meritfunmonnonincr}.
            Moreover, under \assuenumicref{assu:ineqfunbounded}{assu:ineqfunboundedvalue} and \assuenumicref{assu:objfunbounded}{assu:objfunboundedbelow}, the sequence $\brc*{\meritfun[{\barrparam[]}]\paren*{\ptinrsuccsubseqiter}}_{\succsubseqidx}$ is bounded below.
            Therefore, \cref{lemm:meritfunconv} follows from the monotone convergence theorem.
        \end{proof}
        }
        We next derive positive lower bounds on $\brc*{\ineqfun[\ineqidx]\paren*{\ptinriter}}_{\inriteridx}$ and its sufficiently small neighborhoods by the upper-boundedness of $\brc*{\meritfun[{\barrparam[]}]\paren*{\ptinriter}}_{\inriteridx}$ and the assumptions.
        \begin{lemma}\label{lemm:ineqfunbound}
            Under \cref{assu:Cauchydecreasing}, \assuenumicref{assu:ineqfunbounded}{assu:ineqfunboundedvalue}, and \assuenumicref{assu:objfunbounded}{assu:objfunboundedbelow},
            the following hold:
            \begin{enumerate}
                \item There exists $\epsfiv > 0$ such that $\ineqfun[\ineqidx]\paren*{\ptinriter} \geq \epsfiv$ holds for any $\ineqidx\in\ineqset$ and all $\inriteridx\in\setNz$.\label{lemm:ienqfunptlowbd}
                \item Under \cref{assu:pullbackineqfuncont}, there exist $\epsthr > 0$ and $\deltathr > 0$ such that, for any $\inriteridx\in\setNz$ and all $\tanvecone[\ptinriter]\in\tanspc[\ptinriter]\mani$ with $\Riemnorm[\ptinriter]{\tanvecone[\ptinriter]}\leq\deltathr$, $\pullback[\ptinriter]{\ineqfun[\ineqidx]}\paren*{\tanvecone[\ptinriter]} \geq \epsthr$ holds.\label{lemm:ienqfununiflowbd}
            \end{enumerate}
        \end{lemma}
        \isextendedversion{        \begin{proof}
            See \cref{appx:proofinnerglobal}.
        \end{proof}}{
        \begin{proof}
            We defer the complete proof to \cite[Appendix~D.2]{Obaraetal2025APrimalDualIPTRMfor2ndOrdStnryPtofRiemIneqCstrOptimProbs} and provide here a brief overview.
            By the definition of the log barrier function, we have $\sum_{\ineqidx\in\ineqset} \log\ineqfun[\ineqidx]\paren*{\ptinriter} = \inv{\barrparam[]} \paren*{\objfun\paren*{\ptinriter} - \meritfun[{\barrparam[]}]\paren*{\ptinriter}}$.
            Since the right-hand side is bounded below by \assuenumicref{assu:objfunbounded}{assu:objfunboundedbelow} and \enumicref{lemm:meritfunanal}{lemm:meritfunmonnonincr}, the equation implies \cref{lemm:ienqfunptlowbd}.
            Combining \cref{lemm:ienqfunptlowbd} with \cref{assu:pullbackineqfuncont} then yields \cref{lemm:ienqfununiflowbd}.
        \end{proof}
        }
        Note that \cref{lemm:ineqfunbound} provides a positive lower bound on the values of the inequality constraints during the inner iterations, while \cref{lemm:ineqfunfeasi} ensures the positivity of the inequality constraints around a fixed point.
        We use \cref{lemm:ineqfunbound} to prove the boundedness of several sequences, one of which is $\brc*{\opnorm{\trsquadinriter}}_{\inriteridx}$ in the following lemma.
        \begin{lemma}\label{lemm:trsquadbounded}
            Under \cref{assu:Cauchydecreasing,assu:objfunbounded,assu:ineqLagmultbounded,assu:ineqfunbounded},
             there exists $\tholdvaltwo > 0$ such that $\opnorm{\trsquadinriter} \leq \tholdvaltwo$ for all $\inriteridx\in\setNz$.
        \end{lemma}
        \isextendedversion{        \begin{proof}
            See \cref{appx:proofinnerglobal}.
        \end{proof}}{
        \begin{proof}    
            By definition, we have 
            \begin{align}
                &\Riemnorm[\ptinriter]{\trsquadinriter\sbra*{\tanvecone[\ptinriter]}}\leq \opnorm{\Hess\objfun\paren*{\ptinriter}}\Riemnorm[\ptinriter]{\tanvecone[\ptinriter]} \\
                &\quad + \sum_{\ineqidx\in\ineqset} \ineqLagmultinriter[\ineqidx]\opnorm{\Hess\ineqfun[\ineqidx]\paren*{\ptinriter}}\Riemnorm[\ptinriter]{\tanvecone[\ptinriter]} + \sum_{\ineqidx\in\ineqset} \frac{\ineqLagmultinriter[\ineqidx]}{\ineqfun[\ineqidx]\paren*{\ptinriter}}\Riemnorm[\ptinriter]{\gradstr\ineqfun[\ineqidx]\paren*{\ptinriter}}^{2}\Riemnorm[\ptinriter]{\tanvecone[\ptinriter]}^{2}
            \end{align}
            for any $\inriteridx\in\setNz$ and all $\tanvecone[\ptinriter]\in\tanspc[\ptinriter]\mani$ with $\Riemnorm[\ptinriter]{\tanvecone[\ptinriter]} \leq 1$.
            We then prove that the right hand is bounded above by analyzing each term.
            See \cite[Appendix D.2]{Obaraetal2025APrimalDualIPTRMfor2ndOrdStnryPtofRiemIneqCstrOptimProbs} for the complete proof.
        \end{proof}
        }
        The following lemma is a crucial ingredient for proving the global convergence.
        \begin{lemma}\label{lemm:diffpredared}
            Suppose that \cref{assu:retrsecondord,assu:radialLCone,assu:radialLCtwo,assu:Cauchydecreasing,assu:ineqfunbounded}, \assuenumiCref{assu:objfunbounded}{assu:objfunboundedbelow}, \Cref{assu:pullbackineqfuncont,assu:ineqLagmultbounded} hold.
            The following hold:
            \begin{enumerate}
                \item There exist $\thldtrradiusfou > 0$ and $\coefftwo > 0$ such that, for all $\inriteridx\in\setNz$, if the search direction $\dirptinriter\in\tanspc[\ptinriter]\mani$ satisfies $\Riemnorm[\ptinriter]{\dirptinriter} \leq \thldtrradiusfou$, then 
                \isextendedversion{
                \begin{equation}\label{ineq:diffpredaredell}
                    \abs*{\predinriter-\aredinriter} \leq \coefftwo \Riemnorm[\ptinriter]{\dirptinriter}^{2}.
                \end{equation}
                }{
                $\abs*{\predinriter-\aredinriter} \leq \coefftwo \Riemnorm[\ptinriter]{\dirptinriter}^{2}$.}
                \label{lemm:diffpredareduniform}
                \item Suppose \cref{assu:ineqLagmultconverged} in addition.
                Choose $\coefffou > 0$ arbitrarily.
                Then, there exist $\thldtrradiussix > 0$ and $\tholdidxtwo \in \setNz$ such that, for all $\inriteridx\geq\tholdidxtwo$, if the search direction $\dirptinriter\in\tanspc[\ptinriter]\mani$ satisfies 
                \isextendedversion{
                \begin{align}\label{ineq:diffpredaredsmallo}
                    \abs*{\predinriter-\aredinriter} \leq \coefffou\Riemnorm[\ptinriter]{\dirptinriter}^{2}.
                \end{align}
                }{
                $\Riemnorm[\ptinriter]{\dirptinriter} \leq \thldtrradiussix$, then 
                $\abs*{\predinriter-\aredinriter} \leq \coefffou\Riemnorm[\ptinriter]{\dirptinriter}^{2}$.}
                \label{lemm:diffpredaredsmallo}
            \end{enumerate}
        \end{lemma}
        \begin{proof}
            See \cref{appx:proofdiffpredared}.
        \end{proof}
        In contrast to Euclidean optimization, we make use of the retraction in the Riemannian setting, which may cause the curve $\tmethr\mapsto\retr[\pt]\paren*{\tmethr\dir[\pt]}$ to reach the boundary of the feasible region and potentially leads to a division by zero.
        Additionally, since the pullbacks of the Riemannian gradient and Hessian are defined on $\tanspc[\pt]\mani$, which varies depending on $\inriteridx$, their norms may diverge as $\inriteridx$ tends to infinity.
        Our analysis in the proof appropriately addresses these concerns using \cref{lemm:ineqfunbound} and \cref{assu:retrsecondord,assu:radialLCone,assu:radialLCtwo}.
        Compared with \cref{lemm:diffpredaredfixedpt}, \cref{lemm:diffpredared} estimates the gap between the predicted and actual reductions during the inner iteration.
        In particular, by additionally supposing \cref{assu:ineqLagmultconverged}, we obtain a bound whose coefficient can be made arbitrarily small as proved in \enumicref{lemm:diffpredared}{lemm:diffpredaredsmallo}.
        We will use \enumicref{lemm:diffpredared}{lemm:diffpredareduniform} and \enumicref{lemm:diffpredared}{lemm:diffpredaredsmallo} to prove the global convergence to an \AKKT{} point and an \SOSP{}, respectively.


        In the following theorem, we analyze the limit inferior of $\brc*{\Riemnorm[\ptinriter]{\trslininriter}}_{\inriteridx}$ that will be used for analyzing the entire sequence.
        \begin{theorem}\label{theo:liminftrslininriter}
            Under \cref{assu:retrsecondord,assu:radialLCone,assu:radialLCtwo,assu:Cauchydecreasing,assu:ineqfunbounded,assu:objfunbounded,assu:pullbackineqfuncont,assu:ineqLagmultbounded}, $\liminf_{\inriteridx\to\infty}\Riemnorm[\ptinriter]{\trslininriter} = 0$ holds.
        \end{theorem}
        \begin{proof}
            We argue by contradiction.
            Suppose that there exists $\epsfou > 0$ and $\tholdidxone\in\setNz$ such that $\Riemnorm[\ptinriter]{\trslininriter} \geq \epsfou$ holds for all $\inriteridx\geq\tholdidxone$.
            It follows from \cref{assu:Cauchydecreasing,lemm:trsquadbounded} that, for all $\inriteridx\geq\tholdidxone$, 
            \begin{equation}\label{ineq:predCauchy}
                \predinriter \geq \constCauchy\Riemnorm[\ptinriter]{\trslininriter}\min\paren*{\trradiusinriter, \frac{\Riemnorm[\ptinriter]{\trslininriter}}{\opnorm{\trsquadinriter}}} \geq \constCauchy\epsfou\min\paren*{\trradiusinriter, \frac{\epsfou}{\tholdvaltwo}} > 0.
            \end{equation}
            Therefore, letting $\thldtrradiusfou > 0$ be the threshold in \enumicref{lemm:diffpredared}{lemm:diffpredareduniform},
            we have
            \isextendedversion{
            \begin{align}
                \abs*{\trratioinriter - 1} 
                = \frac{\abs*{\predinriter - \aredinriter}}{\abs*{\predinriter}} \leq \frac{\coefftwo\Riemnorm[\ptinriter]{\dirptinriter}^{2}}{\constCauchy\epsfou\min\paren*{\trradiusinriter, \frac{\epsfou}{\tholdvaltwo}}} \leq \frac{\coefftwo\paren*{\trradiusinriter}^{2}}{\constCauchy\epsfou\min\paren*{\trradiusinriter, \frac{\epsfou}{\tholdvaltwo}}}
            \end{align}}{
            \begin{align}
                \abs*{\trratioinriter - 1} 
                = \frac{\abs*{\predinriter - \aredinriter}}{\abs*{\predinriter}} \leq \frac{\coefftwo\Riemnorm[\ptinriter]{\dirptinriter}^{2}}{\constCauchy\epsfou\min\paren*{\trradiusinriter, \frac{\epsfou}{\tholdvaltwo}}} \leq \frac{\coefftwo\paren*{\trradiusinriter}^{2}}{\constCauchy\epsfou\min\paren*{\trradiusinriter, \frac{\epsfou}{\tholdvaltwo}}}
            \end{align}
                }
            whenever $\inriteridx \geq \tholdidxone$ and $\trradiusinriter \leq \thldtrradiusfou$, where the first inequality follows from \enumicref{lemm:diffpredared}{lemm:diffpredareduniform} and \cref{ineq:predCauchy} and the second one from $\Riemnorm[\ptinriter]{\dirptinriter} \leq \trradiusinriter$.
            Define $\thldtrradiusfiv\coloneqq\min\paren*{\thldtrradiusfou, \frac{\constCauchy\epsfou}{2\coefftwo}, \frac{\epsfou}{\tholdvaltwo}}$.
            If $\trradiusinriter \leq \thldtrradiusfiv$, we have $\min\paren*{\trradiusinriter, \frac{\epsfou}{\tholdvaltwo}}=\trradiusinriter$ and hence $\abs*{\trratioinriter - 1} \leq \frac{\coefftwo\trradiusinriter}{\constCauchy\epsfou} \leq \frac{1}{2}$.
            This implies $\trratioinriter\geq\frac{1}{2}>\frac{1}{4}>\trratiothold$, and thus the update $\trradiusinriterp=\frac{1}{4}\trradiusinriter$ can occur only if $\trradiusinriter > \thldtrradiusfiv$.
            Therefore, for all $\inriteridx\geq\tholdidxone$, we have
            \begin{equation}\label{ineq:trradiuslowboundCauchy}
                \trradiusinriter\geq\min\paren*{\trradiusinritertholdone, \frac{\thldtrradiusfiv}{4}}.
            \end{equation} 
            On the other hand, recall that $\succset$ is the set of the successful iterations, and $\brc*{\ptinrsuccsubseqiter}_{\succsubseqidx}$ is the ordered sequence of $\succset$.
            When $\abs*{\succset}$ is finite, all the sufficiently large iterations become unsuccessful.
            Thus, it follows that $\lim_{\inriteridx\to\infty}\trradiusinriter = 0$, which contradicts \cref{ineq:trradiuslowboundCauchy}.
            In the following, we consider the case where $\card[\succset]$ is infinite.
            Recall that $\pt^{1+\inrsuccsubseqiter}=\ptinrsuccsubseqiterp$ holds since the $\paren*{1+\inrsuccsubseqiter}$-, $\paren*{2+\inrsuccsubseqiter}$-, $\ldots$, $\paren*{\inrsuccsubseqiterp -1}$-th iterations are all unsuccessful.
            For all $\succsubseqidx\in\setNz$ satisfying $\inrsuccsubseqiter\geq\tholdidxone$, we obtain
            \begin{equation}\label{ineq:diffmeritfunCauchy}
                \meritfun[{\barrparam[]}]\paren*{\ptinrsuccsubseqiter} - \meritfun[{\barrparam[]}]\paren*{\ptinrsuccsubseqiterp} = \meritfun[{\barrparam[]}]\paren*{\ptinrsuccsubseqiter} - \meritfun[{\barrparam[]}]\paren*{\pt^{1+\inrsuccsubseqiter}} \geq \trratiothold\pred[\inrsuccsubseqiter]
                \geq \trratiothold\constCauchy\epsfou\min\paren*{\trradiusinrsuccsubseqiter, \frac{\epsfou}{\tholdvaltwo}} \geq 0,
            \end{equation}
            where the equality follows from $\ptinrsuccsubseqiterp=\pt^{1+\inrsuccsubseqiter}$, the first inequality from the fact that the $\inrsuccsubseqiter$-th iterate is successful, and the second one from \cref{ineq:predCauchy}.
            From \enumicref{lemm:meritfunanal}{lemm:meritfunconv}, $\brc*{\meritfun[{\barrparam[]}]\paren*{\ptinriter}}_{\inriteridx}$ is a Cauchy sequence.
            Thus, \cref{ineq:diffmeritfunCauchy} implies $\liminf_{\inrsuccsubseqiter\to\infty}\trradiusinrsuccsubseqiter = 0$,
            which contradicts \cref{ineq:trradiuslowboundCauchy}.
            
            Now, we have that, for any  $\epsfou > 0$ and any index $\tholdidxone\in\setNz$, there exists $\inritertholdoneidx\geq\tholdidxone$ such that $\Riemnorm[\ptinritertholdone]{\trslininritertholdone} < \epsfou$ holds.
            This implies $\liminf_{\inriteridx\to\infty}\Riemnorm[\ptinriter]{\trslininriter} = 0$.
        \end{proof} 

        Now, we present the following theorem on the limit of $\brc*{\Riemnorm[\ptinriter]{\trslininriter}}_{\inriteridx}$.
        Recall that $\succset$ is the set of the successful iterations.
        \begin{theorem}\label{theo:limtrslininriter}
            Under \cref{assu:retrsecondord,assu:radialLCone,assu:radialLCtwo,assu:Cauchydecreasing,assu:ineqfunbounded,assu:objfunbounded,assu:pullbackineqfuncont,assu:ineqLagmultbounded,assu:retrdistbound},
            $\lim_{\inriteridx\to\infty}\Riemnorm[\ptinriter]{\trslininriter} = 0$ holds.
        \end{theorem}
        \begin{proof}
            If $\abs*{\succset}$ is finite, the point $\ptinriter\in\strictfeasirgn$ remains the same for all $\inriteridx\in\setNz$ sufficiently large.
            Hence, from \cref{theo:liminftrslininriter}, it follows that $\lim_{\inriteridx\to\infty}\Riemnorm[\ptinriter]{\trslininriter} = 0$.
            In the following, we consider the case where $\abs*{\succset}$ is infinite.
            Suppose that there exists an infinite subsequence of $\succset$, denoted by $\brc*{\ptinrsuccsubseqtwoiter}_{\succsubseqidxtwo}$, and a constant $\epssix > 0$ such that $\Riemnorm[\ptinrsuccsubseqtwoiter]{\trslininrsuccsubseqtwoiter} \geq 3 \epssix$ holds for all $\succsubseqidxtwo\in\setNz$.
            From \cref{theo:liminftrslininriter}, for any $\succsubseqidxtwo\in\setNz$, there exists the first index $\inrsuccsubseqthriter\in\setNz$ satisfying $\inrsuccsubseqthriter > \inrsuccsubseqtwoiter$ and $\Riemnorm[\ptinrsuccsubseqthriter]{\trslininrsuccsubseqthriter} < \epssix$.
            Define $\succprocset[\succsubseqidxtwo]\coloneqq \brc*{\inrsuccsubseqfouiter\in\succset \colon \inrsuccsubseqtwoiter \leq \inrsuccsubseqfouiter < \inrsuccsubseqthriter}$ for each $\succsubseqidxtwo\in\setNz$ and $\succprocset\coloneqq\bigcup_{\succsubseqidxtwo\in\setNz}\succprocset[\succsubseqidxtwo]$.
            Notice $\succprocset[\succsubseqidxtwo]$ is nonempty since $\inrsuccsubseqtwoiter\in\succprocset[\succsubseqidxtwo]\subseteq\succset$ holds.
            Note also that, for any $\inrsuccsubseqfouiter\in\succprocset$, $\Riemnorm[\ptinrsuccsubseqfouiter]{\trslininrsuccsubseqfouiter}\geq\epssix$ holds by definition.
            For any $\inrsuccsubseqfouiter\in\succprocset$, since the $\inrsuccsubseqfouiter$-th iterate is successful, we have
            \begin{align}
                \begin{split}\label{ineq:meritfundiffCauchysuccsubseqiter}
                    &\meritfun[{\barrparam[]}]\paren*{\ptinrsuccsubseqfouiter} - \meritfun[{\barrparam[]}]\paren*{\ptinrsuccsubseqfouiterp} \geq \trratiothold\paren*{\modelfuninrsuccsubseqfouiter\paren*{\zerovec[\ptinrsuccsubseqfouiter]} - \modelfuninrsuccsubseqfouiter\paren*{\dirptinrsuccsubseqfouiter}}\\
                    &\geq \trratiothold\constCauchy\Riemnorm[\ptinrsuccsubseqfouiter]{\trslininrsuccsubseqfouiter}\min\paren*{\trradiusinrsuccsubseqfouiter, \frac{\Riemnorm[\ptinrsuccsubseqfouiter]{\trslininrsuccsubseqfouiter}}{\opnorm{\trsquadinrsuccsubseqfouiter}}} \geq \trratiothold\constCauchy\epssix\min\paren*{\trradiusinrsuccsubseqfouiter, \frac{\epssix}{\tholdvaltwo}} \geq 0,
                \end{split}
            \end{align}
            where the second inequality follows from \cref{assu:Cauchydecreasing} and the third one from \cref{lemm:trsquadbounded} and the definition of the set $\succprocset$.
            Thus, it holds by \enumicref{lemm:meritfunanal}{lemm:meritfunconv} that the term $\meritfun[{\barrparam[]}]\paren*{\ptinrsuccsubseqfouiter} - \meritfun[{\barrparam[]}]\paren*{\ptinrsuccsubseqfouiterp}$
            converges to zero as $\inrsuccsubseqfouiter\in\succprocset\to\infty$ and hence $\lim_{\inrsuccsubseqfouiter\in\succprocset\to\infty}\trradiusinrsuccsubseqfouiter=0$.
            For all $\inrsuccsubseqfouiter\in\succprocset$ sufficiently large, since $\min\paren*{\trradiusinrsuccsubseqfouiter, \frac{\epssix}{\tholdvaltwo}}=\trradiusinrsuccsubseqfouiter$ holds, \cref{ineq:meritfundiffCauchysuccsubseqiter} implies
            \begin{equation}\label{ineq:trradiusboundmeritfundiff}
                \trradiusinrsuccsubseqfouiter \leq \frac{1}{\trratiothold\constCauchy\epssix}\paren*{\meritfun[{\barrparam[]}]\paren*{\ptinrsuccsubseqfouiter} - \meritfun[{\barrparam[]}]\paren*{\ptinrsuccsubseqfouiterp}}.
            \end{equation}
            Additionally, note that \cref{ineq:retrdistbound} holds for all $\inrsuccsubseqfouiter\in\succprocset$ sufficiently large from \cref{assu:retrdistbound}.
            We have
            \isextendedversion{
            \begin{align}
                &\Riemdist{\ptinrsuccsubseqtwoiter}{\ptinrsuccsubseqthriter} \leq \sum_{\dummyinriteridx=\inrsuccsubseqtwoiter}^{\inrsuccsubseqthriter-1}\Riemdist{\ptdummyinriter}{\ptdummyinriterp} = \sum_{\inrsuccsubseqfouiter\in\succprocset[\succsubseqidxtwo]}\Riemdist{\ptinrsuccsubseqfouiter}{\retr[\ptinrsuccsubseqfouiter]\paren*{\dirptinrsuccsubseqfouiter}}\\
                &\quad\leq \sum_{\inrsuccsubseqfouiter\in\succprocset[\succsubseqidxtwo]}\deltatwo\trradiusinrsuccsubseqfouiter \leq \sum_{\inrsuccsubseqfouiter\in\succprocset[\succsubseqidxtwo]}\frac{\deltatwo}{\trratiothold\constCauchy\epssix}\paren*{\meritfun[{\barrparam[]}]\paren*{\ptinrsuccsubseqfouiter} - \meritfun[{\barrparam[]}]\paren*{\ptinrsuccsubseqfouiterp}}\\
                &\quad= \sum_{\dummyinriteridx=\inrsuccsubseqtwoiter}^{\inrsuccsubseqthriter-1}\frac{\deltatwo}{\trratiothold\constCauchy\epssix}\paren*{\meritfun[{\barrparam[]}]\paren*{\ptdummyinriter} - \meritfun[{\barrparam[]}]\paren*{\ptdummyinriterp}} = \frac{\deltatwo}{\trratiothold\constCauchy\epssix}\paren*{\meritfun[{\barrparam[]}]\paren*{\ptinrsuccsubseqtwoiter} - \meritfun[{\barrparam[]}]\paren*{\ptinrsuccsubseqthriter}}
            \end{align}
            }{\begin{align}
                &\Riemdist{\ptinrsuccsubseqtwoiter}{\ptinrsuccsubseqthriter} \leq \sum_{\dummyinriteridx=\inrsuccsubseqtwoiter}^{\inrsuccsubseqthriter-1}\Riemdist{\ptdummyinriter}{\ptdummyinriterp} = \sum_{\inrsuccsubseqfouiter\in\succprocset[\succsubseqidxtwo]}\Riemdist{\ptinrsuccsubseqfouiter}{\retr[\ptinrsuccsubseqfouiter]\paren*{\dirptinrsuccsubseqfouiter}} \leq \sum_{\inrsuccsubseqfouiter\in\succprocset[\succsubseqidxtwo]} \deltatwo\trradiusinrsuccsubseqfouiter\\ 
                &\leq \sum_{\inrsuccsubseqfouiter\in\succprocset[\succsubseqidxtwo]}\frac{\deltatwo}{\trratiothold\constCauchy\epssix}\paren*{\meritfun[{\barrparam[]}]\paren*{\ptinrsuccsubseqfouiter} - \meritfun[{\barrparam[]}]\paren*{\ptinrsuccsubseqfouiterp}} = \sum_{\dummyinriteridx=\inrsuccsubseqtwoiter}^{\inrsuccsubseqthriter-1}\frac{\deltatwo}{\trratiothold\constCauchy\epssix}\paren*{\meritfun[{\barrparam[]}]\paren*{\ptdummyinriter} - \meritfun[{\barrparam[]}]\paren*{\ptdummyinriterp}}\\[-1em]
                &= \frac{\deltatwo}{\trratiothold\constCauchy\epssix}\paren*{\meritfun[{\barrparam[]}]\paren*{\ptinrsuccsubseqtwoiter} - \meritfun[{\barrparam[]}]\paren*{\ptinrsuccsubseqthriter}}
            \end{align}}
            for any $\succsubseqidxtwo\in\setNz$ sufficiently large, where the first and the second equalities follow from the definition of $\succprocset[\succsubseqidxtwo]$,
            the second inequality from \cref{ineq:retrdistbound}, and the third one from \cref{ineq:trradiusboundmeritfundiff}.
            From \enumicref{lemm:meritfunanal}{lemm:meritfunconv}, the right-hand side converges to zero as $\succsubseqidxtwo\to\infty$, which means $\lim_{\succsubseqidxtwo\to\infty}\Riemdist{\ptinrsuccsubseqtwoiter}{\ptinrsuccsubseqthriter}=0$.
            Therefore, from the continuities of $\trslin[{\barrparam[]}]$ 
            \isextendedversion{
            defined in \cref{eq:trslindef}}{defined below \cref{eq:dualNewtoneq} in \cref{subsec:inneriter}}
            and the norm, we obtain
            \isextendedversion{
            \begin{align}
                \abs*{\Riemnorm[\ptinrsuccsubseqtwoiter]{\trslininrsuccsubseqtwoiter} - \Riemnorm[\ptinrsuccsubseqthriter]{\trslininrsuccsubseqthriter}} \leq \epssix
            \end{align}}{$\abs*{\Riemnorm[\ptinrsuccsubseqtwoiter]{\trslininrsuccsubseqtwoiter} - \Riemnorm[\ptinrsuccsubseqthriter]{\trslininrsuccsubseqthriter}} \leq \epssix$}
            for all $\succsubseqidxtwo\in\setNz$ sufficiently large, which,
            together with the definitions of $\inrsuccsubseqtwoiter$ and $\inrsuccsubseqthriter$, implies that
            \isextendedversion{
            \begin{align}
                2\epssix = 3\epssix - \epssix < \Riemnorm[\ptinrsuccsubseqtwoiter]{\trslininrsuccsubseqtwoiter} - \Riemnorm[\ptinrsuccsubseqthriter]{\trslininrsuccsubseqthriter} \leq \abs*{\Riemnorm[\ptinrsuccsubseqtwoiter]{\trslininrsuccsubseqtwoiter} - \Riemnorm[\ptinrsuccsubseqthriter]{\trslininrsuccsubseqthriter}} \leq \epssix
            \end{align}}{
            \begin{align}
                2\epssix = 3\epssix - \epssix < \Riemnorm[\ptinrsuccsubseqtwoiter]{\trslininrsuccsubseqtwoiter} - \Riemnorm[\ptinrsuccsubseqthriter]{\trslininrsuccsubseqthriter} \leq \abs*{\Riemnorm[\ptinrsuccsubseqtwoiter]{\trslininrsuccsubseqtwoiter} - \Riemnorm[\ptinrsuccsubseqthriter]{\trslininrsuccsubseqthriter}} \leq \epssix
            \end{align}
            }
            holds for some $\succsubseqidxtwo\in\setNz$.
            This is, however, a contradiction.

            Now, we verify that, for any infinite subsequence of $\succset$, denoted by $\brc*{\ptinrsuccsubseqtwoiter}_{\succsubseqidxtwo}$, and any $\epssix > 0$, there exists $\succsubseqidxtwo\in\setNz$ such that $\Riemnorm[\ptinrsuccsubseqtwoiter]{\trslininrsuccsubseqtwoiter} < \epssix$.
            Consider a subsequence of $\succset$ that realizes $\limsup_{\inriteridx\to\infty}\Riemnorm[\ptinriter]{\trslininriter}$.
            For any $\succsubseqidxtwo\in\setNz$, there exists an index $\inrsuccsubseqtwoiter\in\succset$ of the subsequence that satisfies $\Riemnorm[\ptinrsuccsubseqtwoiter]{\trslininrsuccsubseqtwoiter} < \inv{\succsubseqidxtwo}$.
            Since the sequence $\brc*{\Riemnorm[\ptinrsuccsubseqtwoiter]{\trslininrsuccsubseqtwoiter}}_{\succsubseqidxtwo}$ converges to zero as $\succsubseqidxtwo\to\infty$, we see that the original subsequence also converges to zero, meaning that $\limsup_{\inriteridx\to\infty}\Riemnorm[\ptinriter]{\trslininriter}=0$. 
            Thus, combining this with \cref{theo:liminftrslininriter} yields $\lim_{\inriteridx\to\infty}\Riemnorm[\ptinriter]{\trslininriter} = 0$.
            The proof is complete.
        \end{proof}
        Using \cref{theo:limtrslininriter}, we derive the bound on $\Riemnorm[\ptinriter]{\gradstr[]\Lagfun\paren*{\allvarinriter}}$: for any $\ptinriter\in\setNz$,  
        \isextendedversion{
        \begin{align}
            \Riemnorm[\ptinriter]{\gradstr[]\Lagfun\paren*{\allvarinriter}} \leq \Riemnorm[\ptinriter]{\trslininriter} + \sum_{\ineqidx\in\ineqset} \abs*{\frac{\barrparam[]}{\ineqfun[\ineqidx]\paren*{\ptinriter}} - \ineqLagmultinriter[\ineqidx]}\Riemnorm[\ptinriter]{\gradstr\ineqfun[\ineqidx]\paren*{\ptinriter}}.
        \end{align}}{$\Riemnorm[\ptinriter]{\gradstr[]\Lagfun\paren*{\allvarinriter}} \leq \Riemnorm[\ptinriter]{\trslininriter} + \sum_{\ineqidx\in\ineqset} \abs*{\frac{\barrparam[]}{\ineqfun[\ineqidx]\paren*{\ptinriter}} - \ineqLagmultinriter[\ineqidx]}\Riemnorm[\ptinriter]{\gradstr\ineqfun[\ineqidx]\paren*{\ptinriter}}$.}
        By \cref{theo:limtrslininriter} and additionally supposing \cref{assu:ineqLagmultconverged}, the right-hand side converges to zero as $\inriteridx\to\infty$.
        We formally state the result in the following corollary.
        \begin{corollary}\label{coro:limgradLaginriter}
            Under \cref{assu:retrsecondord,assu:radialLCone,assu:radialLCtwo,assu:Cauchydecreasing,assu:ineqfunbounded,assu:objfunbounded,assu:pullbackineqfuncont,assu:ineqLagmultbounded,assu:retrdistbound,assu:ineqLagmultconverged}, 
            \isextendedversion{
            \begin{align}
                \lim_{\inriteridx\to\infty}\Riemnorm[\ptinriter]{\gradstr[]\Lagfun\paren*{\allvarinriter}} = 0.
            \end{align}
            }{
            $\lim_{\inriteridx\to\infty}\Riemnorm[\ptinriter]{\gradstr[]\Lagfun\paren*{\allvarinriter}} = 0$.}
        \end{corollary}
        
        Next, we provide the second-order analysis.
        Recall that $\mineigval\sbra*{\trsquadinriter}\in\setR[]$ is the minimum eigenvalue of $\trsquadinriter$.
        \begin{theorem}\label{theo:limsupsosp}
            Under \cref{assu:retrsecondord,assu:radialLCone,assu:radialLCtwo,assu:Cauchydecreasing,assu:ineqfunbounded,assu:pullbackineqfuncont,assu:ineqLagmultbounded,assu:ineqLagmultconverged}, \assuenumiCref{assu:objfunbounded}{assu:objfunboundedbelow}, and \Cref{assu:Eigendecreasing},
            \begin{equation}\label{ineq:limsupmineigvalinriter}
                \limsup_{\inriteridx\to\infty}\mineigval\sbra*{\trsquadinriter} \geq 0.
            \end{equation}
        \end{theorem}
    \begin{proof}
        \isextendedversion{
        See \cref{appx:proofinnerglobal}.}{We defer the complete proof to \cite[Appendix D.2]{Obaraetal2025APrimalDualIPTRMfor2ndOrdStnryPtofRiemIneqCstrOptimProbs} and here provide an overview.
        We argue by contradiction.
        Suppose that there exist $\constEigen > 0$ and an index $\tholdidxthr\in\setNz$ such that $\mineigval\sbra*{\trsquadinriter} < - \constEigen$ for all $\inriteridx\geq\tholdidxthr$.
        The remaining proof follows the same progression as the proof of \cref{theo:liminftrslininriter}. 
        Using \enumicref{lemm:diffpredared}{lemm:diffpredaredsmallo} and \cref{assu:Eigendecreasing}, we obtain the same inequality $\abs*{\trratioinriter - 1} \leq \frac{1}{2}$.
        Furthermore, under \cref{assu:Eigendecreasing}, instead of \cref{ineq:diffmeritfunCauchy}, we derive $\meritfun[{\barrparam[]}]\paren*{\ptinrsuccsubseqiter} - \meritfun[{\barrparam[]}]\paren*{\ptinrsuccsubseqiterp} \geq \trratiothold\constEigen\paren*{\trradiusinrsuccsubseqiter}^{2}\epseigen \geq 0$.
        Then, we reach a contradiction to the convergence of $\brc*{\trradiusinriter}_{\inriteridx}$.}
    \end{proof}

    Based on the analyses above, we prove the global convergence of \cref{algo:RIPTRMInner} in the following theorem.
    \isextendedversion{}{Recall that $\mineigval\sbra*{\trsquadinriter}$ denotes the minimum eigenvalue of $\trsquadinriter$.}
    \begin{theorem}\label{theo:innerglobconv}
        Let $\epsgradLag, \epscompl, \epssosp \in\setRpp[]$ be any positive scalars.
        Under \cref{assu:retrsecondord,assu:radialLCone,assu:radialLCtwo,assu:Cauchydecreasing,assu:ineqfunbounded,assu:objfunbounded,assu:pullbackineqfuncont,assu:ineqLagmultbounded,assu:retrdistbound,assu:ineqLagmultconverged}, the following hold after a finite number of iterations of \cref{algo:RIPTRMInner}: \isextendedversion{\begin{align}
            &\Riemnorm[\ptinriter]{\gradstr\Lagfun\paren*{\allvarinriter}} \leq \epsgradLag,\label{ineq:gradLaginriterconverged}\\
            &\norm*{\Ineqfunmat[]\paren*{\ptinriter}\ineqLagmultinriter[] - \barrparam[]\onevec} \leq \epscompl,\label{ineq:complinriterconverged}\\
            &\ineqLagmultinriter[] > 0, \ineqfun[]\paren*{\ptinriter} > 0.\label{ineq:ineqfunineqLagmultinriterfeasi}
        \end{align}}{
        \vspace{-1\baselineskip}
        \begin{align}
            \Riemnorm[\ptinriter]{\gradstr\Lagfun\paren*{\allvarinriter}} \leq \epsgradLag, \; \norm*{\Ineqfunmat[]\paren*{\ptinriter}\ineqLagmultinriter[] - \barrparam[]\onevec} \leq \epscompl, \; \ineqLagmultinriter[] > 0, \; \ineqfun[]\paren*{\ptinriter} > 0.\label{ineq:gradLagcomplprimaldualfeasiinriterconverged}
        \end{align}
        }
        Additionally, suppose \cref{assu:Eigendecreasing}.
        Then, \isextendedversion{\eqscref{ineq:gradLaginriterconverged,ineq:complinriterconverged,ineq:ineqfunineqLagmultinriterfeasi},}{\eqscref{ineq:gradLagcomplprimaldualfeasiinriterconverged}} and 
        \begin{align}\label{ineq:sospinriterconverged}
            \mineigval\sbra*{\trsquadinriter} \geq - \epssosp
        \end{align}
        hold after a finite number of iterations of \cref{algo:RIPTRMInner}.
    \end{theorem}
    \begin{proof}
        
    \isextendedversion{\Eqcref{ineq:gradLaginriterconverged}}{The first condition in \cref{ineq:gradLagcomplprimaldualfeasiinriterconverged}} follows from \cref{coro:limgradLaginriter}.
    For \isextendedversion{\cref{ineq:complinriterconverged}}{the second condition in \cref{ineq:gradLagcomplprimaldualfeasiinriterconverged}}, we have
    \isextendedversion{\begin{align}
        \norm*{\Ineqfunmat[]\paren*{\ptinriter}\ineqLagmultinriter[] - \barrparam[]\onevec} \leq \norm*{\Ineqfunmat[]\paren*{\ptinriter}}\norm*{\ineqLagmultinriter[] - \barrparam[]\inv{\Ineqfunmat[]\paren*{\ptinriter}}\onevec}.
    \end{align}}{
    $\norm*{\Ineqfunmat[]\paren*{\ptinriter}\ineqLagmultinriter[] - \barrparam[]\onevec} \leq \norm*{\Ineqfunmat[]\paren*{\ptinriter}}\norm*{\ineqLagmultinriter[] - \barrparam[]\inv{\Ineqfunmat[]\paren*{\ptinriter}}\onevec}$.
    }
    Under \assuenumicref{assu:ineqfunbounded}{assu:ineqfunboundedvalue} and  \cref{assu:ineqLagmultconverged}, the right-hand side tends to zero as $\inriteridx\to\infty$, which implies that \isextendedversion{\cref{ineq:complinriterconverged}}{the second condition in \cref{ineq:gradLagcomplprimaldualfeasiinriterconverged}} holds for all $\inriteridx\in\setNz$ sufficiently large.
    Notice that \isextendedversion{\cref{ineq:ineqfunineqLagmultinriterfeasi} holds}{the third and fourth conditions in  \cref{ineq:gradLagcomplprimaldualfeasiinriterconverged} hold} for every $\inriteridx\in\setNz$ by lines 15 and 20 of \cref{algo:RIPTRMInner}.
    \isextendedversion{Therefore, we conclude that \cref{ineq:gradLaginriterconverged,ineq:complinriterconverged,ineq:ineqfunineqLagmultinriterfeasi} hold under \cref{assu:retrsecondord,assu:radialLCone,assu:radialLCtwo,assu:Cauchydecreasing,assu:ineqfunbounded,assu:objfunbounded,assu:pullbackineqfuncont,assu:ineqLagmultbounded,assu:retrdistbound,assu:ineqLagmultconverged}.}{}
    \isextendedversion{Next, suppose \cref{assu:Eigendecreasing} additionally. Considering}{Next, considering}
    a subsequence that realizes \cref{ineq:limsupmineigvalinriter}, we guarantee \cref{ineq:sospinriterconverged} after finitely many iterations.

    \end{proof}
    

        \stepcounter{assumptionsection}

    \subsection{Global convergence of \texorpdfstring{\cref{algo:RIPTRMOuter}}{outer iterations}}\label{subsec:globconvouter}
        By \cref{theo:innerglobconv}, \cref{algo:RIPTRMInner} serves as an admissible inner iteration of \cref{algo:RIPTRMOuter} to compute the next iterate $\allvarotriterp$.
        In this subsection, we establish the global convergence of \cref{algo:RIPTRMOuter}. 
        We then show the convergence to an \AKKT{} point and an \SOSP{} when the generated sequence accumulates at some point.
        
        We first define a subset of the subsequences.
        This concept was originally introduced in the Euclidean setting by Conn et al.~\cite[Section~4.4]{Connetal2000PrimalDualTRAlgoforNonconvexNLP}.
        \begin{definition}
            Let $\brc*{\ptotrsubseqiter}_{\subseqidx}$ be a subsequence of \cref{algo:RIPTRMOuter}.
            We say that $\brc*{\ptotrsubseqiter}_{\subseqidx}$ is \asymptoticallyconsistent{} if, for each $\ineqidx\in\ineqset$, either
            \isextendedversion{
            \begin{align}
                \lim_{\subseqidx\to\infty}\ineqfun[\ineqidx]\paren*{\ptotrsubseqiter} = 0 \text{ or } \liminf_{\subseqidx\to\infty}\ineqfun[\ineqidx]\paren*{\ptotrsubseqiter} > 0
            \end{align}}{$\lim_{\subseqidx\to\infty}\ineqfun[\ineqidx]\paren*{\ptotrsubseqiter}$ $=0$ or \\ $\liminf_{\subseqidx\to\infty}\ineqfun[\ineqidx]\paren*{\ptotrsubseqiter}$ $>0$}
            holds.
            The former and the latter constraints are said to be asymptotically active and inactive, respectively. 
            For such $\brc*{\ptotrsubseqiter}_{\subseqidx}$, we define $\asymptotactsubset\coloneqq\brc*{\ineqidx\in\ineqset\colon\lim_{\subseqidx\to\infty}\ineqfun[\ineqidx]\paren*{\ptotrsubseqiter}=0}$.
        \end{definition}
        Then, we assume the following:        \begin{assumption}\label{assu:allinriterwellposed}
            For any $\otriteridx\in\setNz$, the point $\ptotriter\in\mani$ satisfies \stopcond{} in \cref{algo:RIPTRMOuter}.
        \end{assumption}
        \begin{assumption}\label{assu:gradineqfunotriterbounded}
            For an \asymptoticallyconsistent{} sequence $\brc*{\ptotrsubseqiter}_{\subseqidx}$, 
            the corresponding sequence
            $\brc*{\Riemnorm[\ptotrsubseqiter]{\gradstr\ineqfun[\ineqidx]\paren*{\ptotrsubseqiter}}}_{\subseqidx}$ is bounded for every $\ineqidx\in\ineqset$.
        \end{assumption}
        Note that, as we analyzed in \cref{subsec:globconvinner}, \cref{assu:allinriterwellposed} is fulfilled if \cref{assu:retrsecondord,assu:radialLCone,assu:radialLCtwo,assu:Cauchydecreasing,assu:ineqfunbounded,assu:objfunbounded,assu:pullbackineqfuncont,assu:ineqLagmultbounded,assu:retrdistbound,assu:ineqLagmultconverged,assu:Eigendecreasing} hold for every inner iteration.
        \cref{assu:gradineqfunotriterbounded} holds if $\brc*{\ptotrsubseqiter}_{\subseqidx}$ or $\brc*{\ptotriter}_{\otriteridx}$ is bounded; the latter
        is assumed in, e.g., 
        \cite[Assumption~(C2)]{LaiYoshise2024RiemIntPtMethforCstrOptimonMani}.

        Now, we prove the global convergence of \cref{algo:RIPTRMOuter} by examining the limiting behavior of $\brc*{\allvarotriter}_{\otriteridx}$.
        The result is a Riemannian extension of an \IPTRM{} by Conn et al.~\cite{Connetal2000PrimalDualTRAlgoforNonconvexNLP}.
        In the theorem, we consider a sequence of tangent vectors $\brc*{\tanvecone[\ptotrsubseqiter]}_{\subseqidx}$, noting that each $\tanvecone[\ptotrsubseqiter]$ lies in a different tangent space $\tanspc[\ptotrsubseqiter]\mani$.
        \begin{theorem}\label{theo:limglobconvotr}
            Suppose \cref{assu:allinriterwellposed}. 
            Then, the following hold:
            \isextendedversion{
            \begin{align}
                &\lim_{\otriteridx\to\infty}\Riemnorm[\ptotriter]{\gradstr \Lagfun\paren*{\allvarotriter}} = 0,\label{eq:limgradLagfunotriter}\\
                &\lim_{\otriteridx\to\infty}\norm*{\IneqLagmultmat[\otriteridx]\ineqfun[]\paren*{\ptotriter}} = 0,\label{eq:limcomplotriter}\\
                &\liminf_{\otriteridx\to\infty}\ineqfun[\ineqidx]\paren*{\ptotriter} \geq 0 \text{ and } \liminf_{\otriteridx\to\infty}\ineqLagmultotriter[\ineqidx] \geq 0 \text{ for all } \ineqidx \in \ineqset.\label{eq:liminffeasiotriter}
            \end{align}
            }{
            \begin{align}
                &\lim_{\otriteridx\to\infty}\Riemnorm[\ptotriter]{\gradstr \Lagfun\paren*{\allvarotriter}} = 0, \lim_{\otriteridx\to\infty}\norm*{\IneqLagmultmat[\otriteridx]\ineqfun[]\paren*{\ptotriter}} = 0,\label{eq:limgradLagfuncomplotriter}\\ &\liminf_{\otriteridx\to\infty}\ineqfun[\ineqidx]\paren*{\ptotriter} \geq 0 \text{ and } \liminf_{\otriteridx\to\infty}\ineqLagmultotriter[\ineqidx] \geq 0 \text{ for all } \ineqidx \in \ineqset.\label{eq:liminfprimaldualfeasiotriter}
            \end{align}
            }
            Let $\brc*{\ptotrsubseqiter}_{\subseqidx}$ be an \asymptoticallyconsistent{} sequence, and let $\brc*{\tanvecone[\ptotrsubseqiter]}_{\subseqidx}$ be any sequence of tangent vectors satisfying $\tanvecone[\ptotrsubseqiter]\in\tanspc[\ptotrsubseqiter]\mani$ for all $\subseqidx\in\setNz$,
            \isextendedversion{
            \begin{subequations}\label{eq:asymptoticwctanvec}  
                \begin{align}
                        &\brc*{\Riemnorm[\ptotrsubseqiter]{\tanvecone[\ptotrsubseqiter]}}_{\subseqidx} \text{ is bounded}, \label{eq:asymptoticwctanvecbounded}\\
                        &\metr[\ptotrsubseqiter]{\gradstr\ineqfun[\ineqidx]\paren*{\ptotrsubseqiter}}{\tanvecone[\ptotrsubseqiter]}=0 \text{ for all } \ineqidx\in\asymptotactsubset \text{ and all } \subseqidx\in\setNz.\label{eq:asymptoticwctanvecorthogonal}
                \end{align}
            \end{subequations}}{
            \begin{equation}\label{eq:asymptoticwctanvec}
                \brc*{\Riemnorm[\ptotrsubseqiter]{\tanvecone[\ptotrsubseqiter]}}_{\subseqidx} \text{ is bounded} \text{ and } \metr[\ptotrsubseqiter]{\gradstr\ineqfun[\ineqidx]\paren*{\ptotrsubseqiter}}{\tanvecone[\ptotrsubseqiter]}=0 \text{ for all } \ineqidx\in\asymptotactsubset.
            \end{equation}}
            Additionally, suppose that \cref{assu:gradineqfunotriterbounded} holds and that \seccondflag{} is \True{} in \cref{algo:RIPTRMOuter}.
            Then, it also follows that
            \isextendedversion{\begin{equation}\label{eq:liminfHessLagfunotriter}
                \liminf_{\subseqidx\to\infty} \metr[\ptotrsubseqiter]{\Hess\Lagfun\paren*{\allvarotrsubseqiter}\sbra*{\tanvecone[\ptotrsubseqiter]}}{\tanvecone[\ptotrsubseqiter]} \geq 0.
            \end{equation}}{
            \begin{equation}\label{eq:liminfHessLagfunotriter}
                \liminf_{\subseqidx\to\infty} \metr[\ptotrsubseqiter]{\Hess[\pt]\Lagfun\paren*{\allvarotrsubseqiter}\sbra*{\tanvecone[\ptotrsubseqiter]}}{\tanvecone[\ptotrsubseqiter]} \geq 0.
            \end{equation}}
        \end{theorem}
        \isextendedversion{\begin{proof}
            See \cref{appx:proofouterglobal}
        \end{proof}}{
        \begin{proof}
            We defer the full proof to \cite[Appendix~D.3]{Obaraetal2025APrimalDualIPTRMfor2ndOrdStnryPtofRiemIneqCstrOptimProbs} and provide here a brief overview.
            The first equality in \cref{eq:limgradLagfuncomplotriter} and the conditions~\cref{eq:liminfprimaldualfeasiotriter} directly follow from \cref{eq:stopcondKKT} and \cref{eq:stopcondstrictfeasi}, respectively.
            For the second equality in \cref{eq:limgradLagfuncomplotriter}, it follows from \cref{eq:stopcondbarrcompl} that
            $\norm*{\IneqLagmultmat[\otriteridx]\ineqfun[]\paren*{\ptotriter}} \leq \norm*{\IneqLagmultmatotriter \ineqfun[]\paren*{\ptotriter} - \barrparamotriterm \onevec} + \barrparamotriterm\norm*{\onevec}  \leq \forcingfuncompl\paren*{\barrparamotriterm} + \ineqdime\barrparamotriterm$,
            where the right-hand side tends to zero as $\otriteridx\to\infty$.
            Hence, the second equality in \cref{eq:limgradLagfuncomplotriter} holds.
            Next, we prove \cref{eq:liminfHessLagfunotriter} by contradiction.
            Assume that there exists a sequence of vectors $\brc*{\tanvecone[\ptotrsubseqiter]}_{\subseqidx}$ satisfying \cref{eq:asymptoticwctanvec} such that
            $\liminf_{\subseqidx\to\infty} \metr[\ptotrsubseqiter]{\Hess\Lagfun\paren*{\allvarotrsubseqiter}\sbra*{\tanvecone[\ptotrsubseqiter]}}{\tanvecone[\ptotrsubseqiter]} < 0$.
            By definition, we have
            $\metr[\ptotrsubseqiter]{\Hess\Lagfun\paren*{\allvarotrsubseqiter}\sbra*{\tanvecone[\ptotrsubseqiter]}}{\tanvecone[\ptotrsubseqiter]} = \metr[\ptotrsubseqiter]{\trsquad\paren*{\allvarotrsubseqiter}\sbra*{\tanvecone[\ptotrsubseqiter]}}{\tanvecone[\ptotrsubseqiter]} - \metr[\ptotrsubseqiter]{\ineqgradopr[\ptotrsubseqiter]\IneqLagmultmat[\otrsubseqiter]\inv{\Ineqfunmat[]\paren*{\ptotrsubseqiter}}\coineqgradopr[\ptotrsubseqiter]\sbra*{\tanvecone[\ptotrsubseqiter]}}{\tanvecone[\ptotrsubseqiter]}$.
            We show that the right-hand side tends to zero as $\subseqidx\to\infty$, which yields a contradiction.
            The proof is complete.
        \end{proof}
        }

        We proceed to consider the global convergence when 
        $\brc*{\ptotriter}_{\otriteridx}$ has an accumulation point. 
        Let $\ptaccum\in\feasirgn$ be any accumulation point, and let $\brc*{\ptotrsubseqiter}_{\subseqidx}$ be a subsequence of the outer iterates that realizes $\ptotrsubseqiter \to \ptaccum$ as $\subseqidx\to\infty$.
        Here, $\brc*{\ptotrsubseqiter}_{\subseqidx}$ is \asymptoticallyconsistent{} by the continuity of the constraint functions, and
        the set of asymptotically active constraints is determined by the values of the inequality constraints at $\ptaccum$, i.e., $\asymptotactsubset=\activeineqset\paren*{\ptaccum}$; see \cref{subsec:optimcond} for the definition of $\activeineqset\paren*{\ptaccum}$.
        In this case, we additionally assume the following conditions:
        \isextendedversion{}{recall that $\partxp[]{\paren*{\cdot}}{\paren*{\cdot}}$ is the parallel transport.}
        \begin{assumption}\label{assu:LICQptaccum}
            The point $\ptaccum\in\mani$ satisfies the \LICQ{}.
        \end{assumption}
        \begin{assumption}\label{assu:LipschitzHess}
            There exist $\constthr[\objfun] > 0$ and a closed neighborhood $\closedsubset_{\objfun}\subseteq\mani$ of $\ptaccum$ such that, for all $\pt[]\in\closedsubset_{\objfun}$,
            \isextendedversion{
            \begin{equation}\label{ineq:LipschitzHessobjfun}
                \opnorm{\Hess\objfun\paren*{\ptaccum} - \partxp[]{\ptaccum}{\pt[]} \circ \Hess\objfun\paren*{\pt[]} \circ \partxp[]{\pt[]}{\ptaccum}} \leq \constthr[\objfun]\Riemdist{\pt[]}{\ptaccum}.
            \end{equation}
            }{
            $\opnorm{\Hess\objfun\paren*{\ptaccum} - \partxp[]{\ptaccum}{\pt[]} \circ \Hess\objfun\paren*{\pt[]} \circ \partxp[]{\pt[]}{\ptaccum}} \leq \constthr[\objfun]\Riemdist{\pt[]}{\ptaccum}$.}
            In addition, for each $\ineqidx\in\ineqset$, there exist $\constthr[{\ineqfun[\ineqidx]}] > 0$ and a closed neighborhood $\closedsubset_{\ineqfun[\ineqidx]}\subseteq\mani$ of $\ptaccum$ such that, for all $\pt[]\in\closedsubset_{\ineqfun[\ineqidx]}$,
            \isextendedversion{
            \begin{equation}\label{ineq:LipschitzHessineqfun}
                \opnorm{\Hess\ineqfun[\ineqidx]\paren*{\ptaccum} - \partxp[]{\ptaccum}{\pt[]} \circ \Hess\ineqfun[\ineqidx]\paren*{\pt[]} \circ \partxp[]{\pt[]}{\ptaccum}} \leq \constthr[{\ineqfun[\ineqidx]}]\Riemdist{\pt[]}{\ptaccum}.
            \end{equation}
            }{
            $\opnorm{\Hess\ineqfun[\ineqidx]\paren*{\ptaccum} - \partxp[]{\ptaccum}{\pt[]} \circ \Hess\ineqfun[\ineqidx]\paren*{\pt[]} \circ \partxp[]{\pt[]}{\ptaccum}} \leq \constthr[{\ineqfun[\ineqidx]}]\Riemdist{\pt[]}{\ptaccum}$.}
        \end{assumption}
        \cref{assu:LICQptaccum} is standard in the literature; for example, \cite[Assumption (C1)]{LaiYoshise2024RiemIntPtMethforCstrOptimonMani} and \cite[Propositions~3.2, 3.4, 4.2]{LiuBoumal2020SimpleAlgoforOptimonRiemManiwithCstr}. 
        We also note that in the equality-constrained setting on $\mani=\setR[\dime]$, \LICQ{} implies that the feasible region is a smooth manifold in a neighborhood of $\ptaccum$; hence it falls within the scope of unconstrained Riemannian optimization~\cite[Section~7.7]{Boumal23IntroOptimSmthMani}.
        However, this does not hold in the presence of inequality constraints, which is the setting considered here.
        \isextendedversion{
        \cref{assu:LipschitzHess} holds if $\objfun$ and $\brc*{\ineqfun[\ineqidx]}_{\ineqidx\in\ineqset}$ are of class $C^{3}$ as in \cref{lemm:diffpartxpHesscont}.}{\cref{assu:LipschitzHess} holds if $\objfun$ and $\brc*{\ineqfun[\ineqidx]}_{\ineqidx\in\ineqset}$ are of class $C^{3}$~\cite[Lemma C.6]{Obaraetal2025APrimalDualIPTRMfor2ndOrdStnryPtofRiemIneqCstrOptimProbs}.}
        In the following theorem, we prove that $\ptaccum$ is an \AKKT{} point and a \wSOSP{} under assumptions.
        \isextendedversion{}{See \cref{subsec:optimcond} for the definition of $\weakcriticalcone$ in the proof.}
        \begin{theorem}\label{theo:globconvotriterptaccum}
            Let $\ptaccum\in\feasirgn$ be any accumulation point of a sequence generated by \cref{algo:RIPTRMOuter}.
            Under \cref{assu:allinriterwellposed}, $\ptaccum\in\feasirgn$ is an \AKKT{} point of \RICO{} \cref{prob:RICO}.
            Suppose additionally that \cref{assu:LICQptaccum,assu:LipschitzHess} hold 
            and \seccondflag{} is \True{} in \cref{algo:RIPTRMOuter}.
            Then, $\ptaccum$ is a \wSOSP{} of \RICO{}~\cref{prob:RICO} with an associated Lagrange multiplier vector $\ineqLagmultaccum[]$.
        \end{theorem}
        \isextendedversion{
        \begin{proof}
            See \cref{appx:proofouterglobal}.
        \end{proof}
        }{
        \begin{proof}
            We defer the complete proof to \cite[Appendix~D.3]{Obaraetal2025APrimalDualIPTRMfor2ndOrdStnryPtofRiemIneqCstrOptimProbs} and here provide an overview.
            By \cref{theo:limglobconvotr}, the sequence $\brc*{\ptotrsubseqiter}_{\subseqidx}$ and the associated sequence $\brc*{\ineqLagmultotrsubseqiter[]}_{\subseqidx}$ satisfy the \isextendedversion{\AKKT{} conditions~\cref{def:AKKTconditions}}{\AKKT{} conditions}.
            For the second-order stationarity, we first note that $\ptaccum$ is a \KKT{} point.
            It follows from the proof of \cite[Theorem~2]{YamakawaSato2022SeqOptimCondforNLOonRiemManiandGlobConvALM} that the sequence $\brc*{\ineqLagmultotrsubseqiter[]}$ is bounded.
            Hence, without loss of generality, by taking a subsequence if necessary, we may assume that $\brc*{\ineqLagmultotrsubseqiter[]}$ converges to the associated Lagrange multipliers $\ineqLagmultaccum[]\in\setR[\ineqdime]$.
            We derive $\metr[\ptaccum]{\Hess[{\pt[]}]\Lagfun\paren*{\allvaraccum}\sbra*{\tanvecone[\ptaccum]}}{\tanvecone[\ptaccum]} \geq 0$ for any $\tanvecone[\ptaccum]\in\weakcriticalcone\paren*{\ptaccum}$ using \cref{eq:liminfHessLagfunotriter} and the parallel transport, which is peculiar to the Riemannian setting.
        \end{proof}
        }
        By \cref{prop:wSOSPequivSOSPudrSC}, we obtain the following corollary. 
        
        \begin{corollary}\label{coro:convSOSPptaccum}
            Suppose that the SC holds at $\ptaccum$ and \seccondflag{} is \True{} in \cref{algo:RIPTRMOuter}.
            Under \cref{assu:allinriterwellposed,assu:gradineqfunotriterbounded,assu:LICQptaccum,assu:LipschitzHess}, $\paren*{\ptaccum, \ineqLagmultaccum[]}$ is an \SOSP{} of \RICO{}~\cref{prob:RICO}.
        \end{corollary}
        Note that \SC{} is a standard assumption for establishing second-order global convergence~\cite{Connetal2000PrimalDualTRAlgoforNonconvexNLP,YabeYamashita2025ConvtoASecondOderCritPtbyATrustRegionSQPMethwithANonsmthMeritFun,Yamashita2022Yamashita2022ConvtoASecOdrCritPtbyaPDIntPtTrustRegionMethforNonlinSemiDefProgram,ArahataOkunoTakeda2023ComplexAnalofIPMsfor2ndOrdSntryPtsofNonlinSemiDefiOptimProbs}.
        In our experiments, we numerically verified that all reported solutions satisfied \SC{}; see \cref{subsec:resultdiscussion} for details.


    \section{Numerical experiments}\label{sec:experiment}
        In this section, we present numerical experiments on stable linear system identification, demonstrating the efficiency of \RIPTRM{}.
        For comparison, we also solve these problems using other Riemannian algorithms.
        All the experiments are implemented in Python with Pymanopt~2.2.0~\cite{TownsendKoepWeichwald2017PymanoptPyToolboxforOptimonManiusingAutoDiff} and executed on a MacBook Pro with an Apple M1 Max chip and 64GB of memory.
        Our code is available on GitHub~\cite{Obara2025CodeRiemIntPtTrustRegionMeth}.

    \subsection{Problem settings: stable linear system identification}\label{subsec:problemsetting}
        Linear system identification is the problem of estimating a linear continuous time-invariant autonomous system $\dot{\state[]}\paren*{t}=A\state[]\paren*{t}$ from observed data, where $A\in\setR[\sysdime\times\sysdime]$ denotes the system matrix.
        We identify a system while enforcing asymptotic stability and incorporating prior knowledge.
        Asymptotic stability is a fundamental property of many systems: the system is said to be asymptotically stable if all eigenvalues of $A$ have strictly negative real parts.
        Identification should also incorporate prior knowledge specific to the system, such as nonnegativity or partial information on certain entries, as discussed in \cite[Chapter~16]{Ljung99SysIdentifBook}.
        Such identification problems arise, for example, in mechanics, electrical engineering, and structural modeling.
        Following \cite{Obaraetal2024StabLinSysIdentifwithPriorKnwlbyRiemSQO}, we model stability using a Riemannian formulation and encode prior knowledge via nonlinear constraints:
        define $\Skew\paren*{\sysdime}\coloneqq\brc*{\matone\in\setR[\sysdime\times\sysdime]\relmiddle{\colon}\matone=-\trsp{\matone}}$ and $\Sympos\paren*{\sysdime}\coloneqq\brc*{\matone\in\setR[\sysdime\times\sysdime]\relmiddle{\colon}\matone\succ0}$.
        Then, the set $\mani_{\text{stab}}\coloneqq\Skew\paren*{\sysdime}\times\Sympos\paren*{\sysdime}\times\Sympos\paren*{\sysdime}$ denotes the product manifold that parameterizes the set of asymptotically stable matrices via $\A=(\J-\R)\Q$ with $(\J,\R,\Q)\in\mani_{\text{stab}}$~\cite{Prajnaetal02LMItoStabLinPortCtrlHamilSys}.
        Given noisy observations $\brc*{\state}_{\stateidx=1}^{\obsnum}\subseteq\setR[\sysdime]$ from the true stable system $\Atrue\in\setR[\sysdime\times\sysdime]$, we estimate a stable system matrix that fits the observed data.
        In the experiments, we consider the case with box-type constraints as prior knowledge.
        Formally, we address the following optimization problem:
        \isextendedversion{
        \begin{mini!}[2]
            {\paren*{\J, \R, \Q} \in \mani_{\text{stab}}}{\frac{1}{\obsnum\norm*{\state[0]}}\sum_{\stateidx=1}^{\obsnum-1}\fronorm*{\state[\stateidx+1] - \paren*{\eyemat + \sampintvl\paren*{\J - \R}\Q}\state[\stateidx]}\label{objfun:stablinsysidentif}}
            {\label{Prob:StabLinSysIdentif}}{}
            \addConstraint{\leftconst \leq \trsp{\stdbasis[\matrowidx]} \paren*{\J-\R}\Q \stdbasis[\matcolidx] \leq \rightconst, \quad\paren*{\matrowidx,\matcolidx} \in \ineqoneboxset\cup\ineqtwoboxset}{\label{ineq:oneboxcstr}}
            \addConstraint{\lengthconst^2 \leq \paren*{\trsp{\stdbasis[\matrowidx]} \paren*{\J-\R}\Q \stdbasis[\matcolidx] - \centerconst}^2, \quad\paren*{\matrowidx,\matcolidx} \in \ineqtwoboxset,}{\label{ineq:twoboxcstr}}
        \end{mini!}}{\begin{mini!}[2]
            {\paren*{\J, \R, \Q} \in \mani_{\text{stab}}}{\frac{1}{\obsnum\norm*{\state[0]}}\sum_{\stateidx=1}^{\obsnum-1}\fronorm*{\state[\stateidx+1] - \paren*{\eyemat + \sampintvl\paren*{\J - \R}\Q}\state[\stateidx]}\label{objfun:stablinsysidentif}}
            {\label{Prob:StabLinSysIdentif}}{}
            \addConstraint{\leftconst \leq \trsp{\stdbasis[\matrowidx]} \paren*{\J-\R}\Q \stdbasis[\matcolidx] \leq \rightconst, \quad\paren*{\matrowidx,\matcolidx} \in \ineqoneboxset\cup\ineqtwoboxset}{\label{ineq:oneboxcstr}}
            \addConstraint{\lengthconst^2 \leq \paren*{\trsp{\stdbasis[\matrowidx]} \paren*{\J-\R}\Q \stdbasis[\matcolidx] - \centerconst}^2, \quad\paren*{\matrowidx,\matcolidx} \in \ineqtwoboxset,}{\label{ineq:twoboxcstr}}
        \end{mini!}}
        where $\sampintvl > 0$ is the sampling interval, 
        $\ineqoneboxset,\ineqtwoboxset\subseteq\brc*{1,\ldots,\sysdime}\times\brc*{1,\ldots,\sysdime}$ are disjoint index sets, $\stdbasis[\matrowidx]\in\setR[\sysdime]$ is the $\matrowidx$-th standard basis, and $\leftconst, \rightconst, \lengthconst, \centerconst \in\setRpp[]$ satisfy $\leftconst \leq \Atrue_{\matrowidx\matcolidx} \leq \rightconst$ for $\paren*{\matrowidx, \matcolidx}\in\ineqoneboxset$ and $\Atrue_{\matrowidx\matcolidx}\in\sbra*{\leftconst, \centerconst - \lengthconst}\cup\sbra*{\centerconst + \lengthconst, \rightconst}$ for $\paren*{\matrowidx, \matcolidx}\in\ineqtwoboxset$.
    
        \textbf{Input.} 
        Our implementation follows \cite{Obaraetal2024StabLinSysIdentifwithPriorKnwlbyRiemSQO} with modifications.
        We consider the case $\sysdime = 5$, $\sampintvl=0.02$, 
        $\abs*{\ineqoneboxset}= \floor{0.2 \sysdime^{2}}$, and $\abs*{\ineqtwoboxset}=\floor{0.1\sysdime^{2}}$. 
        We draw $\paren*{\Jtrue,\Rtrue,\Qtrue}\in\mani$ at random using the Pymanopt's \texttt{Manifold.random\_point()}, and set the true system to $\Atrue \coloneqq \paren*{\Jtrue - \Rtrue} \Qtrue \in \setR[\dime\times\dime]$. 
        We generate $5$ independent observation sequences $\brc*{\state[\stateidx]}_{\stateidx=0}^{19}$ and use them in the objective function; thus, the total number of observations is $\obsnum=100$ in this experiment.
        For the initial points, we use interior points obtained numerically in advance: we solve 
        \isextendedversion{
        $\minimize_{\paren*{\J, \R, \Q}\ \in\mani_{\text{stab}}} 0 \ \subjectto $ \cref{ineq:oneboxcstr} and \cref{ineq:twoboxcstr}
        }{
        $\minimize_{\paren*{\J, \R, \Q}\ \in\mani_{\text{stab}}} 0$ subject to the original constraints}
        using \RALM{}~\cite[Algorithm~5]{LiuBoumal2020SimpleAlgoforOptimonRiemManiwithCstr}, until obtaining a solution with a residual of $10^{-2}$ or until reaching $4$ iterations.
        Finally, we fix one problem instance and draw $20$ initial interior points for that instance at random.      
    
    \subsection{Experimental environment}\label{subsec:experimentenv}
        We compare the following algorithms:
        \begin{itemize}
            \item \textbf{\RIPTRM{} (\tCG{})}: Riemannian interior point trust region method (\cref{algo:RIPTRMOuter,algo:RIPTRMInner}) with the search direction obtained by the truncated conjugate gradient (\tCG{}) method.
            \item \textbf{\RIPTRM{} (\Exact{})}: Riemannian interior point trust region method (\cref{algo:RIPTRMOuter,algo:RIPTRMInner}) with the \exactsolution{}.
            \item \RIPM{}: Riemannian interior point method~\cite[Algorithm~5]{LaiYoshise2024RiemIntPtMethforCstrOptimonMani}.
            \item \RSQO{}: Riemannian sequential quadratic optimization~~\cite{ObaraOkunoTakeda2021SQOforNLOonRiemMani}.
            \item \RALM{}: Riemannian augmented Lagrangian method~\cite[Algorithm~5]{LiuBoumal2020SimpleAlgoforOptimonRiemManiwithCstr}.
        \end{itemize}
        Here, we use the globally convergent \RIPM{}~\cite[Algorithm 5]{LaiYoshise2024RiemIntPtMethforCstrOptimonMani} rather than locally convergent one~\cite[Algorithm 2]{LaiYoshise2024RiemIntPtMethforCstrOptimonMani} for comparison, since the algorithms are initialized at points that may lie far from the solution set.
        In addition to our \RIPTRM{}s, we implemented \RIPM{}, \RSQO{}, and \RALM{} in Python following their original MATLAB implementations.
        We will explain the details of the search direction by \tCG{} and the \Exact{} steps in \RIPTRM{} at the end of this subsection.
        To measure the deviation of an iterate from the set of \KKT{} points, we introduce the residual defined as
        \isextendedversion{
        \begin{align}
            &\residual\paren*{\allvar}\\
            &\coloneqq\sqrt{\Riemnorm[\pt]{\gradstr[\pt]\Lagfun\paren*{\allvar}}^{2} + \sum_{\ineqidx\in\ineqset}\paren*{\min\paren*{0, \ineqLagmult[\ineqidx]}^{2} + \min\paren*{0, \ineqfun[\ineqidx]\paren*{\pt}}^{2} + \paren*{\ineqLagmult[\ineqidx]\ineqfun[\ineqidx]\paren*{\pt}}^{2} } + \Manvio\paren*{\pt}^{2}},
        \end{align}}{
        \begin{align}
            &\residual\paren*{\allvar}\\[-5mm]
            &\coloneqq\sqrt{\Riemnorm[\pt]{\gradstr[\pt]\Lagfun\paren*{\allvar}}^{2} + \sum_{\ineqidx\in\ineqset}\paren*{\min\paren*{0, \ineqLagmult[\ineqidx]}^{2} + \min\paren*{0, \ineqfun[\ineqidx]\paren*{\pt}}^{2} + \paren*{\ineqLagmult[\ineqidx]\ineqfun[\ineqidx]\paren*{\pt}}^{2} } + \Manvio\paren*{\pt}^{2}},
        \end{align}}
        where the first two terms are from the \isextendedversion{\KKT{}  conditions~\cref{eq:KKTconditions}}{\KKT{}  conditions}, and $\Manvio\colon\mani\to\setR[]$ measures the manifold constraint violations, that is, $\Manvio\paren*{\J, \R, \Q} = \fronorm{\J+\trsp{\J}} + \fronorm{\R-\trsp{\R}} + \fronorm{\Q - \trsp{Q}} + \posdefmanvio\paren*{R} + \posdefmanvio\paren*{Q}$ for $\mani_{\text{stab}}$, where $\posdefmanvio\paren*{\matone} \coloneqq +\infty$ if the matrix $\matone$ has a negative real eigenvalue and otherwise $0$.
        In practice, these violations are negligible in Riemannian optimization due to the use of the retractions.
        Each algorithm is run for 240 seconds.

        For \RIPTRM{}, we set the parameters to $\inittrradius_{0} = \frac{\typicaldist}{8}, \mininittrradius=10^{-15}, \maxtrradius=10, \forcingfungradLag\paren*{\barrparam[]}=\barrparam[],  \forcingfuncompl\paren*{\barrparam[]}=10^{-3}\,\barrparam[], \forcingfunsosp\paren*{\barrparam[]}=\barrparam[], \trratiothold=0.1$, and $\contcoeff=0.25$.
        Here, $\typicaldist>0$ denotes the Pymanopt's parameter \texttt{typical\_dist}, which specifies a scale for each manifold in the library.
        We apply the update rule for the barrier parameter in \cite{Obaraetal2025LocalConvofRiemIPMs} to have local convergence property: $\barrparam[0]=0.1$ and $\barrparam[\otriteridx+1] \leftarrow 0.5  \barrparam[\otriteridx]^{1.01}$. 
        \isextendedversion{
        We employ the clipping in the update of the dual variables with $\constfou=0.5$ and $\constfiv=10^{20}$; see \cref{subsec:computdualvar} for details.
        }{For the update of the dual variables, we use the clipping \cref{def:clippingupdate} with $\constfou=0.5$ and $\constfiv=10^{20}$.}
        For \RIPM{}, we use its default parameters.        
        For \RSQO{}, we use $\delta\in\{10^{-4},10^{-2}\}$, where the parameter $\delta$ controls the strength of the regularization added to the Hessian of the Lagrangian to enforce positive definiteness.
        Tuning $\delta$ affects the performance of \RSQO{}; in our experiments, these settings yielded better empirical performance.
        For \RALM{}, we use the default settings except that we employ the gradient descent as the subsolver of \RALM{} instead of limited-memory BFGS, the subsolver used in the original MATLAB implementation, as limited-memory BFGS is not implemented in Pymanopt~2.2.0.
    
    \subsubsection*{Implementation details of \RIPTRM{}}
        In this paper, our \RIPTRM{} employs the two search directions.
        The first is obtained by the \tCG{} method~\cite[Algorithm 6.4]{Boumal23IntroOptimSmthMani}.
        The quality of the search direction is guaranteed to be better than that of the \Cauchypoint{}; that is, the direction obtained by \tCG{} satisfies \cref{assu:Cauchydecreasing} with $\constCauchy = \frac{1}{2}$~\cite[Exercise 6.26]{Boumal23IntroOptimSmthMani}.
        Nevertheless, this direction may not satisfy \cref{assu:Eigendecreasing}~\cite[Exercise 6.27]{Boumal23IntroOptimSmthMani}.
        We use the \tCG{} method implementation available in the Riemannian trust region method for unconstrained optimization in Pymanopt with its default setting.
        We also set \seccondflag{}=\False{} for this setting.
        On the other hand, computing the \exactsolution{}, the global optimum of \cref{prob:TRS}, is challenging since the trust region subproblem~\cref{prob:TRS} is generally nonconvex.
        To compute the \exactsolution{}, we employ the algorithm by Adachi et al.~\cite{Adachietal2017SolvingTRSbyGenEigenProb}, which is based on the generalized eigenvalue problem.
        This direction satisfies \cref{assu:Cauchydecreasing,assu:Eigendecreasing} with $\constCauchy = \frac{1}{2}$ and $\constEigen = \frac{1}{2}$.
        In addition, the use of the \exactsolution{} is required for the local convergence, because it eventually coincides with the Newton step near a solution~\cite[Section~5.1]{Obaraetal2025LocalConvofRiemIPMs}.
        We implemented a Python version of the algorithm since the original code is written in MATLAB.\footnote{\href{https://people.maths.ox.ac.uk/nakatsukasa/codes.htm}{https://people.maths.ox.ac.uk/nakatsukasa/codes.htm}}
        We also set \seccondflag{}=\True{} for this setting.

    \subsection{Results and discussion}\label{subsec:resultdiscussion}
        We applied the algorithms to the problem, starting from 20 initial points.
        \begin{figure}[t]
            \centering
            \begin{minipage}{0.49\linewidth}
                \centering
                \includegraphics[width=\linewidth]{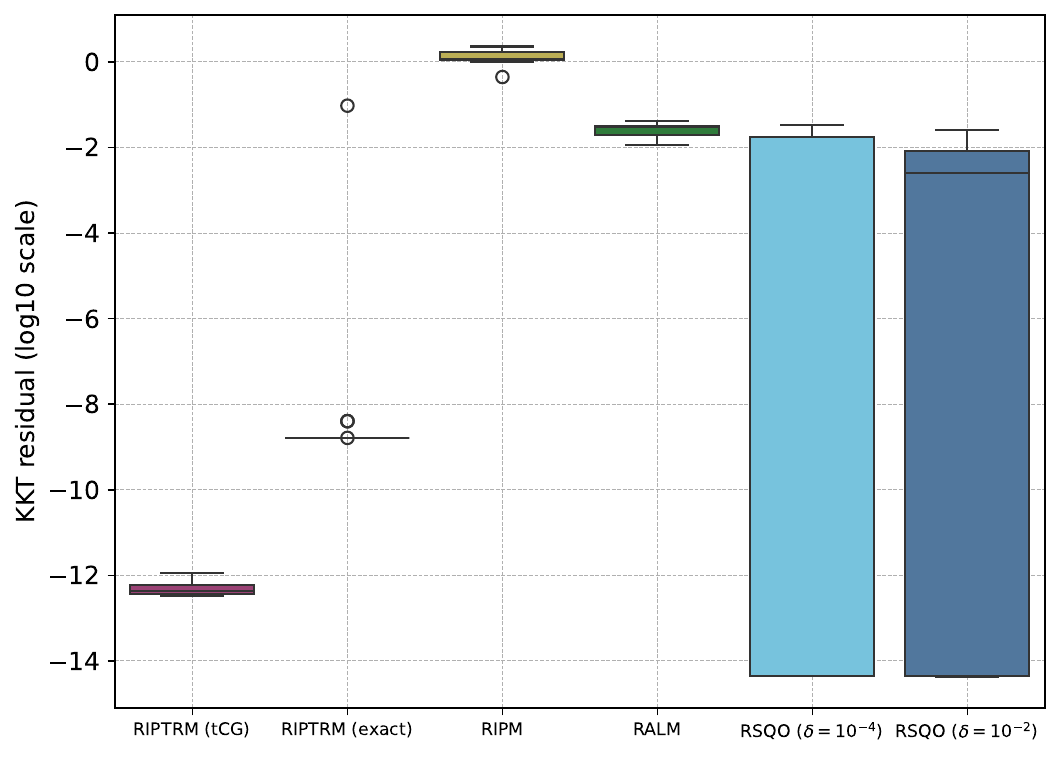}
                \caption{Box plots of minimum residuals among 20 instances 
                }
                \label{fig:boxplotstablinsysidentif}
            \end{minipage}
            \begin{minipage}{0.49\linewidth}
                \centering
                \includegraphics[width=\linewidth]{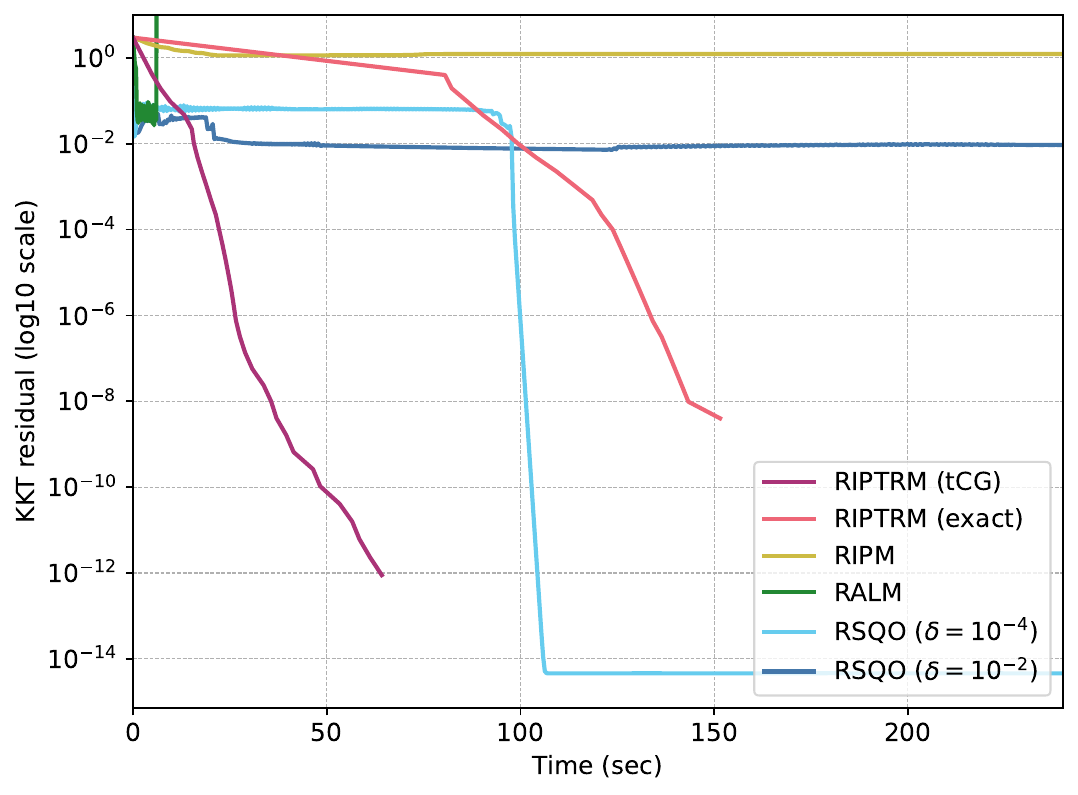}
                \caption{Residuals over time on a representative case 
                }
                \label{fig:residtimestablinsysidentif}
            \end{minipage}
        \end{figure}
        \cref{fig:boxplotstablinsysidentif} shows box plots of the minimum residuals computed by the solvers across 20 outputs.
        \cref{fig:residtimestablinsysidentif} illustrates a representative example of a residual over time for the algorithms.
        Here, we omit the inner iterations of \RIPTRM{}s and only plot the outer iterations for clarity.
        We also verified in our experiments that all reported solutions satisfy the strict complementarity condition; 
        specifically, for every output $(\ptaccum,\ineqLagmultaccum[])$ and every $\ineqidx\in\ineqset$, we never observed $|\ineqLagmultaccum[\ineqidx]|\leq 10^{-8}$ and $\abs*{\ineqfun[\ineqidx]\paren*{\ptaccum}}\leq 10^{-8}$ simultaneously.

        \RIPTRM{}s consistently achieved high-accuracy solutions across instances.
        In \cref{fig:boxplotstablinsysidentif}, the median residuals for \RIPTRM{} (\tCG{}) and \RIPTRM{} (\Exact{}) were approximately $10^{-12.368}$ and $10^{-8.787}$, respectively.
        While \RSQO{}s attained precise solutions to a certain extent, as indicated by its first quartile being below $10^{-14}$, its performance was quite sensitive to the choice of $\delta$.
        For $\delta=10^{-4}$, the median residual was approximately $ 10^{-14.346}$, whereas it was approximately $ 10^{-2.595}$ for $\delta=10^{-2}$, highlighting substantial sensitivity to the choice of $\delta$.
        Moreover, the third quartile is around $10^{-2}$ for both settings.
        This accuracy is still comparable to that of \RALM{} (see \cref{fig:boxplotstablinsysidentif}).
        Nevertheless, these results suggest that, in these experiments,
        \RSQO{} was more sensitive across problem instances than \RIPTRM{}s.
        For \RALM{}, the residual initially decreased but later increased substantially, reaching $61.648$, in \cref{fig:residtimestablinsysidentif}.
        This subsequent increase appears to have been driven by the divergence of the penalty coefficient in \RALM{}.
        \RIPM{} yielded only a slight decrease in the residual and did not reach high-accuracy solutions.
        As the iterations of \RIPM{} progressed, the norm of the search directions increased substantially, and the line-search procedure selected very small step sizes to ensure decrease of the merit function;
        for example, in \cref{fig:residtimestablinsysidentif}, the residual decreased from the initial value of $2.971$ to $1.130$.
        The norm of the search directions increased to around $10^{8}$, while the step size decreased to around $10^{-16}$.
        These phenomena may have slowed the progress of the iterates.
        Overall, in our experiments, the trust-region-based \RIPTRM{} produced more accurate solutions than the line-search-based \RIPM{}.
        
        \subsection{Auxiliary experiment}
        To support the efficiency of our algorithm, we conducted an additional study.
        \isextendedversion{In what follows, we give a summary of some auxiliary numerical results. The full details are provided in \cref{sec:additionalexperiment}.}{In what follows, we give a summary of some auxiliary numerical results.\footnote{The full details are provided in the extended version of this paper~\cite[Appendix F]{Obaraetal2025APrimalDualIPTRMfor2ndOrdStnryPtofRiemIneqCstrOptimProbs}.}}
        We consider a synthetic optimization problem: the minimization of the Rosenbrock function over the Grassmann manifold with simple inequality constraints, where the Hessian of the Lagrangian has a large negative eigenvalue at the initial point.
        We applied the aforementioned algorithms and measured both the second-order stationarity and the \KKT{} residual.

        The results show that \RIPTRM{} (\Exact{}) rapidly reduced the \KKT{} residual.
        Additionally, it achieved and maintained nonnegative second-order stationarity, indicating the successful computation of an \SOSP{}, whereas the other methods did not.
        This suggests that \exactsolution{}s are particularly robust and efficient for problem instances where the Hessian of the Lagrangian exhibits large negative eigenvalues, as they effectively incorporate information about these eigenvalues.
        

    \section{Concluding remarks}\label{sec:conclusion}
        In this paper, we proposed \RIPTRM{}, composed of \cref{algo:RIPTRMOuter,algo:RIPTRMInner}, for solving \RICO{}~\cref{prob:RICO}.
        For \cref{algo:RIPTRMInner}, we proved the global convergence to a point satisfying the stopping conditions in  \cref{theo:innerglobconv}. 
        We analyzed the limiting behavior of \cref{algo:RIPTRMOuter} in \cref{theo:limglobconvotr}, and established the global convergence of \cref{algo:RIPTRMOuter} to an \AKKT{} point and an \SOSP{} in \cref{theo:globconvotriterptaccum,coro:convSOSPptaccum}, respectively.
        In \cref{sec:experiment}, we presented numerical experiments.
        The results indicate that 
        \RIPTRM{}s consistently find solutions with high accuracy.
        \isextendedversion{
        We also introduced an eigenvalue-based subsolver for \RIPTRM{} to obtain the exact search direction by solving the trust region subproblems. 
        We observed that its performance is promising in an instance where the Hessian of the Lagrangian has a large negative eigenvalue.
        }{}

        In closing, we discuss future directions for further development of \RIPTRM{}:
        \begin{enumerate}
            \item \textbf{Treatment of equality constraints and slack variables:} 
            In Euclidean optimization, several \IPTRM{}s have been designed to manage both inequality and equality constraints although their guarantees of the global convergence guarantees are first-order~\cite{ByrdGilbertNocedal2000TRMethBsdIntPtTechniqueforNLP,ByrdLiuNocedal1997LocalBehavIntPtMethforNLP,YamashitaYabeTanabe2005GlobalSuperlinConvPrimalDualIntPtTRNethforLgeScaleCstrOptim}.
            Furthermore, these methods incorporate slack variables, which allow the algorithms to start from infeasible initial points.
            Integrating such techniques into our method could broaden its applicability.
            
            \item \textbf{Use of filter or funnel method:} 
            Our \RIPTRM{} currently uses the log barrier function as a merit function to ensure the global convergence.
            In the Euclidean optimization, filter~\cite{UlbrichUlbrichVicente2004GlobalConvPrimalDualIntPtFltrMethforNLP,Silvaetal2008GlobalConvPrimalDualIntPtFltrMethforNLPNewFltrOptimMeasComputRslt} and funnel methods~\cite{Curtisetal2017IntPtTrustFunnelAlgoforNLP} are also employed as alternative globalization strategies, demonstrating improved practical performance~\cite{BensonVanderbeiShanno2002IntPtMethforNonconvexNLPFltrMethMeritFun,KiesslingLeyfferVanaret2024UnfFuRestorSQPAlgo}.
            Extending these strategies to the Riemannian setting would enhance the efficiency and robustness of our method.

            \item \textbf{Complexity analysis:} 
            Non-asymptotic complexity analysis for algorithms in constrained Euclidean optimization is an active area of research~\cite{CartisGouldToint2022EvalComplexAlgoforNonconvOptim}.
            However, complexity analyses for \IPM{}s are still limited to the cases where the barrier parameter is fixed due to the non-Lipschitzness of the log barrier function~\cite{HinderYe2023WorstCaseIterBdforLogBarrMethProbwithNonconvexCstr,ArahataOkunoTakeda2023ComplexAnalofIPMsfor2ndOrdSntryPtsofNonlinSemiDefiOptimProbs}.
            The overall complexity of \IPM{}s for nonlinear optimization problems remains an open question even in the Euclidean context.
            It would be an significant contribution to provide the overall complexity for \IPM{}s and extend it to the Riemannian case.
                
        \end{enumerate}


    \section*{Acknowledgments}
        The authors are grateful to the associate editor and the two anonymous referees for their insightful comments and constructive suggestions.


    \appendix
    \setcounter{section}{0}
    \renewcommand{\thesection}{\Alph{section}}
    \setcounter{equation}{0}
    \renewcommand{\theequation}{\thesection-\arabic{equation}}
    \bibliographystyle{siamplain}
    \bibliography{reference.bib}

    \section{Proofs of \texorpdfstring{\cref{lemm:diffpredaredfixedpt,lemm:diffpredared}}{lemmas estimating the gaps between the reductions}}\label{appx:proofdiffpredared}
        In this section, we derive the proofs of the key lemmas, \cref{lemm:diffpredaredfixedpt,lemm:diffpredared}, that estimate the gaps between the predicted and actual reductions.
        \isextendedversion{}{See \cref{subsec:inneriter} for the definitions of $\pred$, $\ared$, $\meritfun[{\barrparam[]}]$, $\modelfun[{\allvar, \barrparam[]}]$, $\trsquad$, and $\trslin[{\barrparam[]}]$.
        Recall also that $\pullback[\pt]{\cdot}$ denotes the restriction of the pullback.
        }

        To prove the lemmas, we first express the gaps in the following form: under $\tanspc[\tmethr\dirpt]\paren*{\tanspc[\pt]\mani}\simeq\tanspc[\pt]\mani$, by \cref{lemm:ineqfunfeasi}, we have that, for all $\dirpt\in\tanspc[\pt]\mani$ sufficiently small,
        \isextendedversion{
        \begin{align}
            \begin{split}\label{ineq:diffpredared}
                &\abs*{\pred[{\allvar, \barrparam[]}]\paren*{\dirpt}-\ared[{\barrparam[]}]\paren*{\dirpt}}\\
                &= \abs*{\pullback[\pt]{\meritfun[{\barrparam[]}]}\paren*{\zerovec[\pt]} - \pullback[\pt]{\meritfun[{\barrparam[]}]}\paren*{\dirpt} - \paren*{\modelfun[{\allvar, \barrparam[]}]\paren*{\zerovec[\pt]} - \modelfun[{\allvar, \barrparam[]}]\paren*{\dirpt}}}\\
                &=\abs*{-\D\pullback[\pt]{\meritfun[{\barrparam[]}]}\paren*{\zerovec[\pt]}\sbra*{\dirpt} - \int_{0}^{1} \paren*{1-\tmethr} \D[2]\pullback[\pt]{\meritfun[{\barrparam[]}]}\paren*{\zerovec[\pt] + \tmethr\dirpt}\sbra*{\dirpt, \dirpt} \dd{\tmethr} \\
                &\quad + \frac{1}{2}\metr[\pt]{\trsquad\paren*{\allvar}\sbra*{\dirpt}}{\dirpt} + \metr[\pt]{\trslin[{\barrparam[]}]\paren*{\pt}}{\dir[\pt]}}\\
                &=\abs*{\int_{0}^{1}\paren*{1-\tmethr} \paren*{ \metr[\pt]{\trsquad\paren*{\allvar}\sbra*{\dirpt}}{\dirpt} -  \D[2]\pullback[\pt]{\meritfun[{\barrparam[]}]}\paren*{\tmethr\dirpt}\sbra*{\dirpt, \dirpt}} \dd{\tmethr}}\\
                &\leq \int_{0}^{1} \paren*{1-\tmethr}\Bigg(\abs*{\metr[\pt]{\paren*{\Hess\objfun\paren*{\pt} - \Hess\pullback[\pt]{\objfun}\paren*{\tmethr\dirpt}}\sbra*{\dirpt}}{\dirpt}} \\
                &\quad + \sum_{\ineqidx\in\ineqset}\abs*{\metr[\pt]{\paren*{-\ineqLagmult[\ineqidx]\Hess\ineqfun[\ineqidx]\paren*{\pt} + \frac{\barrparam[]}{\pullback[\pt]{{\ineqfun[\ineqidx]}}\paren*{\tmethr\dirpt}}\Hess\pullback[\pt]{\ineqfun[\ineqidx]}\paren*{\tmethr\dirpt}}\sbra*{\dirpt}}{\dirpt}} \\
                &\quad + \sum_{\ineqidx\in\ineqset}\abs*{\frac{\ineqLagmult[\ineqidx]}{\ineqfun[\ineqidx]\paren*{\pt}}\metr[\pt]{\gradstr\ineqfun[\ineqidx]\paren*{\pt}}{\dirpt}^{2} - \frac{\barrparam[]}{{\pullback[\pt]{\ineqfun[\ineqidx]}\paren*{\tmethr\dirpt}}^{2}} \metr[\pt]{\gradstr\pullback[\pt]{\ineqfun[\ineqidx]}\paren*{\tmethr\dirpt}}{\dirpt}^{2}}\Bigg) \dd{\tmethr},
            \end{split}
        \end{align}}{
        \vspace{-1\baselineskip}
        \begin{align}
            &\label{ineq:diffpredared}\\
            &\abs*{\pred[{\allvar, \barrparam[]}]\paren*{\dirpt}-\ared[{\barrparam[]}]\paren*{\dirpt}} = \abs*{\pullback[\pt]{\meritfun[{\barrparam[]}]}\paren*{\zerovec[\pt]} - \pullback[\pt]{\meritfun[{\barrparam[]}]}\paren*{\dirpt} - \paren*{\modelfun[{\allvar, \barrparam[]}]\paren*{\zerovec[\pt]} - \modelfun[{\allvar, \barrparam[]}]\paren*{\dirpt}}}\\
            &=\abs*{-\D\pullback[\pt]{\meritfun[{\barrparam[]}]}\paren*{\zerovec[\pt]}\sbra*{\dirpt} - \int_{0}^{1} \paren*{1-\tmethr} \D[2]\pullback[\pt]{\meritfun[{\barrparam[]}]}\paren*{\zerovec[\pt] + \tmethr\dirpt}\sbra*{\dirpt, \dirpt} \dd{\tmethr} + \frac{1}{2}\metr[\pt]{\trsquad\paren*{\allvar}\sbra*{\dirpt}}{\dirpt}\\
            &+ \metr[\pt]{\trslin[{\barrparam[]}]\paren*{\pt}}{\dir[\pt]}} = \abs*{\int_{0}^{1}\paren*{1-\tmethr} \paren*{ \metr[\pt]{\trsquad\paren*{\allvar}\sbra*{\dirpt}}{\dirpt} -  \D[2]\pullback[\pt]{\meritfun[{\barrparam[]}]}\paren*{\tmethr\dirpt}\sbra*{\dirpt, \dirpt}} \dd{\tmethr}}\\
            &\leq \int_{0}^{1} \paren*{1-\tmethr}\Bigg(\abs*{\metr[\pt]{\paren*{\Hess\objfun\paren*{\pt} - \Hess\pullback[\pt]{\objfun}\paren*{\tmethr\dirpt}}\sbra*{\dirpt}}{\dirpt}}\\
            &+ \sum_{\ineqidx\in\ineqset}\abs*{\metr[\pt]{\paren*{-\ineqLagmult[\ineqidx]\Hess\ineqfun[\ineqidx]\paren*{\pt} + \frac{\barrparam[]}{\pullback[\pt]{{\ineqfun[\ineqidx]}}\paren*{\tmethr\dirpt}}\Hess\pullback[\pt]{\ineqfun[\ineqidx]}\paren*{\tmethr\dirpt}}\sbra*{\dirpt}}{\dirpt}}\\
            &+ \sum_{\ineqidx\in\ineqset}\abs*{\frac{\ineqLagmult[\ineqidx]}{\ineqfun[\ineqidx]\paren*{\pt}}\metr[\pt]{\gradstr\ineqfun[\ineqidx]\paren*{\pt}}{\dirpt}^{2} - \frac{\barrparam[]}{{\pullback[\pt]{\ineqfun[\ineqidx]}\paren*{\tmethr\dirpt}}^{2}} \metr[\pt]{\gradstr\pullback[\pt]{\ineqfun[\ineqidx]}\paren*{\tmethr\dirpt}}{\dirpt}^{2}}\Bigg) \dd{\tmethr},
        \end{align}}
        where the second equality follows from Taylor's theorem, the third one \isextendedversion{
        from \cref{eq:dirderivmeritfun},}{from \\ \cref{lemm:firstdirderivmeritfun},}
        and the inequality from 
        \isextendedversion{
        \cref{eq:twicedirderivmeritfun}.
        }{
        \cref{lemm:twicedirderivmeritfun}.}
        In the following, we will derive the bounds on the right-hand side of \cref{ineq:diffpredared} under the settings of \cref{lemm:diffpredaredfixedpt}, \enumicref{lemm:diffpredared}{lemm:diffpredareduniform}, and \enumicref{lemm:diffpredared}{lemm:diffpredaredsmallo}.

        \begin{proof}[Proof of \cref{lemm:diffpredaredfixedpt}]
            It follows from \cref{ineq:diffpredared} that, for all $\dirpt\in\tanspc[\pt]\mani$ sufficiently small,
            \isextendedversion{
            \begin{align}
                \begin{split}
                    &\abs*{\pred[{\allvar, \barrparam[]}]\paren*{\dirpt}-\ared[{\barrparam[]}]\paren*{\dirpt}} \leq \int_{0}^{1} \paren*{1-\tmethr}\Bigg(\opnorm{\Hess\objfun\paren*{\pt}} + \opnorm{\Hess\pullback[\pt]{\objfun}\paren*{\tmethr\dirpt}}\\
                    &\quad + \sum_{\ineqidx\in\ineqset}\paren*{\ineqLagmult[\ineqidx]\opnorm{\Hess\ineqfun[\ineqidx]\paren*{\pt}} + \frac{\barrparam[]}{\pullback[\pt]{{\ineqfun[\ineqidx]}}\paren*{\tmethr\dirpt}}\opnorm{\Hess\pullback[\pt]{\ineqfun[\ineqidx]}\paren*{\tmethr\dirpt}}\\
                    &\quad + \frac{\ineqLagmult[\ineqidx]}{\ineqfun[\ineqidx]\paren*{\pt}}\opnorm{\gradstr\ineqfun[\ineqidx]\paren*{\pt}} + \frac{\barrparam[]}{{\pullback[\pt]{\ineqfun[\ineqidx]}\paren*{\tmethr\dirpt}}^{2}} \opnorm{\gradstr\pullback[\pt]{\ineqfun[\ineqidx]}\paren*{\tmethr\dirpt}}}\Bigg){\Riemnorm[\pt]{\dirpt}^{2}} \dd{\tmethr}.
                \end{split}
            \end{align}}{
            \vspace{-1\baselineskip}
            \begin{align}
                \begin{split}
                    &\abs*{\pred[{\allvar, \barrparam[]}]\paren*{\dirpt}-\ared[{\barrparam[]}]\paren*{\dirpt}} \leq \int_{0}^{1} \paren*{1-\tmethr}\Bigg(\opnorm{\Hess\objfun\paren*{\pt}} + \opnorm{\Hess\pullback[\pt]{\objfun}\paren*{\tmethr\dirpt}}\\
                    &\quad + \sum_{\ineqidx\in\ineqset}\paren*{\ineqLagmult[\ineqidx]\opnorm{\Hess\ineqfun[\ineqidx]\paren*{\pt}} + \frac{\barrparam[]}{\pullback[\pt]{{\ineqfun[\ineqidx]}}\paren*{\tmethr\dirpt}}\opnorm{\Hess\pullback[\pt]{\ineqfun[\ineqidx]}\paren*{\tmethr\dirpt}}\\
                    &\quad + \frac{\ineqLagmult[\ineqidx]}{\ineqfun[\ineqidx]\paren*{\pt}}\opnorm{\gradstr\ineqfun[\ineqidx]\paren*{\pt}} + \frac{\barrparam[]}{{\pullback[\pt]{\ineqfun[\ineqidx]}\paren*{\tmethr\dirpt}}^{2}} \opnorm{\gradstr\pullback[\pt]{\ineqfun[\ineqidx]}\paren*{\tmethr\dirpt}}}\Bigg){\Riemnorm[\pt]{\dirpt}^{2}} \dd{\tmethr}.
                \end{split}
            \end{align}
            }
            Let $\coeffone > 0$ be sufficiently large.
            Since the sequences $\opnorm{\Hess\pullback[\pt]{\objfun}\paren*{\cdot}}$, $\brc*{\pullback[\pt]{{\ineqfun[\ineqidx]}}\paren*{\cdot}}_{\ineqidx\in\ineqset}$, $\brc*{\opnorm{\Hess\pullback[\pt]{\ineqfun[\ineqidx]}\paren*{\cdot}}}_{\ineqidx\in\ineqset}$,  and $\brc*{\norm[\cdot]{\gradstr\pullback[\pt]{\ineqfun[\ineqidx]}\paren*{\cdot}}}_{\ineqidx\in\ineqset}$ are all continuous, the right-hand side can be bounded above by $\coeffone{\Riemnorm[\pt]{\dirpt}}^{2}$ for all  $\dirpt\in\tanspc[\pt]\mani$ sufficiently small.
            The proof is complete.
        \end{proof}

        We proceed to bound \cref{ineq:diffpredared}.
        Under \cref{assu:retrsecondord,assu:radialLCtwo}, the first term of the right-hand side of \cref{ineq:diffpredared} is bounded as follows: 
        if $\dirpt\in\tanspc[\pt]\mani$ satisfies $\Riemnorm[\pt]{\dirpt}\leq\deltaRLtwo[\objfun]$, then we have
        \begin{equation}\label{ineq:boundhessobjfun}
            \abs*{\metr[\pt]{\paren*{\Hess\objfun\paren*{\pt} - \Hess\pullback[\pt]{\objfun}\paren*{\tmethr\dirpt}}\sbra*{\dirpt}}{\dirpt}} \leq \betaRLtwo[\objfun] \Riemnorm[\ptinriter]{\dirpt}^{3}
        \end{equation}
        for $\tmefour \in \sbra*{0, 1}$, where the inequality follows from \cref{prop:secodrretrHess}, $\tmefour\Riemnorm[\pt]{\dirpt}\leq\Riemnorm[\pt]{\dirpt}\leq\deltaRLtwo[\objfun]$, and \cref{def:radialLCtwo} with $\tmefour \leq 1$.
        For the third term of the right-hand side of \cref{ineq:diffpredared}, we obtain the following bound: for each $\ineqidx\in\ineqset$, all $\pt\in\strictfeasirgn$, $\tmethr\in\sbra*{0, 1}$, and all $\dirpt\in\tanspc[\pt]\mani$ sufficiently small, 
        \isextendedversion{
        \begin{align}
            \begin{split}\label{ineq:basicbounddiffsquaregradineqfun}
                &\abs*{\frac{\ineqLagmult[\ineqidx]}{\ineqfun[\ineqidx]\paren*{\pt}}\metr[\pt]{\gradstr\ineqfun[\ineqidx]\paren*{\pt}}{\dirpt}^{2} - \frac{\barrparam[]}{{\pullback[\pt]{\ineqfun[\ineqidx]}\paren*{\tmethr\dirpt}}^{2}} \metr[\pt]{\gradstr\pullback[\pt]{\ineqfun[\ineqidx]}\paren*{\tmethr\dirpt}}{\dirpt}^{2}}\\
                &\leq \abs*{\frac{\ineqLagmult[\ineqidx]}{\ineqfun[\ineqidx]\paren*{\pt}} - \frac{\barrparam[]}{{\pullback[\pt]{\ineqfun[\ineqidx]}\paren*{\tmethr\dirpt}}^{2}}} \metr[\pt]{\gradstr\ineqfun[\ineqidx]\paren*{\pt}}{\dirpt}^{2}\\
                &\quad + \frac{\barrparam[]}{{\pullback[\pt]{\ineqfun[\ineqidx]}\paren*{\tmethr\dirpt}}^{2}} \abs*{\metr[\pt]{\gradstr\ineqfun[\ineqidx]\paren*{\pt}}{\dirpt}^{2} -\metr[\pt]{\gradstr\pullback[\pt]{\ineqfun[\ineqidx]}\paren*{\tmethr\dirpt}}{\dirpt}^{2}}\\
                &= \abs*{\frac{\ineqLagmult[\ineqidx]}{\ineqfun[\ineqidx]\paren*{\pt}} - \frac{\barrparam[]}{{\pullback[\pt]{\ineqfun[\ineqidx]}\paren*{\tmethr\dirpt}}^{2}}} \metr[\pt]{\gradstr\ineqfun[\ineqidx]\paren*{\pt}}{\dirpt}^{2}\\
                &\quad + \frac{\barrparam[]}{{\pullback[\pt]{\ineqfun[\ineqidx]}\paren*{\tmethr\dirpt}}^{2}}\paren*{\abs*{\metr[\pt]{\gradstr\ineqfun[\ineqidx]\paren*{\pt}}{\dirpt} -\metr[\pt]{\gradstr\pullback[\pt]{\ineqfun[\ineqidx]}\paren*{\tmethr\dirpt}}{\dirpt}} \\
                &\qquad \cdot \abs*{\metr[\pt]{\gradstr\pullback[\pt]{\ineqfun[\ineqidx]}\paren*{\tmethr\dirpt} - \gradstr\ineqfun[\ineqidx]\paren*{\pt} + 2 \gradstr\ineqfun[\ineqidx]\paren*{\pt}}{\dirpt}}}, 
            \end{split}
        \end{align}}{\begin{align}
            &\label{ineq:basicbounddiffsquaregradineqfun}\\
            &\abs*{\frac{\ineqLagmult[\ineqidx]}{\ineqfun[\ineqidx]\paren*{\pt}}\metr[\pt]{\gradstr\ineqfun[\ineqidx]\paren*{\pt}}{\dirpt}^{2} - \frac{\barrparam[]}{{\pullback[\pt]{\ineqfun[\ineqidx]}\paren*{\tmethr\dirpt}}^{2}} \metr[\pt]{\gradstr\pullback[\pt]{\ineqfun[\ineqidx]}\paren*{\tmethr\dirpt}}{\dirpt}^{2}} \leq \abs*{\frac{\ineqLagmult[\ineqidx]}{\ineqfun[\ineqidx]\paren*{\pt}} - \frac{\barrparam[]}{{\pullback[\pt]{\ineqfun[\ineqidx]}\paren*{\tmethr\dirpt}}^{2}}} \\
            &\metr[\pt]{\gradstr\ineqfun[\ineqidx]\paren*{\pt}}{\dirpt}^{2} + \frac{\barrparam[]}{{\pullback[\pt]{\ineqfun[\ineqidx]}\paren*{\tmethr\dirpt}}^{2}} \abs*{\metr[\pt]{\gradstr\ineqfun[\ineqidx]\paren*{\pt}}{\dirpt}^{2} -\metr[\pt]{\gradstr\pullback[\pt]{\ineqfun[\ineqidx]}\paren*{\tmethr\dirpt}}{\dirpt}^{2}}\\
            &= \abs*{\frac{\ineqLagmult[\ineqidx]}{\ineqfun[\ineqidx]\paren*{\pt}} - \frac{\barrparam[]}{{\pullback[\pt]{\ineqfun[\ineqidx]}\paren*{\tmethr\dirpt}}^{2}}} \metr[\pt]{\gradstr\ineqfun[\ineqidx]\paren*{\pt}}{\dirpt}^{2} + \frac{\barrparam[]}{{\pullback[\pt]{\ineqfun[\ineqidx]}\paren*{\tmethr\dirpt}}^{2}}\paren*{\abs*{\metr[\pt]{\gradstr\ineqfun[\ineqidx]\paren*{\pt}}{\dirpt}\\
            &-\metr[\pt]{\gradstr\pullback[\pt]{\ineqfun[\ineqidx]}\paren*{\tmethr\dirpt}}{\dirpt}} \abs*{\metr[\pt]{\gradstr\pullback[\pt]{\ineqfun[\ineqidx]}\paren*{\tmethr\dirpt} - \gradstr\ineqfun[\ineqidx]\paren*{\pt} + 2 \gradstr\ineqfun[\ineqidx]\paren*{\pt}}{\dirpt}}}, 
        \end{align}}
        where the inequality follows from $\barrparam[] > 0$.
        In the following, we provide the proofs of \cref{lemm:diffpredared} by developing these bounds, together with that of the second term of the right-hand side of \cref{ineq:diffpredared}, under the settings of \cref{lemm:diffpredareduniform,lemm:diffpredaredsmallo}, respectively.
        
        \begin{proof}[Proof of \enumicref{lemm:diffpredared}{lemm:diffpredareduniform}]
            Define 
            \isextendedversion{
            \begin{equation}
                \thldtrradiusfou\coloneqq\min\brc*{1, \deltathr, \deltaRLtwo[\objfun],  \brc*{\deltaRL[{\ineqfun[\ineqidx]}]}_{\ineqidx\in\ineqset}, \brc*{\deltaRLtwo[{\ineqfun[\ineqidx]}]}_{\ineqidx\in\ineqset}} > 0,    
            \end{equation}}{
            $\thldtrradiusfou\hspace{-1mm}\coloneqq\hspace{-1mm}\min\brc*{1, \deltathr, \deltaRLtwo[\objfun],  \brc*{\deltaRL[{\ineqfun[\ineqidx]}]}_{\ineqidx\in\ineqset}, \\ \brc*{\deltaRLtwo[{\ineqfun[\ineqidx]}]}_{\ineqidx\in\ineqset}} > 0$,}
            where $\deltathr > 0$ is the threshold associated with $\epsthr > 0$ in \enumicref{lemm:ineqfunbound}{lemm:ienqfununiflowbd} and $\deltaRLtwo[\objfun] > 0$, $\brc*{\deltaRL[{\ineqfun[\ineqidx]}]}_{\ineqidx\in\ineqset}$, and $\brc*{\deltaRLtwo[{\ineqfun[\ineqidx]}]}_{\ineqidx\in\ineqset}$ are the ones for \cref{def:radialLCone} and \cref{def:radialLCtwo}, respectively.
            Let $\dirptinriter\in\tanspc[\ptinriter]\mani$ be a search direction satisfying $\Riemnorm[\ptinriter]{\dirptinriter} \leq \thldtrradiusfou$.

            As for the second term of the right-hand side of \cref{ineq:diffpredared}, we have that, for each $\ineqidx\in\ineqset$, 
            \isextendedversion{
            \begin{equation}
                \begin{aligned}\label{ineq:uniboundhessineqfun}
                    &\abs*{\metr[\ptinriter]{\paren*{-\ineqLagmultinriter[\ineqidx]\Hess\ineqfun[\ineqidx]\paren*{\ptinriter} + \frac{\barrparam[]}{\pullback[\ptinriter]{\ineqfun[\ineqidx]}\paren*{\tmethr\dirptinriter}}\Hess\pullback[\ptinriter]{\ineqfun[\ineqidx]}\paren*{\tmethr\dirptinriter}}\sbra*{\dirptinriter}}{\dirptinriter}}\\
                    &\leq \abs*{\paren*{-\ineqLagmultinriter[\ineqidx] + \frac{\barrparam[]}{\pullback[\ptinriter]{\ineqfun[\ineqidx]}\paren*{\tmethr\dirptinriter}}}  \metr[\ptinriter]{\Hess\ineqfun[\ineqidx]\paren*{\ptinriter}\sbra*{\dirptinriter}}{\dirptinriter}}\\
                    &\quad + \abs*{\frac{\barrparam[]}{\pullback[\ptinriter]{{\ineqfun[\ineqidx]}}\paren*{\tmethr\dirptinriter}}\metr[\ptinriter]{\paren*{- \Hess\pullback[\ptinriter]{\ineqfun[\ineqidx]}\paren*{\zerovec[\ptinriter]} + \Hess\pullback[\ptinriter]{\ineqfun[\ineqidx]}\paren*{\tmethr\dirptinriter}}\sbra*{\dirptinriter}}{\dirptinriter}}\\
                    &\leq \paren*{\ineqLagmultinriter[\ineqidx] + \frac{\barrparam[]}{\epsthr}}\opnorm{\Hess\ineqfun[\ineqidx]\paren*{\ptinriter}} \Riemnorm[\ptinriter]{\dirptinriter}^{2} + \frac{\barrparam[]\betaRLtwo[{\ineqfun[\ineqidx]}]}{\epsthr}\Riemnorm[\ptinriter]{\dirptinriter}^{3},
                \end{aligned}
            \end{equation}}{\begin{align}
                &\label{ineq:uniboundhessineqfun}\\
                &\abs*{\metr[\ptinriter]{\paren*{\hspace{-1mm}-\ineqLagmultinriter[\ineqidx]\Hess\ineqfun[\ineqidx]\paren*{\ptinriter} \hspace{-1mm} + \hspace{-1mm}\frac{\barrparam[]}{\pullback[\ptinriter]{\ineqfun[\ineqidx]}\paren*{\tmethr\dirptinriter}}\Hess\pullback[\ptinriter]{\ineqfun[\ineqidx]}\paren*{\tmethr\dirptinriter}}\sbra*{\dirptinriter}}{\dirptinriter\hspace{-1mm}}} \leq \abs*{\paren*{\hspace{-1mm}-\ineqLagmultinriter[\ineqidx] \hspace{-1mm} + \hspace{-1mm}\frac{\barrparam[]}{\pullback[\ptinriter]{\ineqfun[\ineqidx]}\paren*{\tmethr\dirptinriter}}} \\
                &\metr[\ptinriter]{\Hess\ineqfun[\ineqidx]\paren*{\ptinriter}\sbra*{\dirptinriter}}{\dirptinriter}} + \abs*{\frac{\barrparam[]}{\pullback[\ptinriter]{{\ineqfun[\ineqidx]}}\paren*{\tmethr\dirptinriter}}
                \ang*{\paren*{- \Hess\pullback[\ptinriter]{\ineqfun[\ineqidx]}\paren*{\zerovec[\ptinriter]} + \Hess\pullback[\ptinriter]{\ineqfun[\ineqidx]}\paren*{\tmethr\dirptinriter}}\sbra*{\dirptinriter}, \\ &\dirptinriter}_{\ptinriter}} \leq \paren*{\ineqLagmultinriter[\ineqidx] + \frac{\barrparam[]}{\epsthr}}\opnorm{\Hess\ineqfun[\ineqidx]\paren*{\ptinriter}} \Riemnorm[\ptinriter]{\dirptinriter}^{2} + \frac{\barrparam[]\betaRLtwo[{\ineqfun[\ineqidx]}]}{\epsthr}\Riemnorm[\ptinriter]{\dirptinriter}^{3},
            \end{align}}
            where the first inequality follows from \cref{prop:secodrretrHess} and the second one from \cref{def:radialLCtwo} with $\tmethr \leq 1$, $\barrparam[], \ineqLagmultinriter[\ineqidx] \in \setRpp[]$, and $\pullback[\ptinriter]{\ineqfun[\ineqidx]}\paren*{\tmethr\dirptinriter} > \epsthr$ by \enumicref{lemm:ineqfunbound}{lemm:ienqfununiflowbd}.
            Using \cref{ineq:basicbounddiffsquaregradineqfun}, the third term of the right-hand side of \cref{ineq:diffpredared} can be bounded as
            \isextendedversion{
            \begin{equation}
                \begin{aligned}\label{ineq:unibounddiffsquaregradineqfun}
                    &\abs*{\frac{\ineqLagmultinriter[\ineqidx]}{\ineqfun[\ineqidx]\paren*{\ptinriter}}\metr[\ptinriter]{\gradstr\ineqfun[\ineqidx]\paren*{\ptinriter}}{\dirptinriter}^{2} - \frac{\barrparam[]}{{\pullback[\ptinriter]{\ineqfun[\ineqidx]}\paren*{\tmethr\dirptinriter}}^{2}} \metr[\ptinriter]{\gradstr\pullback[\ptinriter]{\ineqfun[\ineqidx]}\paren*{\tmethr\dirptinriter}}{\dirptinriter}^{2}}\\
                    &\leq \paren*{\frac{\ineqLagmultinriter[\ineqidx]}{\ineqfun[\ineqidx]\paren*{\ptinriter}} + \frac{\barrparam[]}{{\pullback[\ptinriter]{\ineqfun[\ineqidx]}\paren*{\tmethr\dirptinriter}}^{2}}} \metr[\ptinriter]{\gradstr\ineqfun[\ineqidx]\paren*{\ptinriter}}{\dirptinriter}^{2}\\
                    &\quad + \frac{\barrparam[]}{{\pullback[\ptinriter]{\ineqfun[\ineqidx]}\paren*{\tmethr\dirptinriter}}^{2}} \abs*{\metr[\ptinriter]{\gradstr\pullback[\ptinriter]{\ineqfun[\ineqidx]}\paren*{\zerovec[\ptinriter]} - \gradstr\pullback[\ptinriter]{\ineqfun[\ineqidx]}\paren*{\tmethr\dirptinriter}}{\dirptinriter}}\\
                    &\qquad\paren*{\abs*{\metr[\ptinriter]{\gradstr\pullback[\ptinriter]{\ineqfun[\ineqidx]}\paren*{\tmethr\dirptinriter} - \gradstr\pullback[\ptinriter]{\ineqfun[\ineqidx]}\paren*{\zerovec[\ptinriter]}}{\dirptinriter}} + 2\abs*{\metr[\ptinriter]{\gradstr\ineqfun[\ineqidx]\paren*{\ptinriter}}{\dirptinriter}}}\\
                    &\leq \paren*{\frac{\ineqLagmultinriter[\ineqidx]}{\epsthr} + \frac{\barrparam[]}{\epsthr^{2}}} \Riemnorm[\ptinriter]{\gradstr\ineqfun[\ineqidx]\paren*{\ptinriter}}^{2} \Riemnorm[\ptinriter]{\dirptinriter}^{2}\\
                    &\quad + \frac{\barrparam[]\paren*{\betaRL[{\ineqfun[\ineqidx]}]}^{2}}{\epsthr^{2}} \Riemnorm[\ptinriter]{\dirptinriter}^{4} + \frac{2\barrparam[]\betaRL[{\ineqfun[\ineqidx]}]}{\epsthr^{2}} \Riemnorm[\ptinriter]{\gradstr\ineqfun[\ineqidx]\paren*{\ptinriter}}\Riemnorm[\ptinriter]{\dirptinriter}^{3},
                \end{aligned}
            \end{equation}}{\begin{align}
                &\label{ineq:unibounddiffsquaregradineqfun}\\
                &\abs*{\frac{\ineqLagmultinriter[\ineqidx]}{\ineqfun[\ineqidx]\paren*{\ptinriter}}\metr[\ptinriter]{\gradstr\ineqfun[\ineqidx]\paren*{\ptinriter}}{\dirptinriter}^{2} - \frac{\barrparam[]}{{\pullback[\ptinriter]{\ineqfun[\ineqidx]}\paren*{\tmethr\dirptinriter}}^{2}} \metr[\ptinriter]{\gradstr\pullback[\ptinriter]{\ineqfun[\ineqidx]}\paren*{\tmethr\dirptinriter}}{\dirptinriter}^{2}} \leq \paren*{\frac{\ineqLagmultinriter[\ineqidx]}{\ineqfun[\ineqidx]\paren*{\ptinriter}} \\
                &+ \frac{\barrparam[]}{{\pullback[\ptinriter]{\ineqfun[\ineqidx]}\paren*{\tmethr\dirptinriter}}^{2}}} \metr[\ptinriter]{\gradstr\ineqfun[\ineqidx]\paren*{\ptinriter}}{\dirptinriter}^{2} + \frac{\barrparam[]}{{\pullback[\ptinriter]{\ineqfun[\ineqidx]}\paren*{\tmethr\dirptinriter}}^{2}} \abs*{
                \ang*{ \gradstr\pullback[\ptinriter]{\ineqfun[\ineqidx]}\paren*{\zerovec[\ptinriter]} - \gradstr\pullback[\ptinriter]{\ineqfun[\ineqidx]}\paren*{\tmethr\dirptinriter}, \\ &\dirptinriter}_{\ptinriter}}
                \paren*{\abs*{\metr[\ptinriter]{\gradstr\pullback[\ptinriter]{\ineqfun[\ineqidx]}\paren*{\tmethr\dirptinriter} - \gradstr\pullback[\ptinriter]{\ineqfun[\ineqidx]}\paren*{\zerovec[\ptinriter]}}{\dirptinriter}} + 2\abs*{\metr[\ptinriter]{\gradstr\ineqfun[\ineqidx]\paren*{\ptinriter}}{\dirptinriter}}} \leq \paren*{\frac{\ineqLagmultinriter[\ineqidx]}{\epsthr}\\
                &+ \frac{\barrparam[]}{\epsthr^{2}}} \Riemnorm[\ptinriter]{\gradstr\ineqfun[\ineqidx]\paren*{\ptinriter}}^{2} \Riemnorm[\ptinriter]{\dirptinriter}^{2} \hspace{-1mm} + \frac{\barrparam[]\paren*{\betaRL[{\ineqfun[\ineqidx]}]}^{2}}{\epsthr^{2}} \Riemnorm[\ptinriter]{\dirptinriter}^{4} \hspace{-1mm} + \frac{2\barrparam[]\betaRL[{\ineqfun[\ineqidx]}]}{\epsthr^{2}} \Riemnorm[\ptinriter]{\gradstr\ineqfun[\ineqidx]\paren*{\ptinriter}}\Riemnorm[\ptinriter]{\dirptinriter}^{3},
            \end{align}}
            where the first inequality follows from $\ineqLagmultinriter[\ineqidx], \barrparam[] \in\setRpp[]$, $\ineqfun[\ineqidx]\paren*{\ptinriter} > 0$, and \cref{eq:pullbackgradzero} and the second one from \cref{def:radialLCone} with $\tmethr \leq 1$ and $\ineqfun[\ineqidx]\paren*{\ptinriter} > \epsthr$ and $\pullback[\ptinriter]{\ineqfun[\ineqidx]}\paren*{\tmethr\dirptinriter} > \epsthr$ by \enumicref{lemm:ineqfunbound}{lemm:ienqfununiflowbd}.
            By combining \cref{ineq:diffpredared} with \cref{ineq:boundhessobjfun,ineq:uniboundhessineqfun,ineq:unibounddiffsquaregradineqfun}, we obtain
            \isextendedversion{
            \begin{align}
                &\abs*{\predinriter-\aredinriter}\\
                &\leq \int_{0}^{1} \paren*{1-\tmethr} \Bigg(\betaRLtwo[\objfun] \Riemnorm[\ptinriter]{\dirptinriter}^{3}\\
                &\quad + \sum_{\ineqidx\in\ineqset} \paren*{\paren*{\ineqLagmultinriter[\ineqidx] + \frac{\barrparam[]}{\epsthr}}\opnorm{\Hess\ineqfun[\ineqidx]\paren*{\ptinriter}} \Riemnorm[\ptinriter]{\dirptinriter}^{2} + \frac{\barrparam[]\betaRLtwo[{\ineqfun[\ineqidx]}]}{\epsthr}\Riemnorm[\ptinriter]{\dirptinriter}^{3}}\\
                &\quad + \sum_{\ineqidx\in\ineqset}\paren*{\paren*{\frac{\ineqLagmultinriter[\ineqidx]}{\epsthr} + \frac{\barrparam[]}{\epsthr^{2}}} \Riemnorm[\ptinriter]{\gradstr\ineqfun[\ineqidx]\paren*{\ptinriter}}^{2}\Riemnorm[\ptinriter]{\dirptinriter}^{2} \\
                &\quad + \frac{\barrparam[]\paren*{\betaRL[{\ineqfun[\ineqidx]}]}^{2}}{\epsthr^{2}} \Riemnorm[\ptinriter]{\dirptinriter}^{4} + \frac{2\barrparam[]\betaRL[{\ineqfun[\ineqidx]}]}{\epsthr^{2}} \Riemnorm[\ptinriter]{\gradstr\ineqfun[\ineqidx]\paren*{\ptinriter}}\Riemnorm[\ptinriter]{\dirptinriter}^{3}}\Bigg) \dd{\tmethr}\\
                &\leq \frac{1}{2}  \paren*{\betaRLtwo[\objfun] + \sum_{\ineqidx\in\ineqset} \paren*{\paren*{\ineqLagmultinriter[\ineqidx] + \frac{\barrparam[]}{\epsthr}}\opnorm{\Hess\ineqfun[\ineqidx]\paren*{\ptinriter}}  + \frac{\barrparam[]\betaRLtwo[{\ineqfun[\ineqidx]}]}{\epsthr}}\\
                &\quad + \sum_{\ineqidx\in\ineqset}\paren*{\paren*{\frac{\ineqLagmultinriter[\ineqidx]}{\epsthr} + \frac{\barrparam[]}{\epsthr^{2}}} \Riemnorm[\ptinriter]{\gradstr\ineqfun[\ineqidx]\paren*{\ptinriter}}^{2}\\
                &\quad + \frac{\barrparam[]\paren*{\betaRL[{\ineqfun[\ineqidx]}]}^{2}}{\epsthr^{2}} + \frac{2\barrparam[]\betaRL[{\ineqfun[\ineqidx]}]}{\epsthr^{2}} \Riemnorm[\ptinriter]{\gradstr\ineqfun[\ineqidx]\paren*{\ptinriter}}}} \Riemnorm[\ptinriter]{\dirptinriter}^{2},
            \end{align}
            }{
            \begin{align}
                &\abs*{\predinriter-\aredinriter} \leq \int_{0}^{1} \paren*{1-\tmethr} \Bigg(\betaRLtwo[\objfun] \Riemnorm[\ptinriter]{\dirptinriter}^{3} + \sum_{\ineqidx\in\ineqset} \paren*{\paren*{\ineqLagmultinriter[\ineqidx] + \frac{\barrparam[]}{\epsthr}}\opnorm{\Hess\ineqfun[\ineqidx]\paren*{\ptinriter}} \Riemnorm[\ptinriter]{\dirptinriter}^{2} \\
                &+ \frac{\barrparam[]\betaRLtwo[{\ineqfun[\ineqidx]}]}{\epsthr}\Riemnorm[\ptinriter]{\dirptinriter}^{3}} + \sum_{\ineqidx\in\ineqset}\paren*{\paren*{\frac{\ineqLagmultinriter[\ineqidx]}{\epsthr} + \frac{\barrparam[]}{\epsthr^{2}}} \Riemnorm[\ptinriter]{\gradstr\ineqfun[\ineqidx]\paren*{\ptinriter}}^{2}\Riemnorm[\ptinriter]{\dirptinriter}^{2} + \frac{\barrparam[]\paren*{\betaRL[{\ineqfun[\ineqidx]}]}^{2}}{\epsthr^{2}} \Riemnorm[\ptinriter]{\dirptinriter}^{4}\\
                &+ \frac{2\barrparam[]\betaRL[{\ineqfun[\ineqidx]}]}{\epsthr^{2}} \Riemnorm[\ptinriter]{\gradstr\ineqfun[\ineqidx]\paren*{\ptinriter}}\Riemnorm[\ptinriter]{\dirptinriter}^{3}}\Bigg) \dd{\tmethr} \leq \frac{1}{2}  \paren*{\betaRLtwo[\objfun] + \sum_{\ineqidx\in\ineqset} \paren*{\paren*{\ineqLagmultinriter[\ineqidx] + \frac{\barrparam[]}{\epsthr}}\opnorm{\Hess\ineqfun[\ineqidx]\paren*{\ptinriter}} \\
                &+ \frac{\barrparam[]\betaRLtwo[{\ineqfun[\ineqidx]}]}{\epsthr}} + \sum_{\ineqidx\in\ineqset}\paren*{\paren*{\frac{\ineqLagmultinriter[\ineqidx]}{\epsthr} + \frac{\barrparam[]}{\epsthr^{2}}} \Riemnorm[\ptinriter]{\gradstr\ineqfun[\ineqidx]\paren*{\ptinriter}}^{2} + \frac{\barrparam[]\paren*{\betaRL[{\ineqfun[\ineqidx]}]}^{2}}{\epsthr^{2}}\\
                &+ \frac{2\barrparam[]\betaRL[{\ineqfun[\ineqidx]}]}{\epsthr^{2}} \Riemnorm[\ptinriter]{\gradstr\ineqfun[\ineqidx]\paren*{\ptinriter}}}} \Riemnorm[\ptinriter]{\dirptinriter}^{2},
            \end{align}
            }
            where the second inequality follows from $\Riemnorm[\ptinriter]{\dirptinriter} \leq \thldtrradiusfou \leq 1$.
            From \assuenumicref{assu:ineqfunbounded}{assu:ineqfunboundedgradhess} and \cref{assu:ineqLagmultbounded}, all the terms in the coefficients of the right-hand side are bounded above.
            Thus, by taking $\coefftwo > 0$ sufficiently large,  we complete the proof of \enumicref{lemm:diffpredared}{lemm:diffpredareduniform}.
        \end{proof}

        \begin{proof}[Proof of \enumicref{lemm:diffpredared}{lemm:diffpredaredsmallo}]
            By \assuenumicref{assu:ineqfunbounded}{assu:ineqfunboundedgradhess} and \cref{assu:ineqLagmultbounded}, there exist positive scalars $\tholdvalthr, \tholdvalfou, \tholdvalfiv \in \setRpp[]$ such that
            \begin{equation}\label{ineq:gradhessineqLagmultbound}
                \Riemnorm[\ptinriter]{\gradstr[]\ineqfun[\ineqidx]\paren*{\ptinriter}} \leq \tholdvalthr, \quad \Riemnorm[\ptinriter]{\Hess[]\ineqfun[\ineqidx]\paren*{\ptinriter}} \leq \tholdvalfou, \quad \ineqLagmultinriter[\ineqidx] \leq \tholdvalfiv \text{ for all } \ineqidx\in\ineqset.
            \end{equation}
            Let $\coefffou > 0$ be any positive scalar, and $\epsthr, \deltathr \in \setRpp[]$ be the positive values from \enumicref{lemm:ineqfunbound}{lemm:ienqfununiflowbd}, respectively.
            Under \cref{assu:pullbackineqfuncont}, we choose $\thldtrradiusthr > 0$ so that, for any $\inriteridx\in\setNz$ and all $\tanvecone[\ptinriter]\in\tanspc[\ptinriter]\mani$ with $\Riemnorm[\ptinriter]{\tanvecone[\ptinriter]} \leq \thldtrradiusthr$,
            \begin{equation}\label{ineq:pullbackineqfuncoefffoubound}
                \abs*{\pullback[\ptinriter]{\ineqfun[\ineqidx]}\paren*{\tanvecone[\ptinriter]} - \ineqfun[\ineqidx]\paren*{\ptinriter}} \leq \frac{2\coefffou\epsthr^{2}}{9\ineqdime}\min\brc*{\frac{1}{\barrparam[]\tholdvalfou}, \frac{1}{3\tholdvalfiv\tholdvalthr^{2}}, \frac{\epsthr}{3\barrparam[]\tholdvalthr^{2}}}.
            \end{equation}
            From \cref{assu:ineqLagmultconverged}, we denote by $\tholdidxtwo\in\setNz$ an index satisfying that, for all $\inriteridx\geq\tholdidxtwo$,
            \begin{equation}\label{ineq:ineqLagmultcoefffoubound}
                \norm{\ineqLagmultinriter[]-\barrparam[]\inv{\Ineqfunmat[]}\paren*{\ptinriter}\onevec} \leq \frac{2\coefffou}{9\ineqdime} \min\brc*{\frac{1}{\tholdvalfou},\frac{\epsthr}{3\tholdvalthr^{2}}}.
            \end{equation}
            Let
            \isextendedversion{
            \begin{align}
                \begin{split}\label{eq:thldtrradiussixdef}
                    &\thldtrradiussix\coloneqq\min\Bigg\{ \deltathr, \thldtrradiusthr, \deltaRLtwo[\objfun], \brc*{\deltaRL[{\ineqfun[\ineqidx]}]}_{\ineqidx\in\ineqset}, \brc*{\deltaRLtwo[{\ineqfun[\ineqidx]}]}_{\ineqidx\in\ineqset},\\
                    &\quad\brc*{\frac{2\coefffou\epsthr}{9\barrparam[]\ineqdime\betaRLtwo[{\ineqfun[\ineqidx]}]}}_{\ineqidx\in\ineqset}, \brc*{\frac{\coefffou\epsthr^{2}}{9\barrparam[]\ineqdime\betaRL[{\ineqfun[\ineqidx]}]\tholdvalthr}}_{\ineqidx\in\ineqset}, \brc*{\frac{\epsthr\sqrt{2\coefffou}}{3\ineqdime\betaRL[{\ineqfun[\ineqidx]}]\sqrt{\barrparam[]}}}_{\ineqidx\in\ineqset}, \frac{2\coefffou}{3\betaRLtwo[\objfun]} \Bigg\} > 0
                \end{split}
            \end{align}}{$\thldtrradiussix\coloneqq\min\Bigg\{ \deltathr, \thldtrradiusthr, \deltaRLtwo[\objfun], \brc*{\deltaRL[{\ineqfun[\ineqidx]}]}_{\ineqidx\in\ineqset}, \brc*{\deltaRLtwo[{\ineqfun[\ineqidx]}]}_{\ineqidx\in\ineqset}, \brc*{\frac{2\coefffou\epsthr}{9\barrparam[]\ineqdime\betaRLtwo[{\ineqfun[\ineqidx]}]}}_{\ineqidx\in\ineqset}, \brc*{\frac{\coefffou\epsthr^{2}}{9\barrparam[]\ineqdime\betaRL[{\ineqfun[\ineqidx]}]\tholdvalthr}}_{\ineqidx\in\ineqset},\\ 
            \brc*{\frac{\epsthr\sqrt{2\coefffou}}{3\ineqdime\betaRL[{\ineqfun[\ineqidx]}]\sqrt{\barrparam[]}}}_{\ineqidx\in\ineqset}, \frac{2\coefffou}{3\betaRLtwo[\objfun]} \Bigg\} > 0$}
            and let $\dirptinriter\in\tanspc[\ptinriter]\mani$ be a search direction satisfying $\Riemnorm[\ptinriter]{\dirptinriter} \leq \thldtrradiussix$.
            We derive the bounds on each term of the right-hand side of \cref{ineq:diffpredared} under the assumptions.
            Note that, under \cref{assu:retrsecondord,assu:radialLCtwo}, the bound \cref{ineq:boundhessobjfun} also holds.
            As for the second term of the right-hand side of \cref{ineq:diffpredared}, we have the following inequality: for each $\ineqidx\in\ineqset$, it holds that
            \isextendedversion{
            \begin{equation}
                \begin{aligned}[t]\label{ineq:coefffouboundhessineqfun}
                    &\abs*{\metr[\ptinriter]{\paren*{-\ineqLagmultinriter[\ineqidx]\Hess\ineqfun[\ineqidx]\paren*{\ptinriter} + \frac{\barrparam[]}{\pullback[\ptinriter]{\ineqfun[\ineqidx]}\paren*{\tmethr\dirptinriter}}\Hess\pullback[\ptinriter]{\ineqfun[\ineqidx]}\paren*{\tmethr\dirptinriter}}\sbra*{\dirptinriter}}{\dirptinriter}}\\
                    &\leq \paren*{\abs*{\frac{\barrparam[]}{\ineqfun[\ineqidx]\paren*{\ptinriter}}-\ineqLagmultinriter[\ineqidx]} + \abs*{\frac{\barrparam[]\paren*{\pullback[\ptinriter]{\ineqfun[\ineqidx]}\paren*{\tmethr\dirptinriter} - \ineqfun[\ineqidx]\paren*{\ptinriter}}}{\pullback[\ptinriter]{\ineqfun[\ineqidx]}\paren*{\tmethr\dirptinriter}\ineqfun[\ineqidx]\paren*{\ptinriter}}}}\abs*{\metr[\ptinriter]{\Hess\ineqfun[\ineqidx]\paren*{\ptinriter}\sbra*{\dirptinriter}}{\dirptinriter}}\\
                    &\quad + \frac{\barrparam[]}{\pullback[\ptinriter]{{\ineqfun[\ineqidx]}}\paren*{\tmethr\dirptinriter}} \abs*{\metr[\ptinriter]{\paren*{- \Hess\pullback[\ptinriter]{\ineqfun[\ineqidx]}\paren*{\zerovec[\ptinriter]} + \Hess\pullback[\ptinriter]{\ineqfun[\ineqidx]}\paren*{\tmethr\dirptinriter}}\sbra*{\dirptinriter}}{\dirptinriter}}\\
                    &\leq \paren*{\frac{2\coefffou}{9\ineqdime\tholdvalfou} + \frac{2\coefffou}{9\ineqdime\tholdvalfou}}\opnorm{\Hess\ineqfun[\ineqidx]\paren*{\ptinriter}} \Riemnorm[\ptinriter]{\dirptinriter}^{2} + \frac{\barrparam[]\betaRLtwo[{\ineqfun[\ineqidx]}]}{\epsthr}\thldtrradiussix\Riemnorm[\ptinriter]{\dirptinriter}^{2} \leq \frac{2\coefffou}{3\ineqdime}\Riemnorm[\ptinriter]{\dirptinriter}^{2},
                \end{aligned}
            \end{equation}}
            {\begin{align}
                &\label{ineq:coefffouboundhessineqfun}\\
                &\abs*{\metr[\ptinriter]{\paren*{-\ineqLagmultinriter[\ineqidx]\Hess\ineqfun[\ineqidx]\paren*{\ptinriter} + \frac{\barrparam[]}{\pullback[\ptinriter]{\ineqfun[\ineqidx]}\paren*{\tmethr\dirptinriter}}\Hess\pullback[\ptinriter]{\ineqfun[\ineqidx]}\paren*{\tmethr\dirptinriter}}\sbra*{\dirptinriter}}{\dirptinriter}}\\
                &\leq \paren*{\abs*{\frac{\barrparam[]}{\ineqfun[\ineqidx]\paren*{\ptinriter}}-\ineqLagmultinriter[\ineqidx]} + \abs*{\frac{\barrparam[]\paren*{\pullback[\ptinriter]{\ineqfun[\ineqidx]}\paren*{\tmethr\dirptinriter} - \ineqfun[\ineqidx]\paren*{\ptinriter}}}{\pullback[\ptinriter]{\ineqfun[\ineqidx]}\paren*{\tmethr\dirptinriter}\ineqfun[\ineqidx]\paren*{\ptinriter}}}}\abs*{\metr[\ptinriter]{\Hess\ineqfun[\ineqidx]\paren*{\ptinriter}\sbra*{\dirptinriter}}{\dirptinriter}}\\
                &\quad + \frac{\barrparam[]}{\pullback[\ptinriter]{{\ineqfun[\ineqidx]}}\paren*{\tmethr\dirptinriter}} \abs*{\metr[\ptinriter]{\paren*{- \Hess\pullback[\ptinriter]{\ineqfun[\ineqidx]}\paren*{\zerovec[\ptinriter]} + \Hess\pullback[\ptinriter]{\ineqfun[\ineqidx]}\paren*{\tmethr\dirptinriter}}\sbra*{\dirptinriter}}{\dirptinriter}}\\
                &\leq \paren*{\frac{2\coefffou}{9\ineqdime\tholdvalfou} + \frac{2\coefffou}{9\ineqdime\tholdvalfou}}\opnorm{\Hess\ineqfun[\ineqidx]\paren*{\ptinriter}} \Riemnorm[\ptinriter]{\dirptinriter}^{2} + \frac{\barrparam[]\betaRLtwo[{\ineqfun[\ineqidx]}]}{\epsthr}\thldtrradiussix\Riemnorm[\ptinriter]{\dirptinriter}^{2} \leq \frac{2\coefffou}{3\ineqdime}\Riemnorm[\ptinriter]{\dirptinriter}^{2},
            \end{align}
            }
            where the first inequality follows from \cref{prop:secodrretrHess}, $\barrparam[] > 0$, and $\pullback[\ptinriter]{\ineqfun[\ineqidx]}\paren*{\tmethr\dirptinriter} > 0$ by \enumicref{lemm:ineqfunbound}{lemm:ienqfununiflowbd}, the second one from \cref{ineq:ineqLagmultcoefffoubound}, \cref{ineq:pullbackineqfuncoefffoubound}, \cref{def:radialLCtwo} with $\tmethr \leq 1$, $\Riemnorm[\ptinriter]{\dirptinriter}\leq\thldtrradiussix$, and  $\ineqfun[\ineqidx]\paren*{\ptinriter} > \epsthr$ and $\pullback[\ptinriter]{\ineqfun[\ineqidx]}\paren*{\tmethr\dirptinriter} > \epsthr$ by \enumicref{lemm:ineqfunbound}{lemm:ienqfununiflowbd} again, and the last one from 
            \isextendedversion{
            from \cref{ineq:gradhessineqLagmultbound,eq:thldtrradiussixdef}.}{\cref{ineq:gradhessineqLagmultbound} and the definition of $\thldtrradiussix$.}
            Next, we provide a bound on the third term of the right-hand side of \cref{ineq:diffpredared}.
            Using \cref{ineq:basicbounddiffsquaregradineqfun}, we have
            \isextendedversion{\begin{equation}
                \begin{aligned}[t]\label{ineq:coefffoubounddiffsquaregradineqfun}
                    &\abs*{\frac{\ineqLagmultinriter[\ineqidx]}{\ineqfun[\ineqidx]\paren*{\ptinriter}}\metr[\ptinriter]{\gradstr\ineqfun[\ineqidx]\paren*{\ptinriter}}{\dirptinriter}^{2} - \frac{\barrparam[]}{{\pullback[\ptinriter]{\ineqfun[\ineqidx]}\paren*{\tmethr\dirptinriter}}^{2}} \metr[\ptinriter]{\gradstr\pullback[\ptinriter]{\ineqfun[\ineqidx]}\paren*{\tmethr\dirptinriter}}{\dirptinriter}^{2}}\\
                    &\quad\leq \abs*{\frac{\ineqLagmultinriter[\ineqidx]}{\ineqfun[\ineqidx]\paren*{\ptinriter}} - \frac{\barrparam[]}{{\pullback[\ptinriter]{\ineqfun[\ineqidx]}\paren*{\tmethr\dirptinriter}}^{2}}} \metr[\ptinriter]{\gradstr\ineqfun[\ineqidx]\paren*{\ptinriter}}{\dirptinriter}^{2}\\
                    &\quad + \frac{\barrparam[]}{{\pullback[\ptinriter]{\ineqfun[\ineqidx]}\paren*{\tmethr\dirptinriter}}^{2}} \abs*{\metr[\ptinriter]{\gradstr\pullback[\ptinriter]{\ineqfun[\ineqidx]}\paren*{\zerovec[\ptinriter]} - \gradstr\pullback[\ptinriter]{\ineqfun[\ineqidx]}\paren*{\tmethr\dirptinriter}}{\dirptinriter}}\\
                    &\qquad\paren*{\abs*{\metr[\ptinriter]{\gradstr\pullback[\ptinriter]{\ineqfun[\ineqidx]}\paren*{\tmethr\dirptinriter} - \gradstr\pullback[\ptinriter]{\ineqfun[\ineqidx]}\paren*{\zerovec[\ptinriter]}}{\dirptinriter}} + 2\abs*{\metr[\ptinriter]{\gradstr\ineqfun[\ineqidx]\paren*{\ptinriter}}{\dirptinriter}}}\\
                    &\quad\leq \paren*{\frac{\ineqLagmultinriter[\ineqidx]\abs*{\pullback[\ptinriter]{\ineqfun[\ineqidx]}\paren*{\tmethr\dirptinriter} - \ineqfun[\ineqidx]\paren*{\ptinriter}}}{\pullback[\ptinriter]{\ineqfun[\ineqidx]}\paren*{\tmethr\dirptinriter}\ineqfun[\ineqidx]\paren*{\ptinriter}} + \frac{1}{\pullback[\ptinriter]{\ineqfun[\ineqidx]}\paren*{\tmethr\dirptinriter}}\abs*{\ineqLagmultinriter[\ineqidx] - \frac{\barrparam[]}{\ineqfun[\ineqidx]\paren*{\ptinriter}}}+ \frac{\barrparam[]\abs*{\pullback[\ptinriter]{\ineqfun[\ineqidx]}\paren*{\tmethr\dirptinriter} - \ineqfun[\ineqidx]\paren*{\ptinriter}}}{\pullback[\ptinriter]{\ineqfun[\ineqidx]}\paren*{\tmethr\dirptinriter}^{2}\ineqfun[\ineqidx]\paren*{\ptinriter}}}\\
                    &\qquad\Riemnorm[\ptinriter]{\gradstr\ineqfun[\ineqidx]\paren*{\ptinriter}}^{2} \Riemnorm[\ptinriter]{\dirptinriter}^{2} + \frac{\barrparam[]\paren*{\betaRL[{\ineqfun[\ineqidx]}]}^{2}}{\epsthr^{2}} \Riemnorm[\ptinriter]{\dirptinriter}^{4} + \frac{2\barrparam[]\betaRL[{\ineqfun[\ineqidx]}]}{\epsthr^{2}} \Riemnorm[\ptinriter]{\gradstr\ineqfun[\ineqidx]\paren*{\ptinriter}}\Riemnorm[\ptinriter]{\dirptinriter}^{3}\\
                    &\quad\leq \paren*{\frac{2\coefffou}{27\ineqdime\tholdvalthr^{2}} + \frac{2\coefffou}{27\ineqdime\tholdvalthr^{2}} + \frac{2\coefffou}{27\ineqdime\tholdvalthr^{2}}}\Riemnorm[\ptinriter]{\gradstr\ineqfun[\ineqidx]\paren*{\ptinriter}}^{2} \Riemnorm[\ptinriter]{\dirptinriter}^{2}\\
                    &\quad + \frac{\barrparam[]\paren*{\betaRL[{\ineqfun[\ineqidx]}]}^{2}}{\epsthr^{2}} \paren*{\thldtrradiussix}^{2}\Riemnorm[\ptinriter]{\dirptinriter}^{2} + \frac{2\barrparam[]\betaRL[{\ineqfun[\ineqidx]}]\Riemnorm[\ptinriter]{\gradstr\ineqfun[\ineqidx]\paren*{\ptinriter}}}{\epsthr^{2}}\thldtrradiussix\Riemnorm[\ptinriter]{\dirptinriter}^{2}\\
                    &\quad\leq \paren*{\frac{2\coefffou}{27\ineqdime} + \frac{2\coefffou}{27\ineqdime} + \frac{2\coefffou}{27\ineqdime}}\Riemnorm[\ptinriter]{\dirptinriter}^{2} + \frac{2\coefffou}{9\ineqdime}\Riemnorm[\ptinriter]{\dirptinriter}^{2} + \frac{2\coefffou}{9\ineqdime}\Riemnorm[\ptinriter]{\dirptinriter}^{2} = \frac{2\coefffou}{3\ineqdime}\Riemnorm[\ptinriter]{\dirptinriter}^{2},
                \end{aligned}
            \end{equation}}{\begin{align}
                &\label{ineq:coefffoubounddiffsquaregradineqfun}\\
                &\abs*{\frac{\ineqLagmultinriter[\ineqidx]}{\ineqfun[\ineqidx]\paren*{\ptinriter}}\metr[\ptinriter]{\gradstr\ineqfun[\ineqidx]\paren*{\ptinriter}}{\dirptinriter}^{2} - \frac{\barrparam[]}{{\pullback[\ptinriter]{\ineqfun[\ineqidx]}\paren*{\tmethr\dirptinriter}}^{2}} \metr[\ptinriter]{\gradstr\pullback[\ptinriter]{\ineqfun[\ineqidx]}\paren*{\tmethr\dirptinriter}}{\dirptinriter}^{2}} \leq \abs*{\frac{\ineqLagmultinriter[\ineqidx]}{\ineqfun[\ineqidx]\paren*{\ptinriter}} \\
                &- \frac{\barrparam[]}{{\pullback[\ptinriter]{\ineqfun[\ineqidx]}\paren*{\tmethr\dirptinriter}}^{2}}} \metr[\ptinriter]{\gradstr\ineqfun[\ineqidx]\paren*{\ptinriter}}{\dirptinriter}^{2} + \frac{\barrparam[]}{{\pullback[\ptinriter]{\ineqfun[\ineqidx]}\paren*{\tmethr\dirptinriter}}^{2}} \abs*{
                \ang*{\gradstr\pullback[\ptinriter]{\ineqfun[\ineqidx]}\paren*{\zerovec[\ptinriter]} - \gradstr\pullback[\ptinriter]{\ineqfun[\ineqidx]}\paren*{\tmethr\dirptinriter},\\ &\dirptinriter}_{\ptinriter}} \paren*{\abs*{\metr[\ptinriter]{\gradstr\pullback[\ptinriter]{\ineqfun[\ineqidx]}\paren*{\tmethr\dirptinriter} - \gradstr\pullback[\ptinriter]{\ineqfun[\ineqidx]}\paren*{\zerovec[\ptinriter]}}{\dirptinriter}} + 2\abs*{\metr[\ptinriter]{\gradstr\ineqfun[\ineqidx]\paren*{\ptinriter}}{\dirptinriter}}}\\
                &\leq \paren*{\frac{\ineqLagmultinriter[\ineqidx]\abs*{\pullback[\ptinriter]{\ineqfun[\ineqidx]}\paren*{\tmethr\dirptinriter} - \ineqfun[\ineqidx]\paren*{\ptinriter}}}{\pullback[\ptinriter]{\ineqfun[\ineqidx]}\paren*{\tmethr\dirptinriter}\ineqfun[\ineqidx]\paren*{\ptinriter}} + \frac{1}{\pullback[\ptinriter]{\ineqfun[\ineqidx]}\paren*{\tmethr\dirptinriter}}\abs*{\ineqLagmultinriter[\ineqidx] - \frac{\barrparam[]}{\ineqfun[\ineqidx]\paren*{\ptinriter}}}+ \frac{\barrparam[]\abs*{\pullback[\ptinriter]{\ineqfun[\ineqidx]}\paren*{\tmethr\dirptinriter} - \ineqfun[\ineqidx]\paren*{\ptinriter}}}{\pullback[\ptinriter]{\ineqfun[\ineqidx]}\paren*{\tmethr\dirptinriter}^{2}\ineqfun[\ineqidx]\paren*{\ptinriter}}}\\
                &\Riemnorm[\ptinriter]{\gradstr\ineqfun[\ineqidx]\paren*{\ptinriter}}^{2} \Riemnorm[\ptinriter]{\dirptinriter}^{2} + \frac{\barrparam[]\paren*{\betaRL[{\ineqfun[\ineqidx]}]}^{2}}{\epsthr^{2}} \Riemnorm[\ptinriter]{\dirptinriter}^{4} + \frac{2\barrparam[]\betaRL[{\ineqfun[\ineqidx]}]}{\epsthr^{2}} \Riemnorm[\ptinriter]{\gradstr\ineqfun[\ineqidx]\paren*{\ptinriter}}\Riemnorm[\ptinriter]{\dirptinriter}^{3}\\
                &\leq \paren*{\frac{2\coefffou}{27\ineqdime\tholdvalthr^{2}} + \frac{2\coefffou}{27\ineqdime\tholdvalthr^{2}} + \frac{2\coefffou}{27\ineqdime\tholdvalthr^{2}}}\Riemnorm[\ptinriter]{\gradstr\ineqfun[\ineqidx]\paren*{\ptinriter}}^{2} \Riemnorm[\ptinriter]{\dirptinriter}^{2} + \frac{\barrparam[]\paren*{\betaRL[{\ineqfun[\ineqidx]}]}^{2}}{\epsthr^{2}} \paren*{\thldtrradiussix}^{2}\Riemnorm[\ptinriter]{\dirptinriter}^{2}\\
                &+ \frac{2\barrparam[]\betaRL[{\ineqfun[\ineqidx]}]\Riemnorm[\ptinriter]{\gradstr\ineqfun[\ineqidx]\paren*{\ptinriter}}}{\epsthr^{2}}\thldtrradiussix\Riemnorm[\ptinriter]{\dirptinriter}^{2} \leq \paren*{\frac{2\coefffou}{27\ineqdime} + \frac{2\coefffou}{27\ineqdime} + \frac{2\coefffou}{27\ineqdime}}\Riemnorm[\ptinriter]{\dirptinriter}^{2}\\[-1.5em]
                &+ \frac{2\coefffou}{9\ineqdime}\Riemnorm[\ptinriter]{\dirptinriter}^{2} + \frac{2\coefffou}{9\ineqdime}\Riemnorm[\ptinriter]{\dirptinriter}^{2} = \frac{2\coefffou}{3\ineqdime}\Riemnorm[\ptinriter]{\dirptinriter}^{2},
            \end{align}}
            where the first inequality follows from \cref{eq:pullbackgradzero} and $\barrparam[] > 0$, the second one from \cref{def:radialLCone} with $\tmethr \leq 1$, $\ineqLagmultinriter[\ineqidx], \barrparam[] \in \setRpp[]$, and $\ineqfun[\ineqidx]\paren*{\ptinriter} > \epsthr > 0$ and $\pullback[\ptinriter]{\ineqfun[\ineqidx]}\paren*{\tmethr\dirptinriter} > \epsthr > 0$ by \enumicref{lemm:ineqfunbound}{lemm:ienqfununiflowbd}, the third one from \cref{ineq:pullbackineqfuncoefffoubound,ineq:gradhessineqLagmultbound,ineq:ineqLagmultcoefffoubound}, \enumicref{lemm:ineqfunbound}{lemm:ienqfununiflowbd} again, and $\Riemnorm[\ptinriter]{\dirptinriter} \leq \thldtrradiussix$, and the last one \isextendedversion{from \cref{ineq:gradhessineqLagmultbound,eq:thldtrradiussixdef}.}{from \cref{ineq:gradhessineqLagmultbound} and the definition of $\thldtrradiussix$.}
            By combining \cref{ineq:diffpredared} with \cref{ineq:boundhessobjfun,ineq:coefffouboundhessineqfun,ineq:coefffoubounddiffsquaregradineqfun}, we obtain
            \isextendedversion{
            \begin{align}
                &\abs*{\predinriter-\aredinriter}\\
                &\leq \int_{0}^{1} \paren*{1-\tmethr} \paren*{\betaRLtwo[\objfun] \thldtrradiussix \Riemnorm[\ptinriter]{\dirptinriter}^{2} + \sum_{\ineqidx\in\ineqset} \frac{2\coefffou}{3\ineqdime}\Riemnorm[\ptinriter]{\dirptinriter}^{2} + \sum_{\ineqidx\in\ineqset}\frac{2\coefffou}{3\ineqdime}\Riemnorm[\ptinriter]{\dirptinriter}^{2}} \dd{\tmethr}\\
                &\leq \frac{1}{2}\paren*{\frac{2\coefffou}{3} + \frac{2\coefffou}{3} + \frac{2\coefffou}{3}}\Riemnorm[\ptinriter]{\dirptinriter}^{2} = \coefffou\Riemnorm[\ptinriter]{\dirptinriter}^{2},
            \end{align}}{
            {\small \begin{align}
                &\abs*{\predinriter-\aredinriter} \leq \int_{0}^{1} \paren*{1-\tmethr} \paren*{\betaRLtwo[\objfun] \thldtrradiussix \Riemnorm[\ptinriter]{\dirptinriter}^{2} + \sum_{\ineqidx\in\ineqset} \frac{2\coefffou}{3\ineqdime}\Riemnorm[\ptinriter]{\dirptinriter}^{2} + \sum_{\ineqidx\in\ineqset}\frac{2\coefffou}{3\ineqdime}\Riemnorm[\ptinriter]{\dirptinriter}^{2}} \dd{\tmethr}\\[-1.5em]
                &\leq \frac{1}{2}\paren*{\frac{2\coefffou}{3} + \frac{2\coefffou}{3} + \frac{2\coefffou}{3}}\Riemnorm[\ptinriter]{\dirptinriter}^{2} = \coefffou\Riemnorm[\ptinriter]{\dirptinriter}^{2},
            \end{align}}
            }
            where the first inequality holds by $\Riemnorm[\ptinriter]{\dirptinriter} \leq \thldtrradiussix$ and the second one 
            \isextendedversion{
            from \cref{eq:thldtrradiussixdef}.}{from the definition of $\thldtrradiussix$.}
            The proof is complete.
        \end{proof}

    \isextendedversion{
    
    \section{Search direction}\label{sec:searchdir}
        In this appendix, we present three search directions obtained from the subproblem \cref{prob:TRS}, namely, the \Cauchypoint{}, the \eigenstep{}, and the \exactsolution{}.
        We first define the \Cauchypoint{} as follows:
        \begin{definition}[{\cite[Equation (7.8)]{Absiletal08OptimAlgoonMatMani}}]
            The \Cauchypoint{} $\dirCauchy[\pt]\in\tanspc[\pt]\mani$ of the subproblem \cref{prob:TRS}
            is defined as
            \begin{align}
                \dirCauchy[\pt] \coloneqq - \coeffCauchyoptim \trslin[{\barrparam[]}]\paren*{\pt},
            \end{align}
            where we define $\coeffCauchyoptim=0$ if $\trslin[{\barrparam[]}]\paren*{\pt} = \zerovec[\pt]$, otherwise,
            \begin{align}
                \coeffCauchyoptim &\coloneqq \argmin[{\coeffCauchy \geq 0}]\brc*{\modelfun[{\allvar,\barrparam[]}]\paren*{-\coeffCauchy\trslin[{\barrparam[]}]\paren*{\pt}} \text{ subject to } \Riemnorm[\pt]{\coeffCauchy\trslin[{\barrparam[]}]\paren*{\pt}} \leq \trradius[]}\\
                &= 
                \begin{cases}
                    \frac{\trradius[]}{\Riemnorm[\pt]{\trslin[{\barrparam[]}]\paren*{\pt}}} & \text{if }
                    \metr[\pt]{\trsquad\paren*{\allvar}\sbra*{\trslin[{\barrparam[]}]\paren*{\pt}}}{\trslin[{\barrparam[]}]\paren*{\pt}} \leq 0, \\
                    \min\paren*{\frac{\Riemnorm[\pt]{\trslin[{\barrparam[]}]\paren*{\pt}}^{2}}{\metr[\pt]{\trsquad\paren*{\allvar}\sbra*{\trslin[{\barrparam[]}]\paren*{\pt}}}{\trslin[{\barrparam[]}]\paren*{\pt}}}, \frac{\trradius[]}{\Riemnorm[\pt]{\trslin[{\barrparam[]}]\paren*{\pt}}}} & \text{otherwise}.\\
                \end{cases}
            \end{align}
        \end{definition}
        The bound on the predicted reduction from the \Cauchypoint{} is shown below.
        \begin{lemma}[{\cite[Equation (7.14)]{Absiletal08OptimAlgoonMatMani}}]\label{lemm:Cauchydecreasing}
            The following holds:
            \begin{align}
                \modelfun[{\allvar, \barrparam[]}]\paren*{\zerovec[\pt]} - \modelfun[{\allvar, \barrparam[]}]\paren*{\dirCauchy[\pt]} \geq \frac{1}{2}\Riemnorm[\pt]{\trslin[{\barrparam[]}]\paren*{\pt}}\min\paren*{\trradius[], \frac{\Riemnorm[\pt]{\trslin[{\barrparam[]}]\paren*{\pt}}}{\opnorm{\trsquad\paren*{\allvar}}}}.
            \end{align}
            If $\opnorm{\trsquad\paren*{\allvar}}=0$, then we regard $\min\paren*{\trradius[], \frac{\Riemnorm[\pt]{\trslin[{\barrparam[]}]\paren*{\pt}}}{\opnorm{\trsquad\paren*{\allvar}}}} = \trradius[]$.
        \end{lemma}
        In \cref{subsec:globconvinner}, we prove that, if a search direction is adopted such that the decrease in the objective function of \cref{prob:TRS} is not less than that achieved by the \Cauchypoint{}, \RIPTRM{} has a global convergence property to an \AKKT{} point under certain assumptions.
        In the experiments described in \cref{sec:experiment,sec:additionalexperiment}, we use the search direction computed by the \tCG{} method, which satisfies the aforementioned condition.
        
        Next, we define the \eigenstep{}.
        Let $\mineigval\sbra*{\trsquad\paren*{\allvar}}\in\setR[]$ be the minimum eigenvalue of $\trsquad\paren*{\allvar}$.
        \begin{definition}[{\cite[Lemma 6.16]{Boumal23IntroOptimSmthMani}}]
            Given $\ineqLagmult[]\in\setRpp[\ineqdime]$ and $\trradius[], \barrparam[] \in \setRpp[]$, let $\diruniteigen[\pt]\in\tanspc[\pt]\mani$ satisfy
            \begin{align}
                \Riemnorm[\pt]{\diruniteigen[\pt]} = 1, \, \metr[\pt]{\trslin[{\barrparam[]}]\paren*{\pt}}{\diruniteigen[\pt]} \leq 0, \text{ and } \metr[\pt]{\trsquad\paren*{\allvar}\sbra*{\diruniteigen[\pt]}}{\diruniteigen[\pt]} < -\epseigen,
            \end{align}
            where $\epseigen\in\setR[]$ is a predefined parameter satisfying $\mineigval\sbra*{\trsquad\paren*{\allvar}}< -\epseigen$.
            We call $\direigen[\pt]\coloneqq\trradius[]\diruniteigen[\pt]$ an \eigenstep{}.
        \end{definition}
        Note that an eigenvector corresponding to $\mineigval\sbra*{\trsquad\paren*{\allvar}}\in\setR[]$ satisfies the conditions on the \eigenstep{} with the appropriate choice of its sign and the scaling.
        The bound on the predicted reduction is also known when using the \eigenstep{}:
        \begin{lemma}[{\cite[Lemma 6.16]{Boumal23IntroOptimSmthMani}}]\label{lemm:Eigendecreasing}
            For any $\epseigen\in\setR[]$, any point $\allvar=\paren*{\pt,\ineqLagmult[]}\in\strictfeasirgn\times\setRpp[\ineqdime]$ with $\mineigval\sbra*{\trsquad\paren*{\allvar}} < -\epseigen$, and $\barrparam[] > 0$, the corresponding \eigenstep{} $\direigen[\pt]\in\tanspc[\pt]\mani$ satisfies that
            \begin{align}\label{eq:boundeigenstep}
                \modelfun[{\allvar, \barrparam[]}]\paren*{\zerovec[\pt]} - \modelfun[{\allvar, \barrparam[]}]\paren*{\direigen[\pt]} \geq \frac{1}{2}\paren*{\trradius[]}^{2} \epseigen.
            \end{align}
        \end{lemma}
        We will prove that, if we adopt a search direction whose decrease in the objective function of \cref{prob:TRS} is not less than those of the \Cauchypoint{} and the \eigenstep{}, \RIPTRM{} has a global convergence property to an \SOSP{} under assumptions.
        The \eigenstep{}, however, is rarely computed in practice as mentioned in \cite[Section~6.4.2]{Boumal23IntroOptimSmthMani}.
        In this paper, the \eigenstep{} mainly serves to show that the \exactsolution{} of the subproblem~\cref{prob:TRS} satisfies \cref{eq:boundeigenstep}.
        
        Let us end the subsection by considering the \exactsolution{}, that is, a global optimum of the subproblem~\cref{prob:TRS}.
        Note that the subproblem~\cref{prob:TRS} has a global optimum since the feasible region is bounded and closed, and the objective function is continuous.
        We provide the necessary and sufficient condition for the global optimality of \cref{prob:TRS}.
        We discuss the computation of $\direxactinriter$ in \cref{subsec:experimentenv}.
        \begin{proposition}[{\cite[Proposition 7.3.1]{Absiletal08OptimAlgoonMatMani}}]\label{prop:trsgloboptimiffcond}
            The vector $\direxact[\pt]\in\tanspc[\pt]\mani$ is a global optimum of \cref{prob:TRS}
            if and only if there exists a scalar $\trsLagmult \geq 0$ such that the following hold:
            \begin{subequations}\label{eq:trsgloboptimiffcond}
                \begin{align}
                    \paren*{\trsquad\paren*{\allvar}+\trsLagmult\id[{\tanspc[\pt]\mani}]}\direxact[\pt] &= - \trslin[{\barrparam[]}]\paren*{\pt},\label{eq:trsgloboptimiffcondquadlineq}\\
                    \trsLagmult\paren*{\trradius[] - \Riemnorm[\pt]{\direxact[\pt]}} &= 0,\label{eq:trsgloboptimiffcondconmpl}\\
                    \Riemnorm[\pt]{\direxact[\pt]} &\leq \trradius[],\label{eq:trsgloboptimiffconddirtrradius}\\
                    \trsquad\paren*{\allvar}+\trsLagmult\id[{\tanspc[\pt]\mani}] &\succeq 0.\label{eq:trsgloboptimiffcondquadpsd}
                \end{align}
            \end{subequations}
        \end{proposition}
        In a companion paper~\cite{Obaraetal2025LocalConvofRiemIPMs}, we prove that, if we adopt the exact solution of the subproblem \cref{prob:TRS} as the search direction, \RIPTRM{} also possesses a local near-quadratic convergence property to a solution of \RICO{}~\cref{prob:RICO} under assumptions.

    \section{Sufficient conditions for assumptions}
    \label{appx:suffcondassus}
        In this appendix, we provide detailed discussions of the sufficient conditions for the assumptions
        in \cref{subsec:globconvouter,subsec:globconvinner}.
    
    \subsection{Sufficient conditions for assumptions in \texorpdfstring{\cref{subsec:globconvinner}}{global convergence of inner iterations}}
        In this subsection, we provide the sufficient conditions for \cref{assu:radialLCone,assu:radialLCtwo,assu:pullbackineqfuncont,assu:retrdistbound}.

        We first consider the sufficient conditions for \cref{assu:radialLCone,assu:radialLCtwo}.
        To this end, we derive the following lemma for the \rLCone{} and \rLCtwo{} properties:
        \begin{lemma}\label{lemm:suffcondradialLConetwo}
            Let $\subsetmani$ be a compact subset of $\mani$ and $\funtwo\colon\mani\to\setR$.
            If $\funtwo$ is of class $C^{2}$, then $\funtwo$ is \rLCone{} on $\subsetmani$.
            Additionally, if $\funtwo$ is of class $C^{3}$, then $\funtwo$ is \rLCtwo{} on $\subsetmani$.
        \end{lemma}
        \begin{proof}
            Let $\deltaRL[\funtwo] > 0$ be any scalar.
            If $\funtwo$ is of class $C^{2}$, it follows from \cite[Lemma~10.57]{Boumal23IntroOptimSmthMani} that there exists $\Lipschitzconsttwo[1] > 0$ such that, for any $\ptone\in\subsetmani$ and all $\tanvectwo[\ptone]\in\tanspc[\ptone]\mani$ with $\Riemnorm[\ptone]{\tanvectwo[\ptone]}\leq \deltaRL[\funtwo]$, $\Riemnorm[\ptone]{\gradstr\pullback[\ptone]{\funtwo}\paren*{\tanvectwo[\ptone]} - \gradstr\pullback[\ptone]{\funtwo}\paren*{\zerovec[\ptone]}} \leq \Lipschitzconsttwo[1]\Riemnorm[\ptone]{\tanvectwo[\ptone]}$ holds.
            Therefore, for any $\ptone\in\subsetmani$ and all $\tmetwo \geq 0, \tanvecone[\ptone]\in\tanspc[\ptone]\mani$ with $\tmetwo \Riemnorm[\ptone]{\tanvecone[\ptone]} \leq \deltaRL[\funtwo]$, it holds that
            \begin{align}
                &\abs*{\metr[\ptone]{\gradstr\pullback[\ptone]{\funtwo}\paren*{\tmetwo\tanvecone[\ptone]} - \gradstr\pullback[\ptone]{\funtwo}\paren*{\zerovec[\ptone]}}{\tanvecone[\ptone]}} \\
                &\quad\leq \Riemnorm[\ptone]{\gradstr\pullback[\ptone]{\funtwo}\paren*{\tmetwo\tanvecone[\ptone]} - \gradstr\pullback[\ptone]{\funtwo}\paren*{\zerovec[\ptone]}} \Riemnorm[\ptone]{\tanvecone[\ptone]} \leq \Lipschitzconsttwo[1]\tmetwo\Riemnorm[\ptone]{\tanvecone[\ptone]}^{2}.
            \end{align}
            By setting $\betaRL[\funtwo]=\Lipschitzconsttwo[1]$, we complete the proof of the sufficient condition for \rLCone{} property.
            
            Similarly, if $\funtwo$ is of class $C^{3}$, it follows again from \cite[Lemma~10.57]{Boumal23IntroOptimSmthMani} that there exists $\Lipschitzconsttwo[2] > 0$ such that, for any $\ptone\in\subsetmani$ and all $\tanvectwo[\ptone]\in\tanspc[\ptone]\mani$ with $\Riemnorm[\ptone]{\tanvectwo[\ptone]}\leq \deltaRLtwo[\funtwo]$, $\opnorm{\Hess\pullback[\ptone]{\funtwo}\paren*{\tanvectwo[\ptone]} - \Hess\pullback[\ptone]{\funtwo}\paren*{\zerovec[\ptone]}} \leq \Lipschitzconsttwo[2]\Riemnorm[\ptone]{\tanvectwo[\ptone]}$ holds.
            Thus, for any $\ptone\in\subsetmani$ and all $\tmetwo \geq 0, \tanvecone[\ptone]\in\tanspc[\ptone]\mani$ with $\tmetwo \Riemnorm[\ptone]{\tanvecone[\ptone]} \leq \deltaRLtwo[\funtwo]$, it holds that
            \begin{align}
                &\abs*{\metr[\ptone]{\paren*{\Hess\pullback[\ptone]{\funtwo}\paren*{\tmetwo\tanvecone[\ptone]} - \Hess\pullback[\ptone]{\funtwo}\paren*{\zerovec[\ptone]}}\sbra*{\tanvecone[\ptone]} }{\tanvecone[\ptone]}} \\
                &\quad\leq \opnorm{\Hess\pullback[\ptone]{\funtwo}\paren*{\tmetwo\tanvecone[\ptone]} - \Hess\pullback[\ptone]{\funtwo}\paren*{\zerovec[\ptone]}} \Riemnorm[\ptone]{\tanvecone[\ptone]}^{2} \leq \Lipschitzconsttwo[2]\tmetwo\Riemnorm[\ptone]{\tanvecone[\ptone]}^{3}.
            \end{align}
            By setting $\betaRLtwo[\funtwo]=\Lipschitzconsttwo[2]$, we complete the proof of the sufficient condition for the \rLCtwo{} property.
        \end{proof}
        Using the lemma, we formally derive the following corollaries:
        \begin{corollary}
            \cref{assu:radialLCone} is fulfilled if the generated sequence $\brc*{\ptinriter}_{\inriteridx}$ is bounded.
        \end{corollary}
        \begin{corollary}
            \cref{assu:radialLCtwo} holds if the generated sequence $\brc*{\ptinriter}_{\inriteridx}$ is \\ bounded and all functions $\objfun, \brc*{\ineqfun[\ineqidx]}_{\ineqidx\in\ineqset}$ are of class $C^3$.
        \end{corollary}
        As for \cref{assu:pullbackineqfuncont}, we provide sufficient conditions as follows:
        \begin{lemma}\label{lemm:suffcondpullbackineqfuncont}
            The following hold:
            \begin{enumerate}
                \item Suppose $\retr[]=\expmap[]$ and that, for every $\ineqidx\in\ineqset$, the function $\ineqfun[\ineqidx]\colon\mani\to\setR[]$ is $\Lipschitzconstone[\ineqidx]$-Lipschitz continuous; that is, there exists $\Lipschitzconstone[\ineqidx] > 0$ such that
                \begin{align}
                    \abs*{\ineqfun[\ineqidx]\paren*{\pt[1]} - \ineqfun[\ineqidx]\paren*{\pt[2]}} \leq \Lipschitzconstone[\ineqidx] \Riemdist{\pt[1]}{\pt[2]} \text{ for all } \pt[1], \pt[2] \in \mani.
                \end{align}
                Then, \cref{assu:pullbackineqfuncont} holds.\label{lemm:suffcondLipschitzcont} 
                \item If $\brc*{\ptinriter}_{\inriteridx}$ is bounded, then \cref{assu:pullbackineqfuncont} holds.\label{lemm:suffcondbounded} 
            \end{enumerate}
        \end{lemma}
        \begin{proof}
            We first prove \cref{lemm:suffcondLipschitzcont}.
            Notice that, due to the completeness of $\mani$, the domain of the exponential map at $\pt$ is the whole of $\tanspc[\pt]\mani$ for all $\pt\in\mani$.
            It follows from \cite[Proposition 10.41]{Boumal23IntroOptimSmthMani} that $\Lipschitzconstone[\ineqidx]$-Lipschitz continuity is equivalent to $\abs*{\ineqfun[\ineqidx]\paren*{\expmap[\pt]\paren*{\tanvecone[\pt]}} - \ineqfun[\ineqidx]\paren*{\pt}} \leq \Lipschitzconstone[\ineqidx] \Riemnorm[\pt]{\tanvecone[\pt]}$ for all $\pt\in\mani$ and all $\tanvecone \in \tanspc[\pt]\mani$.
            Therefore, for any $\epsone > 0$, we attain the conclusion by setting $\thldtrradiusthr\coloneqq\min_{\ineqidx\in\ineqset} \frac{\epsone}{\Lipschitzconstone[\ineqidx]}$.

            Next, we consider \cref{lemm:suffcondbounded}.
            We regard the tangent bundle $\tanspc[]\mani$ as a Riemannian manifold endowed with the Riemannian distance 
            \begin{align}\label{eq:Riemdisttanbddef}
                \Riemdist[{\tanspc[]\mani}]{\paren*{\pt[1], \tanvecone[{\pt[1]}]}}{\paren*{\pt[2], \tanvecone[{\pt[2]}]}}\coloneqq \inf_{\pwisecurve\in\pwisecurveset{\pt[1]}{\pt[2]}} \brc*{\sqrt{{\length\paren*{\pwisecurve}}^{2} + \Riemnorm[{\pt[1]}]{\partxp[\pwisecurve]{\pt[1]}{\pt[2]}\sbra*{\tanvecone[{\pt[2]}]} - \tanvecone[{\pt[1]}]}^{2}}}
            \end{align}
            for any $\paren*{\pt[1], \tanvecone[{\pt[1]}]}, \paren*{\pt[2], \tanvecone[{\pt[2]}]} \in \tanspc[]\mani$, where $\pwisecurveset{\pt[1]}{\pt[2]}$ denotes the set of all piecewise regular curve segments in $\mani$ joining $\pt[2]$ to $\pt[1]$, $\length\paren*{\pwisecurve}$ is the length of $\pwisecurve\in\pwisecurveset{\pt[1]}{\pt[2]}$, and $\partxp[\pwisecurve]{\pt[1]}{\pt[2]}$ is the parallel transport from $\tanspc[{\pt[2]}]\mani$ to $\tanspc[{\pt[1]}]\mani$ along $\pwisecurve$; see \cite[pp.33—34,108]{Lee18IntrotoRiemManibook2ndedn} and \cite[Section 10.1]{Boumal23IntroOptimSmthMani} for the topics related to curve segments and \cite[Section 2]{deOliveiraFerreira2020NewrtonMethforFindSingofSplClofLocLipSchitzContVecFldonRiemMani} and \cite[Appendix II.A.2]{CanaryEpsteinMarden2006FundHyperbolicMani} for the Riemannian distance on tangent bundles.
            Let $\compactsubset\subseteq\mani$ be a compact subset with $\brc*{\ptinriter}_{\inriteridx}\subseteq\compactsubset$ and $\compactsubsettanbd\coloneqq\brc*{\paren*{\pt, \tanvecone[\pt]}\in\tanspc[]\mani\colon \pt\in\compactsubset \text{ and } \Riemnorm[\pt]{\tanvecone[\pt]}\leq\thldval}$ with some $\thldval > 0$.
            Note that $\compactsubsettanbd$ is also compact~\cite[Exercise 10.31]{Boumal23IntroOptimSmthMani}.
            For each $\ineqidx\in\ineqset$, we consider a composite function $\ineqfun[\ineqidx]\circ\retr[]\colon\tanspc[]\mani\to\setR[]$.
            From the continuities of $\ineqfun[\ineqidx]$ and $\retr[]$, the restriction of $\ineqfun[\ineqidx]\circ\retr[]$ to $\compactsubsettanbd$ is uniformly continuous by the Heine-Cantor theorem.
            Namely, for any $\epsone > 0$, there exists $\thldtrradiusthr > 0$ such that $\thldtrradiusthr \leq \thldval$ and, for all $\paren*{\pt[1], \tanvecone[{\pt[1]}]}, \paren*{\pt[2], \tanvecone[{\pt[2]}]} \in \compactsubsettanbd$ with $\Riemdist[{\tanspc[]\mani}]{\paren*{\pt[1], \tanvecone[{\pt[1]}]}}{\paren*{\pt[2], \tanvecone[{\pt[2]}]}}\leq\thldtrradiusthr$, $\abs*{\ineqfun[\ineqidx]\circ\retr[]\paren*{\pt[1], \tanvecone[1]} -  \ineqfun[\ineqidx]\circ\retr[]\paren*{\pt[2], \tanvecone[2]}}\leq\epsone$ holds.
            Therefore, for all $\inriteridx\in\setNz$ and any $\tanvecone[\ptinriter]\in\tanspc[\ptinriter]\mani$ with $\Riemnorm[\ptinriter]{\tanvecone[\ptinriter]} \leq \thldtrradiusthr$, we conclude that the statement is true by substituting $\paren*{\ptinriter, \tanvecone[\ptinriter]},\paren*{\ptinriter, \zerovec[\ptinriter]}$ for $\paren*{\pt[1], \tanvecone[{\pt[1]}]}, \paren*{\pt[2], \tanvecone[{\pt[2]}]}$, together with $\Riemdist[{\tanspc[]\mani}]{\paren*{\ptinriter, \tanvecone[\ptinriter]}}{\paren*{\ptinriter, \zerovec[\ptinriter]}} = \Riemnorm[\ptinriter]{\tanvecone[\ptinriter]}$. 
        \end{proof}        
        We also provide sufficient conditions for \cref{assu:retrdistbound} in the following lemma:
        \begin{lemma}\label{lemm:retrdistboundsuffcond}
            The following hold:
            \begin{enumerate}
                \item If $\brc*{\ptinriter}_{\inriteridx}$ is bounded, then \cref{assu:retrdistbound} holds.\label{lemm:retrdistboundsuffcondseqbound}
                \item If $\mani$ is compact, then \cref{assu:retrdistbound} holds.\label{lemm:retrdistboundsuffcondcompact}
            \end{enumerate}
        \end{lemma}
        \begin{proof}
            Under \cref{lemm:retrdistboundsuffcondseqbound} or \cref{lemm:retrdistboundsuffcondcompact}, let $\compactsubset\subseteq\mani$ be a compact subset with $\brc*{\ptinriter}_{\inriteridx}\subseteq\compactsubset$.
            For any $\epstwo > 0$, define $\compactsubsettanbd\coloneqq\brc*{\paren*{\pt, \tanvecone[\pt]}\in\tanspc[]\mani\colon \pt\in\compactsubset \text{ and } \Riemnorm[\pt]{\tanvecone[\pt]}\leq\epstwo}$.
            Note that $\compactsubsettanbd$ is compact~\cite[Exercise 10.31]{Boumal23IntroOptimSmthMani}.
            Thus, from the smoothness of $\retr[]$ and the continuity of the operator norm, there exists $\deltatwo > 0$ such that $\opnorm{\D\retr[\pt]\paren*{\tanvecone[\pt]}} \leq \deltatwo$ for all $\paren*{\pt,\tanvecone[\pt]}\in\compactsubsettanbd$.
            For any $\paren*{\pt,\tanvecone[\pt]}\in\compactsubsettanbd$, consider the curve $\curve\paren*{\tmefiv}=\retr[\ptone]\paren*{\tmefiv\tanvecone[\ptone]}$ and let $\length\paren*{\curve} \coloneqq \int_{0}^{1} \Riemnorm[\curve\paren*{\tmefiv}]{\curve[\prime]\paren*{\tmefiv}} \dd{\tmefiv}$ be the length of the curve $\curve$ on the interval $\sbra*{0, 1}$.
            Then, we have
            \begin{align}
                &\Riemdist{\pt}{\retr[\pt]\paren*{\tanvecone[\pt]}} \leq \length\paren*{\curve}  = \int_{0}^{1} \Riemnorm[\curve\paren*{\tmefiv}]{\D\retr[\pt]\paren*{\tmefiv\tanvecone[\pt]}\sbra*{\tanvecone[\pt]}} \dd{\tmefiv}\\
                &\leq \int_{0}^{1} \opnorm{\D\retr[\pt]\paren*{\tmefiv\tanvecone[\pt]}}\Riemnorm[\pt]{\tanvecone[\pt]} \dd{\tmefiv} \leq \int_{0}^{1} \deltatwo\Riemnorm[\pt]{\tanvecone[\pt]} \dd{\tmefiv} = \deltatwo\Riemnorm[\pt]{\tanvecone[\pt]}.
            \end{align}
            The proof is complete.
        \end{proof}
        As in \cref{lemm:Eigendecreasing}, the \eigenstep{} satisfies \cref{assu:Eigendecreasing} with $\constEigen=\frac{1}{2}$, and so does the \exactsolution{}~\cref{eq:trsgloboptimiffcond}.

        Let us end this subsection by noting that all the assumptions in \cref{subsec:globconvinner} are standard in the literature; for example, \cref{assu:retrsecondord} is not restrictive as in \cref{subsec:notationRiemggeo}.
        \cref{assu:radialLCone} is made in \cite[Section~7.4.1]{Absiletal08OptimAlgoonMatMani}.
        The smoothness of the functions and the boundedness of the generated sequences, which are stronger than 
        \cref{assu:radialLCone,assu:radialLCtwo,assu:ineqfunbounded,assu:objfunbounded,assu:pullbackineqfuncont,assu:ineqLagmultbounded,assu:retrdistbound},
        are assumed in \cite[Assumptions~(C1), (C2)]{LaiYoshise2024RiemIntPtMethforCstrOptimonMani}.
        The lower boundedness of $\objfun$, which is a sufficient condition of \assuenumicref{assu:objfunbounded}{assu:objfunboundedbelow}, is made in \cite[Assumption~A.1]{Agarwaletal2021AdaptRegwithCubiconMani} and \cite[A6.5]{Boumal23IntroOptimSmthMani}.
        Assumptions similar to \cref{assu:ineqLagmultbounded,assu:ineqLagmultconverged} are made in \cite[AS.6, AS.10]{Connetal2000PrimalDualTRAlgoforNonconvexNLP}.
        \cref{assu:retrdistbound} is made in \cite[Equation~(7.25)]{Absiletal08OptimAlgoonMatMani}.
        \cref{assu:Eigendecreasing} is made in \cite[A6.4]{Boumal23IntroOptimSmthMani}.

    \subsection{Sufficient condition for assumption in \texorpdfstring{\cref{subsec:globconvouter}}{global convergence of outer iteration}}
        In this subsection, we provide the sufficient condition for \cref{assu:LipschitzHess}.
        To this end, we first provide the following lemma:
        \begin{lemma}\label{lemm:diffpartxpHesscont}
            Let $\funtwo\colon\mani\to\setR[]$ be of class $C^{3}$.
            Given $\ptaccum\in\mani$, there exist $\constthr[\funtwo] > 0$ and a closed neighborhood $\closedsubset_{\funtwo}\subseteq\mani$ of $\ptaccum$ such that, for any $\pt[]\in\closedsubset_{\funtwo}$, 
            $\opnorm{\Hess\funtwo\paren*{\ptaccum} - \partxp[]{\ptaccum}{\pt[]} \circ \Hess\funtwo\paren*{\pt[]} \circ \partxp[]{\pt[]}{\ptaccum}} \leq \constthr[\funtwo]\Riemdist{\pt[]}{\ptaccum}.$
        \end{lemma}
        \begin{proof}
            Let $\closedsubset_{\funtwo}\coloneqq\brc*{\pt\in\mani \relmiddle{\colon} \Riemdist{\pt}{\ptaccum} \leq \frac{1}{2}\injradius\paren*{\ptaccum}}$ be the ball centered at $\ptaccum$,
            and define the product distance on $\mani\times\tanspc[\ptaccum]\mani$ as 
            \begin{align}
                \Riemdist[{\mani\times\tanspc[\ptaccum]\mani}]{\paren*{\ptone[1],\tanvecone[\ptaccum]}}{\paren*{\ptone[2],\tanvectwo[\ptaccum]}} \coloneqq \sqrt{\Riemdist{\ptone[1]}{\ptone[2]}^{2} + \Riemnorm[\ptaccum]{\tanvecone[\ptaccum] - \tanvectwo[\ptaccum]}^{2}}
            \end{align}
            for $\ptone[1], \ptone[2] \in \mani$ and $\tanvecone[\ptaccum], \tanvectwo[\ptaccum]\in\tanspc[\ptaccum]\mani$.
            Let $\hesscomposfun\colon\closedsubset_{\funtwo}\times\tanspc[\ptaccum]\mani\to\tanspc[\ptaccum]\mani\colon\paren*{\pt,\tanvecone[\ptaccum]}\mapsto\partxp[]{\ptaccum}{\pt[]} \circ \Hess\funtwo\paren*{\pt[]} \circ \partxp[]{\pt[]}{\ptaccum}\sbra*{\tanvecone[\ptaccum]}$.
            Since $\funtwo$ is of class $C^{3}$, $\hesscomposfun$ is of class $C^{1}$.
            Hence, there exists $\constthr[\funtwo] > 0$ such that, for any $\pt\in\closedsubset_{\funtwo}$ and  $\tanvecone[\ptaccum]\in\tanspc[\ptaccum]\mani$ with $\Riemnorm[\ptaccum]{\tanvecone[\ptaccum]} \leq 1$, 
            \begin{align}
                \begin{split}\label{ineq:diffhessnormLipschitz}
                    &\Riemnorm[\ptaccum]{\paren*{\Hess\funtwo\paren*{\ptaccum} - \partxp[]{\ptaccum}{\pt[]} \circ \Hess\funtwo\paren*{\pt[]} \circ \partxp[]{\pt[]}{\ptaccum}}\sbra*{\tanvecone[\ptaccum]}}\\
                    &= \Riemnorm[\ptaccum]{\hesscomposfun\paren*{\pt,\tanvecone[\ptaccum]} - \hesscomposfun\paren*{\ptaccum,\tanvecone[\ptaccum]}} \leq \constthr[\funtwo] \Riemdist[{\mani\times\tanspc[\ptaccum]\mani}]{\paren*{\pt,\tanvecone[\ptaccum]}}{\paren*{\ptaccum,\tanvecone[\ptaccum]}}\\
                    &= \constthr[\funtwo]\Riemdist{\pt}{\ptaccum}.
                \end{split}
            \end{align}
            For every $\pt\in\closedsubset_{\funtwo}$, there exists a vector $\tanvectwo[\ptaccum]$ with $\Riemnorm[\ptaccum]{\tanvectwo[\ptaccum]} \leq 1$ such that 
            \begin{align}\label{eq:diffhessopnormriemnorm}
                &\Riemnorm[\ptaccum]{\paren*{\Hess\funtwo\paren*{\ptaccum} - \partxp[]{\ptaccum}{\pt[]} \circ \Hess\funtwo\paren*{\pt[]} \circ \partxp[]{\pt[]}{\ptaccum}}\sbra*{\tanvectwo[\ptaccum]}}\\
                &= \opnorm{\Hess\funtwo\paren*{\ptaccum} - \partxp[]{\ptaccum}{\pt[]} \circ \Hess\funtwo\paren*{\pt[]} \circ \partxp[]{\pt[]}{\ptaccum}}
            \end{align}
            since the set $\brc*{\tanvecone[\ptaccum]\in\tanspc[\ptaccum]\mani\colon\Riemnorm[\ptaccum]{\tanvecone[\ptaccum]} \leq 1}$ is compact.
            Combining \cref{eq:diffhessopnormriemnorm} with \cref{ineq:diffhessnormLipschitz} yields 
            \begin{align}
                \opnorm{\Hess\funtwo\paren*{\ptaccum} - \partxp[]{\ptaccum}{\pt[]} \circ \Hess\funtwo\paren*{\pt[]} \circ \partxp[]{\pt[]}{\ptaccum}} \leq \constthr[\funtwo]\Riemdist{\pt[]}{\ptaccum}
            \end{align}
            for every $\pt\in\closedsubset_{\funtwo}$.
            The proof is complete.
        \end{proof}
        We end the subsection by formally stating the sufficient condition for \cref{assu:LipschitzHess} as follows:
        \begin{corollary}
             \cref{assu:LipschitzHess} holds if $\objfun$ and $\brc*{\ineqfun[\ineqidx]}_{\ineqidx\in\ineqset}$ are of class $C^{3}$.
        \end{corollary}

    \section{Complete proofs in \texorpdfstring{\cref{sec:globconv}}{global convergence}}
    In this appendix, we provide the complete proofs for the lemmas and theorem in the global convergence analysis.
        
    \subsection{Proofs of lemmas  \texorpdfstring{in \cref{subsec:wellposedness}}{for consistency of inner iteration}}\label{appx:proofconsistency}
        In this subsection, we provide the proofs of \cref{lemm:ineqfunfeasi,lemm:firstdirderivmeritfun,lemm:twicedirderivmeritfun}.
        
    \begin{proof}[Proof of \cref{lemm:ineqfunfeasi}]
        Define $\constone\coloneqq\min_{\ineqidx\in\ineqset} \frac{1}{2} \ineqfun[\ineqidx]\paren*{\pt} > 0$.
        It follows from the continuity of $\Retr$ and $\brc*{\ineqfun[\ineqidx]}_{\ineqidx\in\ineqset}$ that $\brc*{\pullback[\pt]{\ineqfun[\ineqidx]}}_{\ineqidx\in\ineqset}$ are all continuous.
        Therefore, for each $\ineqidx\in\ineqset$ and any $\pt\in\strictfeasirgn$, there exists $\deltaone{\ineqidx} > 0$ such that $\abs{\pullback[\pt]{\ineqfun[\ineqidx]}\paren*{\tanvecone[\pt]} - \pullback[\pt]{\ineqfun[\ineqidx]}\paren*{\zerovec[\pt]}} \leq \constone$ for $\tanvecone[\pt]\in\tanspc[\pt]\mani$ with $\Riemnorm[\pt]{\tanvecone[\pt]}\leq\deltaone{\ineqidx}$,
        which implies
        \begin{align}
            \pullback[\pt]{\ineqfun[\ineqidx]}\paren*{\tanvecone[\pt]} \geq \pullback[\pt]{\ineqfun[\ineqidx]}\paren*{\zerovec[\pt]} - \constone
            \geq \frac{1}{2} \ineqfun[\ineqidx]\paren*{\pt} > 0,
        \end{align}
        where the second inequality follows from \cref{eq:retrzero} and the definition of $\constone$.
        Define $\deltaone{} \coloneqq \min_{\ineqidx\in\ineqset} \deltaone{\ineqidx}$. 
        Then, the statement holds for any $\Riemnorm[\pt]{\tanvecone[\pt]}\leq\deltaone{}$, which completes the proof. 
    \end{proof}

    \begin{proof}[Proof of \cref{lemm:firstdirderivmeritfun}]
        \isextendedversion{\Eqcref{eq:dirderivmeritfun}}{The equation} 
        directly follows from
        \begin{align}
            \begin{split}\label{eq:meritfundirderiv}
                &\D\pullback[\pt]{\meritfun[{\barrparam[]}]}\paren*{\zerovec[\pt]}\sbra*{\tanvecone[\pt]} =\D\meritfun[{\barrparam[]}]\paren*{\retr[\pt]\paren*{\zerovec[\pt]}}\sbra*{\D\retr[\pt]\paren*{\zerovec[\pt]}\sbra*{\tanvecone[\pt]}} =\D\meritfun[{\barrparam[]}]\paren*{\pt}\sbra*{\tanvecone[\pt]}\\
                &=\D\objfun\paren*{\pt}\sbra*{\tanvecone[\pt]} - \barrparam[]\sum_{\ineqidx\in\ineqset}\frac{1}{{\ineqfun[\ineqidx]}\paren*{\pt}}\D\ineqfun[\ineqidx]\paren*{\pt}\sbra*{\tanvecone[\pt]}
                =\metr[\pt]{\gradstr\objfun\paren*{\pt} - \barrparam[]\ineqgradopr[\pt]\sbra*{\inv{\Ineqfunmat[]\paren*{\pt}}\onevec}}{\tanvecone[\pt]}\\
                &=\metr[\pt]{\trslin[{\barrparam[]}]\paren*{\pt}}{\tanvecone[\pt]},
            \end{split}
        \end{align}
        where the first equality follows from the chain rule, the second one from \isextendedversion{\cref{eq:retrdef}}{the definition of the retraction}, the fourth one from \isextendedversion{\cref{eq:riemgrad}}{the definition of the Riemannian gradient}, and the last one from 
        \isextendedversion{\cref{eq:trslindef}}{the definition of $\trslin[{\barrparam[]}]$}.
    \end{proof}

    \begin{proof}[Proof of \cref{lemm:twicedirderivmeritfun}]
        For each $\ineqidx\in\ineqset$, all $\tanvectwo[\pt]\in\tanspc[\pt]\mani$, and any $\tanvecone[\pt]\in\tanspc[\pt]\mani$ with $\pullback[\pt]{{\ineqfun[\ineqidx]}}\paren*{\tanvecone[\pt]} \neq 0$, it follows that
        \begin{align}
            &\metr[\pt]{\gradstr\log\pullback[\pt]{\ineqfun[\ineqidx]}\paren*{\tanvecone[\pt]}}{\tanvectwo[\pt]} = \D\paren*{\log\pullback[\pt]{\ineqfun[\ineqidx]}}\paren*{\tanvecone[\pt]}\sbra*{\tanvectwo[\pt]}\\
            &=\D\log\paren*{\pullback[\pt]{\ineqfun[\ineqidx]}\paren*{\tanvecone[\pt]}}\sbra*{\D\pullback[\pt]{\ineqfun[\ineqidx]}\paren*{\tanvecone[\pt]}\sbra*{\tanvectwo[\pt]}} =\metr[\pt]{\frac{1}{\pullback[\pt]{\ineqfun[\ineqidx]}\paren*{\tanvecone[\pt]}}\gradstr\pullback[\pt]{\ineqfun[\ineqidx]}\paren*{\tanvecone[\pt]}}{\tanvectwo[\pt]},
        \end{align}
        where the second equality holds by the chain rule on $\tanspc[\pt]\mani$ and the third one by \isextendedversion{\cref{eq:riemgrad}}{the definition of the Riemannian gradient}.
        Therefore, we have 
        \begin{align}
            \gradstr\log\pullback[\pt]{\ineqfun[\ineqidx]}\paren*{\tanvecone[\pt]} = \frac{1}{\pullback[\pt]{\ineqfun[\ineqidx]}\paren*{\tanvecone[\pt]}}\gradstr\pullback[\pt]{\ineqfun[\ineqidx]}\paren*{\tanvecone[\pt]},
        \end{align}
        which implies
        \begin{align}
            \begin{split}\label{eq:Hesslogineqfun}
                &\Hess\paren*{\log\pullback[\pt]{\ineqfun[\ineqidx]}}\paren*{\tanvecone[\pt]} =\D\paren*{\gradstr\log\pullback[\pt]{\ineqfun[\ineqidx]}}\paren*{\tanvecone[\pt]} = \D\paren*{\frac{1}{\pullback[\pt]{\ineqfun[\ineqidx]}\paren*{\cdot}}\gradstr\pullback[\pt]{\ineqfun[\ineqidx]}\paren*{\cdot}}\paren*{\tanvecone[\pt]}\\
                &= -\frac{\gradstr\pullback[\pt]{\ineqfun[\ineqidx]}\paren*{\tanvecone[\pt]}}{{\pullback[\pt]{\ineqfun[\ineqidx]}\paren*{\tanvecone[\pt]}}^{2}}\D\pullback[\pt]{\ineqfun[\ineqidx]}\paren*{\tanvecone[\pt]} + \frac{1}{\pullback[\pt]{\ineqfun[\ineqidx]}\paren*{\tanvecone[\pt]}}\D\paren*{\gradstr\pullback[\pt]{\ineqfun[\ineqidx]}}\paren*{\tanvecone[\pt]}\\
                &=-\frac{\gradstr\pullback[\pt]{\ineqfun[\ineqidx]}\paren*{\tanvecone[\pt]}}{{\pullback[\pt]{\ineqfun[\ineqidx]}\paren*{\tanvecone[\pt]}}^{2}}\D\pullback[\pt]{\ineqfun[\ineqidx]}\paren*{\tanvecone[\pt]} + \frac{1}{\pullback[\pt]{\ineqfun[\ineqidx]}\paren*{\tanvecone[\pt]}}\Hess\pullback[\pt]{\ineqfun[\ineqidx]}\paren*{\tanvecone[\pt]}.
            \end{split}
        \end{align}
        Using \cref{eq:Hesslogineqfun}, we derive 
        \begin{align}
            &\D[2]\pullback[\pt]{\meritfun[{\barrparam[]}]}\paren*{\tanvecone[\pt]}\sbra*{\tanvectwo[\pt],\tanvecthr[\pt]} =\D[2]\pullback[\pt]{\objfun}\paren*{\tanvecone[\pt]}\sbra*{\tanvectwo[\pt],\tanvecthr[\pt]} - \barrparam[]\sum_{\ineqidx\in\ineqset}\D[2]\paren*{\log\pullback[\pt]{\ineqfun[\ineqidx]}}\paren*{\tanvecone[\pt]}\sbra*{\tanvectwo[\pt],\tanvecthr[\pt]}\\
            &=\metr[\pt]{\paren*{\Hess\pullback[\pt]{\objfun}\paren*{\tanvecone[\pt]} - \barrparam[] \sum_{\ineqidx\in\ineqset} \Hess\paren*{\log\pullback[\pt]{\ineqfun[\ineqidx]}}\paren*{\tanvecone[\pt]}}\sbra*{\tanvectwo[\pt]}}{\tanvecthr[\pt]}\\
            &=\metr[\pt]{\paren*{\Hess\pullback[\pt]{\objfun}\paren*{\tanvecone[\pt]} -  \sum_{\ineqidx\in\ineqset}\frac{\barrparam[]}{\pullback[\pt]{\ineqfun[\ineqidx]}\paren*{\tanvecone[\pt]}}\Hess\pullback[\pt]{\ineqfun[\ineqidx]}\paren*{\tanvecone[\pt]}}\sbra*{\tanvectwo[\pt]}}{\tanvecthr[\pt]}\\
            &\quad + \sum_{\ineqidx\in\ineqset}\frac{\barrparam[]}{{\pullback[\pt]{\ineqfun[\ineqidx]}\paren*{\tanvecone[\pt]}}^{2}}\metr[\pt]{\gradstr\pullback[\pt]{\ineqfun[\ineqidx]}\paren*{\tanvecone[\pt]}}{\tanvectwo[\pt]}\metr[\pt]{\gradstr\pullback[\pt]{\ineqfun[\ineqidx]}\paren*{\tanvecone[\pt]}}{\tanvecthr[\pt]},
        \end{align}
        where the third equality follows from \isextendedversion{\cref{eq:riemgrad}}{the definition of the Riemannian gradient}.
        The proof is complete.
    \end{proof}

    \subsection{Proofs of lemmas and theorem in \texorpdfstring{\cref{subsec:globconvinner}}{global convergence of inner iteration}}\label{appx:proofinnerglobal}
        In this subsection, we provide the proofs of \cref{lemm:meritfunanal,lemm:ineqfunbound,theo:limsupsosp,lemm:trsquadbounded,theo:innerglobconv}.

    \begin{proof}[Proof of \cref{lemm:meritfunanal}]
        We first consider \cref{lemm:meritfunmonnonincr}. 
        When $\card[\succset]<\infty$, the entire sequence $\brc*{\ptinriter}_{\inriteridx}$ reduces to the finite sequence of the successful iterates $\brc*{\ptinrsuccsubseqiter}_{\succsubseqidx}$, and we can easily obtain the result in this case.
        In the following, we consider the case where $\card[\succset]=\infty$.
        Note that, for any $\succsubseqidx\in\setNz$, ${\pt[]}^{1+\inrsuccsubseqiter}={\pt[]}^{2+\inrsuccsubseqiter}=\cdots=\ptinrsuccsubseqiterp$ holds: once the iterate is updated and set to ${\pt[]}^{1+\inrsuccsubseqiter}$ at the $\inrsuccsubseqiter$-th iteration, it remains unchanged until it is updated again and set to ${\pt[]}^{1+\inrsuccsubseqiterp}$ at the $\inrsuccsubseqiterp$-th iteration.
        Thus, we have
        \begin{align}
            &\meritfun[{\barrparam[]}]\paren*{\ptinrsuccsubseqiter} - \meritfun[{\barrparam[]}]\paren*{\ptinrsuccsubseqiterp} =\meritfun[{\barrparam[]}]\paren*{\ptinrsuccsubseqiter} - \meritfun[{\barrparam[]}]\paren*{{\pt[]}^{1+\inrsuccsubseqiter}}\\ 
            &\geq \trratiothold\paren*{\modelfuninrsuccsubseqiter\paren*{\zerovec[\ptinrsuccsubseqiter]} - \modelfuninrsuccsubseqiter\paren*{\dirinrsuccsubseqiter}} \geq \trratiothold\constCauchy\Riemnorm[\ptinrsuccsubseqiter]{\trslininrsuccsubseqiter}\min\paren*{\trradiusinrsuccsubseqiter, \frac{\Riemnorm[\ptinrsuccsubseqiter]{\trslininrsuccsubseqiter}}{\opnorm{\trsquadinrsuccsubseqiter}}} \geq 0,
        \end{align}
        where the first inequality holds since the $\inrsuccsubseqiter$-th iteration is successful and the second one follows from \cref{assu:Cauchydecreasing}.
        Therefore, $\brc*{\meritfun[{\barrparam[]}]\paren*{\ptinrsuccsubseqiter}}_{\succsubseqidx}$ is monotonically non-increasing.
        Since, again, all iterates between $\inrsuccsubseqiter$ and $\inrsuccsubseqiterp$ are unsuccessful, the lower bound extends to all the iterates, and the proof is complete.

        Next, we consider \cref{lemm:meritfunconv}. 
        By \assuenumicref{assu:ineqfunbounded}{assu:ineqfunboundedvalue} and \assuenumicref{assu:objfunbounded}{assu:objfunboundedbelow}, the sequence $\brc*{\meritfun[{\barrparam[]}]\paren*{\ptinrsuccsubseqiter}}_{\succsubseqidx}$ is bounded below.
        Thus, \cref{lemm:meritfunconv} follows by the monotone convergence theorem.
    \end{proof}

    \begin{proof}[Proof of \cref{lemm:ineqfunbound}]
        We first consider \cref{lemm:ienqfunptlowbd}.
        When $\card[\succset] < \infty$, 
        the entire sequence $\brc*{\ptinriter}_{\inriteridx}$ reduces to the finite one of the successful iterates $\brc*{\ptinrsuccsubseqiter}_{\succsubseqidx}$.
        In this case, we can easily obtain the conclusion by letting $\epsfiv\coloneqq\min_{\ineqidx, \succsubseqidx} \frac{1}{2} \ineqfun[\ineqidx]\paren*{\ptinrsuccsubseqiter}$ and $\deltathr \coloneqq \min_{\succsubseqidx}\thldtrradius[\ptinrsuccsubseqiter]$ in \cref{lemm:ineqfunfeasi}.
        In the following, we consider the case 
        where  $\card[\succset] = \infty$. 
        Recall that, from \enumicref{lemm:meritfunanal}{lemm:meritfunmonnonincr}, the sequence $\brc*{\meritfun[{\barrparam[]}]\paren*{\ptinriter}}_{\inriteridx}$ is bounded above.
        Combining this with 
        \cref{eq:logbarrfun} 
        and \assuenumicref{assu:objfunbounded}{assu:objfunboundedbelow}, we obtain that $\sum_{\ineqidx\in\ineqset} \log\ineqfun[\ineqidx]\paren*{\ptinriter} = \inv{\barrparam[]} \paren*{\objfun\paren*{\ptinriter} - \meritfun[{\barrparam[]}]\paren*{\ptinriter}}$ is bounded below.
        This, together with \assuenumicref{assu:ineqfunbounded}{assu:ineqfunboundedvalue}, implies that there exists $\epsfiv > 0$ such that $\ineqfun\paren*{\ptinriter} \geq \epsfiv$ for all $\inriteridx\in\setNz$.
        Thus, \cref{lemm:ienqfunptlowbd} is established.
        
        Next, we consider \cref{lemm:ienqfununiflowbd}.
        Let $\epsthr \coloneqq \frac{3}{4}\epsfiv > 0$.
        By \cref{assu:pullbackineqfuncont}, for each $\ineqidx\in\ineqset$, there exists $\thldtrradiusthr[\ineqidx] > 0$ such that, for any $\inriteridx\in\setNz$ and all $\tanvecone[\ptinriter]\in\tanspc[\ptinriter]\mani$ with $\Riemnorm[\ptinriter]{\tanvecone[\ptinriter]}\leq\thldtrradiusthr[\ineqidx]$, it follows that $\abs*{\pullback[\ptinriter]{{\ineqfun[\ineqidx]}}\paren*{\tanvecone[\ptinriter]} - \ineqfun[\ineqidx]\paren*{\ptinriter}} \leq \frac{1}{3}\epsthr$, which implies $\pullback[\ptinriter]{{\ineqfun[\ineqidx]}}\paren*{\tanvecone[\ptinriter]} \geq \ineqfun[\ineqidx]\paren*{\ptinriter} - \frac{1}{3}\epsthr \geq \frac{4}{3}\epsthr - \frac{1}{3}\epsthr = \epsthr$.
        By letting $\deltathr\coloneqq \min_{\ineqidx\in\ineqset}\thldtrradiusthr[\ineqidx] > 0$, we complete the proof of \cref{lemm:ienqfununiflowbd}.
    \end{proof}

    \begin{proof}[Proof of \cref{lemm:trsquadbounded}]
        For any $\inriteridx\in\setNz$ and all $\tanvecone[\ptinriter]\in\tanspc[\ptinriter]\mani$ with $\Riemnorm[\ptinriter]{\tanvecone[\ptinriter]} \leq 1$, it holds that
        \begin{align}
            &\Riemnorm[\ptinriter]{\trsquadinriter\sbra*{\tanvecone[\ptinriter]}} \leq \opnorm{\Hess\objfun\paren*{\ptinriter}}\Riemnorm[\ptinriter]{\tanvecone[\ptinriter]}\\
            &\quad + \sum_{\ineqidx\in\ineqset} \ineqLagmultinriter[\ineqidx]\opnorm{\Hess\ineqfun[\ineqidx]\paren*{\ptinriter}}\Riemnorm[\ptinriter]{\tanvecone[\ptinriter]} + \sum_{\ineqidx\in\ineqset} \frac{\ineqLagmultinriter[\ineqidx]}{\ineqfun[\ineqidx]\paren*{\ptinriter}}\Riemnorm[\ptinriter]{\gradstr\ineqfun[\ineqidx]\paren*{\ptinriter}}^{2}\Riemnorm[\ptinriter]{\tanvecone[\ptinriter]}^{2}\\
            &\leq \opnorm{\Hess\objfun\paren*{\ptinriter}} + \sum_{\ineqidx\in\ineqset} \ineqLagmultinriter[\ineqidx]\opnorm{\Hess\ineqfun[\ineqidx]\paren*{\ptinriter}} + \sum_{\ineqidx\in\ineqset} \frac{\ineqLagmultinriter[\ineqidx]}{\ineqfun[\ineqidx]\paren*{\ptinriter}}\Riemnorm[\ptinriter]{\gradstr\ineqfun[\ineqidx]\paren*{\ptinriter}}^{2},
        \end{align}
        where the first inequality follows from 
        \isextendedversion{\cref{eq:trsquaddef}}{the definition of $\trsquadinriter$}
        and the second one from $\Riemnorm[\ptinriter]{\tanvecone[\ptinriter]} \leq 1$.
        Under \cref{assu:ineqfunbounded}, \cref{assu:ineqLagmultbounded}, \assuenumicref{assu:objfunbounded}{assu:objfunHessbounded},  and \enumicref{lemm:ineqfunbound}{lemm:ienqfunptlowbd}, the right-hand side is bounded above.
        The proof is complete.
    \end{proof}

    \begin{proof}[Proof of \cref{theo:limsupsosp}]
        We argue by contradiction.
        Suppose that there exist $\constEigen > 0$ and an index $\tholdidxthr\in\setNz$ such that $\mineigval\sbra*{\trsquadinriter} < - \constEigen$ for all $\inriteridx\geq\tholdidxthr$.
        Let $\coefffou\coloneqq\frac{1}{2}\constEigen\epseigen > 0$ and let $\thldtrradiussix > 0$ be the associated threshold value in \enumicref{lemm:diffpredared}{lemm:diffpredaredsmallo}.
        Then, since $\predinriter > 0$ holds by 
        \cref{ineq:Eigendecreasing}, 
        we have
        \begin{align}
            \abs*{\trratioinriter - 1} 
            = \frac{\abs*{\predinriter - \aredinriter}}{\abs*{\predinriter}} \leq \frac{\coefffou\Riemnorm[\ptinriter]{\dirptinriter}^{2}}{\constEigen\epseigen\paren*{\trradiusinriter}^{2}} \leq \frac{1}{2}
        \end{align}
        whenever $\trradiusinriter \leq \thldtrradiussix$,  where the first inequality follows from \enumicref{lemm:diffpredared}{lemm:diffpredaredsmallo} and 
        \cref{ineq:Eigendecreasing}
        again and the second one holds by the definition of $\coefffou$ and $\Riemnorm[\ptinriter]{\dirptinriter} \leq \trradiusinriter$.
        Thus, if $\trradiusinriter \leq \thldtrradiussix$, we have $\trratioinriter\geq\frac{1}{2}>\frac{1}{4}>\trratiothold$, which implies that the update $\trradiusinriterp=\frac{1}{4}\trradiusinriter$ can occur only if $\trradiusinriter > \thldtrradiussix$.
        Hence,
        \begin{equation}\label{ineq:trradiuslowboundEigen}
            \trradiusinriter\geq\min\paren*{\trradiusinritertholdthr, \frac{\thldtrradiussix}{4}}
        \end{equation}
        holds for all $\inriteridx\geq\tholdidxthr$.
        
        On the other hand, recall that $\succset$ is the set of the successful iterations and $\brc*{\ptinrsuccsubseqiter}_{\succsubseqidx}$ is the ordered sequence of $\succset$.
        When $\abs*{\succset}$ is finite, for all $\inriteridx\in\setNz$, the iteration is unsuccessful.
        Thus, it follows that $\lim_{\inriteridx\to\infty}\trradiusinriter = 0$, which contradicts \cref{ineq:trradiuslowboundEigen}.
        In the following, we consider the case where $\abs*{\succset}$ is infinite.
        For all $\succsubseqidx\in\setNz$, we obtain
        \begin{align}
            \begin{split}\label{ineq:diffmeritfunEigen}
                &\meritfun[{\barrparam[]}]\paren*{\ptinrsuccsubseqiter} - \meritfun[{\barrparam[]}]\paren*{\ptinrsuccsubseqiterp} = \meritfun[{\barrparam[]}]\paren*{\ptinrsuccsubseqiter} - \meritfun[{\barrparam[]}]\paren*{\pt^{1 + \inrsuccsubseqiter}}\\
                &\geq \trratiothold\paren*{\modelfuninrsuccsubseqiter\paren*{\zerovec[\ptinrsuccsubseqiter]} - \modelfuninrsuccsubseqiter\paren*{\dirinrsuccsubseqiter}} \geq \trratiothold\constEigen\paren*{\trradiusinrsuccsubseqiter}^{2}\epseigen\geq 0,
            \end{split}
        \end{align}
        where the equality follows since $\ptinrsuccsubseqiterp=\pt^{1 + \inrsuccsubseqiter}$ holds, the first inequality does since the $\inrsuccsubseqiter$-th iterate is successful, and the second one follows from \cref{assu:Eigendecreasing}.
        From \enumicref{lemm:meritfunanal}{lemm:meritfunconv}, $\brc*{\meritfun[{\barrparam[]}]\paren*{\ptinriter}}_{\inriteridx}$ is a Cauchy sequence.
        Therefore, it follows from \cref{ineq:diffmeritfunEigen} that 
        $\lim_{\succsubseqidx\to\infty}\trradiusinrsuccsubseqiter = 0$, which contradicts \cref{ineq:trradiuslowboundEigen}.

        Now, we have that, for any  $\constEigen[\prime] > 0$ and any $\tholdidxthr\in\setNz$, there exists $\inritertholdthridx\geq\tholdidxthr$ such that $\mineigval\sbra*{\trsquadinritertholdthr} \geq - \constEigen[\prime]$ holds.
        This implies the limit superior of the sequence $\mineigval\sbra*{\trsquadinriter}$ is not less than zero.
        Indeed, letting $\constEigen[\prime] =\inv{\tholdidxthr}$ for all $\tholdidxthr\in\setN$, we obtain $\sup_{\inriteridx\geq\tholdidxthr} \mineigval\sbra*{\trsquadinriter} \geq \mineigval\sbra*{\trsquadinritertholdthr} \geq -\inv{\tholdidxthr}$ for all $\tholdidxthr\in\setN$.
        Taking the limit as $\tholdidxthr\to\infty$ leads to $\limsup_{\inriteridx\to\infty}\mineigval\sbra*{\trsquadinriter} \geq 0$.
    \end{proof}

    \subsection{Proof of theorems in \texorpdfstring{\cref{subsec:globconvouter}}{global convergence of outer iterations}}\label{appx:proofouterglobal}
        In this section, we provide the proofs of \cref{theo:limglobconvotr} and a lemma that complements the proof of \cref{theo:globconvotriterptaccum}.

        We first provide the proof of \cref{theo:limglobconvotr}.
        
    \begin{proof}[Proof of \cref{theo:limglobconvotr}]
        We can easily obtain \isextendedversion{\cref{eq:limgradLagfunotriter,eq:liminffeasiotriter}}{the first equality in \cref{eq:limgradLagfuncomplotriter} and the conditions~\cref{eq:liminfprimaldualfeasiotriter}} from \cref{eq:stopcondKKT} and \cref{eq:stopcondstrictfeasi}, respectively.
        As for \isextendedversion{\cref{eq:limcomplotriter}}{the second equality in \cref{eq:limgradLagfuncomplotriter}}, it follows from \cref{eq:stopcondbarrcompl} that
        \begin{align}
            \norm*{\IneqLagmultmat[\otriteridx]\ineqfun[]\paren*{\ptotriter}} \leq \norm*{\IneqLagmultmatotriter \ineqfun[]\paren*{\ptotriter} - \barrparamotriterm \onevec} + \barrparamotriterm\norm*{\onevec}  \leq \forcingfuncompl\paren*{\barrparamotriterm} + \ineqdime\barrparamotriterm,
        \end{align}
        where the right-hand side converges to zero as $\otriteridx\to\infty$.
        Thus, \isextendedversion{\eqcref{eq:limcomplotriter}}{the second equality in \cref{eq:limgradLagfuncomplotriter}} holds.

        Next, additionally supposing that \cref{assu:gradineqfunotriterbounded} holds and \seccondflag{} is \True{} in \cref{algo:RIPTRMOuter}, we show \cref{eq:liminfHessLagfunotriter} by contradiction.
        Assume that there exists 
        a sequence of tangent vectors $\brc*{\tanvecone[\ptotrsubseqiter]}_{\subseqidx}$ satisfying \cref{eq:asymptoticwctanvec} 
        such that 
        \begin{equation}\label{ieq:contraHessLagneg}
            \liminf_{\subseqidx\to\infty} \metr[\ptotrsubseqiter]{\Hess\Lagfun\paren*{\allvarotrsubseqiter}\sbra*{\tanvecone[\ptotrsubseqiter]}}{\tanvecone[\ptotrsubseqiter]} < 0.
        \end{equation}
        Define $\consttwo \coloneqq \min_{\ineqidx\in\ineqset\backslash\asymptotactsubset}\liminf_{\subseqidx\to\infty} \ineqfun[\ineqidx]\paren*{\ptotrsubseqiter} > 0$.
        Note that $\ineqfun[\ineqidx]\paren*{\ptotrsubseqiter} \geq \frac{1}{2}\consttwo > 0$ holds for every $\ineqidx\in\ineqset\backslash\asymptotactsubset$ and all $\subseqidx\in\setNz$ sufficiently large.
        Therefore, together with \isextendedversion{\cref{eq:limcomplotriter}}{the second equality in \cref{eq:limgradLagfuncomplotriter}}, it follows that $\lim_{\subseqidx\to\infty}\sbra*{\ineqLagmultotrsubseqiter[]}_{\ineqidx} = 0$ for every $\ineqidx\in\ineqset\backslash\asymptotactsubset$.
        For all $\subseqidx\in\setNz$ sufficiently large, we write $\IneqLagmultmat[\otrsubseqiter]$ for $\diag\paren*{\ineqLagmultotrsubseqiter[]}$ and have
        \begin{align}
            &\metr[\ptotrsubseqiter]{\Hess\Lagfun\paren*{\allvarotrsubseqiter}\sbra*{\tanvecone[\ptotrsubseqiter]}}{\tanvecone[\ptotrsubseqiter]}\\
            &= \metr[\ptotrsubseqiter]{\trsquad\paren*{\allvarotrsubseqiter}\sbra*{\tanvecone[\ptotrsubseqiter]}}{\tanvecone[\ptotrsubseqiter]} - \metr[\ptotrsubseqiter]{\ineqgradopr[\ptotrsubseqiter]\IneqLagmultmat[\otrsubseqiter]\inv{\Ineqfunmat[]\paren*{\ptotrsubseqiter}}\coineqgradopr[\ptotrsubseqiter]\sbra*{\tanvecone[\ptotrsubseqiter]}}{\tanvecone[\ptotrsubseqiter]}\\
            &\geq - \forcingfunsosp\paren*{\barrparamotrsubseqiterm} - \sum_{\ineqidx\in\ineqset}\frac{\sbra*{\ineqLagmultotrsubseqiter[]}_{\ineqidx}}{\ineqfun[\ineqidx]\paren*{\ptotrsubseqiter}} \metr[\ptotrsubseqiter]{\gradstr\ineqfun[\ineqidx]\paren*{\ptotrsubseqiter}}{\tanvecone[\ptotrsubseqiter]}^{2}\\
            &\geq - \forcingfunsosp\paren*{\barrparamotrsubseqiterm} - \sum_{\ineqidx\in\ineqset\backslash\asymptotactsubset}\frac{2\sbra*{\ineqLagmultotrsubseqiter[]}_{\ineqidx}}{\consttwo} \Riemnorm[\ptotrsubseqiter]{\gradstr\ineqfun[\ineqidx]\paren*{\ptotrsubseqiter}}^{2}\Riemnorm[\ptotrsubseqiter]{\tanvecone[\ptotrsubseqiter]}^{2},
        \end{align}
        where the equality follows from 
        \isextendedversion{\cref{eq:trsquaddef}}{the definition of $\trsquad$}, 
        the first inequality follows from \cref{eq:stopcondsecondord}, and the second one from $\ineqfun[\ineqidx]\paren*{\ptotrsubseqiter} \geq \frac{1}{2}\consttwo$ and 
        \isextendedversion{\cref{eq:asymptoticwctanvecorthogonal}}{\cref{eq:asymptoticwctanvec}}.
        Since $\brc*{\Riemnorm[\ptotrsubseqiter]{\gradstr\ineqfun[\ineqidx]\paren*{\ptotrsubseqiter}}}_{\subseqidx}$ and $\brc*{\Riemnorm[\ptotrsubseqiter]{\tanvecone[\ptotrsubseqiter]}}_{\subseqidx}$ are bounded by 
        \isextendedversion{\cref{assu:gradineqfunotriterbounded,eq:asymptoticwctanvecbounded}}{\cref{assu:gradineqfunotriterbounded,eq:asymptoticwctanvec}},
        $\lim_{\subseqidx\to\infty}\forcingfunsosp\paren*{\barrparamotrsubseqiterm} = 0$ holds, and we have $\lim_{\subseqidx\to\infty}\sbra*{\ineqLagmultotrsubseqiter[]}_{\ineqidx} = 0$ for any $\ineqidx\in\ineqset\backslash\asymptotactsubset$, the right-hand side converges to zero, which contradicts \cref{ieq:contraHessLagneg}.
        The proof is complete.
    \end{proof}

        Next, we derive the following lemma used in the proof of \cref{theo:globconvotriterptaccum}:
        
    \begin{lemma}\label{lemm:LICQardptaccum}
        Under \cref{assu:LICQptaccum}, there exists a neighborhood $\subsetmanitwo\subseteq\mani$ of $\ptaccum\in\mani$ such that $\brc*{\gradstr\ineqfun[\ineqidx]\paren*{\pt}}_{\ineqidx\in\activeineqset\paren*{\ptaccum}}$ are linearly independent for any $\pt[]\in\subsetmanitwo$. 
    \end{lemma}
    \begin{proof}
        We argue by contradiction.
        Suppose that, for any $\subsetmanitwo\subseteq\mani$ with $\ptaccum\in\subsetmanitwo$, there exists $\pt[]\in\subsetmanitwo$ such that $\brc*{\gradstr\ineqfun[\ineqidx]\paren*{\pt}}_{\ineqidx\in\activeineqset\paren*{\ptaccum}}$ are linearly dependent.
        For any $\radiusone\in\setNz$, we let $\subsetmanitwo[\frac{1}{\radiusone}]\coloneqq\brc*{\pt\in\mani\colon\Riemdist{\pt}{\ptaccum} < \frac{1}{\radiusone}}$.
        There exist $\pt[\radiusone]\in\subsetmanitwo[\frac{1}{\radiusone}]$ and a nonzero vector $\vectwo[\radiusone]\in\setR[\abs*{\activeineqset\paren*{\ptaccum}}]$ such that $\norm{\vectwo[\radiusone]}=1$ and $\sum_{\ineqidx\in\activeineqset\paren*{\ptaccum}}\vectwo[\radiusone]_{\ineqidx}\gradstr\ineqfun[\ineqidx]\paren*{\pt[\radiusone]}=0$.
        Considering $\radiusone=1, 2, 3,\ldots$ yields that $\brc*{\pt[\radiusone]}_{\radiusone}$ converges to $\ptaccum$ as $\radiusone\to\infty$, and there exists $\vectwo[\accumsymbol]\in\setR[\abs*{\activeineqset\paren*{\ptaccum}}]$ such that $\norm{\vectwo[\accumsymbol]}=1$ and $\sum_{\ineqidx\in\activeineqset\paren*{\ptaccum}}\vectwo[\accumsymbol]_{\ineqidx}\gradstr\ineqfun[\ineqidx]\paren*{\ptaccum}=0$, which contradicts \cref{assu:LICQptaccum}.
        The proof is complete.
    \end{proof}

        Now, we provide the proof of \cref{theo:globconvotriterptaccum}.

    \newcommand{\Akn}[1][]{A_{#1}}
    \newcommand{\Bkn}[1][]{B_{#1}}
    \newcommand{\Ckn}[1][]{C_{#1}}

    \begin{proof}[Proof of \cref{theo:globconvotriterptaccum}]
        \cref{theo:limglobconvotr} directly implies that the sequence $\brc*{\ptotrsubseqiter}_{\subseqidx}$ and the associated sequence $\brc*{\ineqLagmultotrsubseqiter[]}_{\subseqidx}$ satisfy the \isextendedversion{\AKKT{} conditions~\cref{def:AKKTconditions}}{\AKKT{} conditions}.
        Note also that, since $\brc*{\ptotrsubseqiter}_{\subseqidx}$ is convergent and hence bounded, \cref{assu:gradineqfunotriterbounded} is fulfilled.
        Hereafter, we consider the weak second-order stationarity when \seccondflag{} is \True{} in \cref{algo:RIPTRMOuter}. 
        To this end, we first analyze the limiting behaviors of the orthogonal projection and the parallel transport.
        Let $\brc*{\wcconebasis[\ineqidx]\paren*{\pt}}_{\ineqidx\in\activeineqset\paren*{\ptaccum}}$ be an orthonormal basis for $\spanstr\brc*{\gradstr\ineqfun[\ineqidx]\paren*{\pt}}_{\ineqidx\in\activeineqset\paren*{\ptaccum}}\subseteq\tanspc[\pt]\mani$.
        Under \cref{assu:LICQptaccum}, we can construct the basis around $\ptaccum$ using the conventional Gram-Schmidt process in the following manner:
        \isextendedversion{
        from \cref{lemm:LICQardptaccum}
        }{
        from \cite[Lemma D.1]{Obaraetal2025APrimalDualIPTRMfor2ndOrdStnryPtofRiemIneqCstrOptimProbs},}
        there exists a neighborhood $\subsetmanitwo\subseteq\mani$ of $\ptaccum$ such that $\brc*{\gradstr\ineqfun[\ineqidx]\paren*{\pt}}_{\ineqidx\in\activeineqset\paren*{\ptaccum}}$ are linearly independent for any $\pt\in\subsetmanitwo$.
        Without loss of generality, we let the indices of $\activeineqset\paren*{\ptaccum} = \brc*{1,2,3,\ldots, \card[\activeineqset\paren*{\ptaccum}]}$; see 
        \isextendedversion{
        \cref{def:activeineqset}
        }{
        \cref{subsec:optimcond}
        }
        for the definition of $\activeineqset\paren*{\ptaccum}$.
        We define the orthonormal basis function $\wcconebasis[\idxtwo]\colon\subsetmanitwo\to\tanspc[]\mani$ as
        $\wcconebasis[\idxtwo]\paren*{\pt} \coloneqq \frac{\basisone[\idxtwo]\paren*{\pt}}{\Riemnorm[\pt]{\basisone[\idxtwo]\paren*{\pt}}}, \text{ where } \basisone[\idxtwo]\paren*{\pt} \coloneqq \gradstr\ineqfun[\idxtwo]\paren*{\pt} - \sum_{\idxthr=1}^{\idxtwo-1} \metr[\pt]{\gradstr\ineqfun[\idxtwo]\paren*{\pt}}{\wcconebasis[\idxthr]\paren*{\pt}}\wcconebasis[\idxthr]\paren*{\pt}$
        for $\idxtwo=1,2,3,\dots,\card[\activeineqset\paren*{\ptaccum}]$.
        Note that, due to the continuity of the Riemannian gradients and the Riemannian metric, $\wcconebasis[\idxtwo]$ is a continuous vector field.
        Using the basis, we define the orthogonal projection operator as 
        \isextendedversion{
        \begin{equation}
            \projwcconeaccum\paren*{\pt}\sbra*{\tanvecone[\pt]}\coloneqq \tanvecone[\pt] - \sum_{\idxtwo\in\activeineqset\paren*{\ptaccum}}\metr[\pt]{\wcconebasis[\idxtwo]\paren*{\pt}}{\tanvecone[\pt]}\wcconebasis[\idxtwo]\paren*{\pt}
        \end{equation}
        }{
        $\projwcconeaccum\paren*{\pt}\sbra*{\tanvecone[\pt]}\coloneqq \tanvecone[\pt] - \sum_{\idxtwo\in\activeineqset\paren*{\ptaccum}}\metr[\pt]{\wcconebasis[\idxtwo]\paren*{\pt}}{\tanvecone[\pt]}\wcconebasis[\idxtwo]\paren*{\pt}$}
        for $\tanvecone[\pt]\in\tanspc[\pt]\mani$.
        By replacing $\subsetmanitwo$ with a sufficiently small neighborhood of $\ptaccum$ if necessary, it holds that,
        for any $\tanvecone[\ptaccum]\in\weakcriticalcone\paren*{\ptaccum}$, the map $\subsetmanitwo\ni\pt\mapsto\projwcconeaccum\paren*{\pt}\sbra*{\partxp{\pt}{\ptaccum}\sbra*{\tanvecone[\ptaccum]}}\in\tanspc[\pt]\mani$ is continuous from the continuities of $\wcconebasis[\idxtwo]$ and the parallel transport~\cite[Lemma A.1]{LiuBoumal2020SimpleAlgoforOptimonRiemManiwithCstr}.
        Therefore,
        \isextendedversion{
        \begin{align}
            \begin{split}\label{eq:limdiffpartxpprojpartxp}
                &\lim_{\pt\to\ptaccum} \Riemnorm[\pt]{\partxp{\pt}{\ptaccum}\sbra*{\tanvecone[\ptaccum]} - \projwcconeaccum\paren*{\pt}\sbra*{\partxp{\pt}{\ptaccum}\sbra*{\tanvecone[\ptaccum]}}}\\ 
                &=\Riemnorm[\ptaccum]{\sum_{\idxtwo\in\activeineqset\paren*{\ptaccum}}\metr[\ptaccum]{\wcconebasis[\idxtwo]\paren*{\ptaccum}}{\tanvecone[\ptaccum]}\wcconebasis[\idxtwo]\paren*{\ptaccum}}=0,
            \end{split}
        \end{align}
        }{
        \begin{equation}\label{eq:limdiffpartxpprojpartxp}
            \lim_{\pt\to\ptaccum} \Riemnorm[\pt]{\partxp{\pt}{\ptaccum}\sbra*{\tanvecone[\ptaccum]} - \projwcconeaccum\paren*{\pt}\sbra*{\partxp{\pt}{\ptaccum}\sbra*{\tanvecone[\ptaccum]}}} = \Riemnorm[\ptaccum]{\sum_{\idxtwo\in\activeineqset\paren*{\ptaccum}}\hspace{-1mm}\metr[\ptaccum]{\wcconebasis[\idxtwo]\paren*{\ptaccum}}{\tanvecone[\ptaccum]}\wcconebasis[\idxtwo]\paren*{\ptaccum}} \hspace{-1mm}= 0,
        \end{equation}}
        where the first equality follows from $\lim_{\pt\to\ptaccum}\partxp{\pt}{\ptaccum}\sbra*{\tanvecone[\ptaccum]} = \partxp{\ptaccum}{\ptaccum}\sbra*{\tanvecone[\ptaccum]} = \tanvecone[\ptaccum]$ and the second one from $\tanvecone[\ptaccum]$ and each $\wcconebasis[\idxtwo]\paren*{\ptaccum}$ being orthogonal.
        We also note that, by the continuity again, $\Riemnorm[\pt]{\partxp{\pt}{\ptaccum}\sbra*{\tanvecone[\ptaccum]}}$ and $\Riemnorm[\pt]{\projwcconeaccum\paren*{\pt}\sbra*{\partxp{\pt}{\ptaccum}\sbra*{\tanvecone[\ptaccum]}}}$ are bounded around $\ptaccum$.

        Notice that, from \cref{prop:eqivKKTAKKT,assu:LICQptaccum}, the point $\ptaccum$ is a \KKT{} point.
        From the proof of \cite[Theorem~2]{YamakawaSato2022SeqOptimCondforNLOonRiemManiandGlobConvALM}, the sequence $\brc*{\ineqLagmultotrsubseqiter[]}$ is bounded.
        Hence, without loss of generality, we may assume that $\brc*{\ineqLagmultotrsubseqiter[]}$ converges to the associated Lagrange multipliers $\ineqLagmultaccum[]\in\setR[\ineqdime]$ by taking a subsequence if necessary.
        Define $\allvaraccum\coloneqq\paren*{\ptaccum, \ineqLagmultaccum[]}\in\feasirgn\times\setRp[\ineqdime]$.
        We next analyze the behavior of the Hessian of the Lagrangian around $\allvaraccum$.
        It follows that, for all $\subseqidx\in\setNz$ sufficiently large, 
        \isextendedversion{
        \begin{align}
            &\opnorm{\Hess[{\pt[]}]\Lagfun\paren*{\allvaraccum} - \partxp[]{\ptaccum}{\ptotrsubseqiter} \circ \Hess[{\pt[]}]\Lagfun\paren*{\allvarotrsubseqiter} \circ \partxp[]{\ptotrsubseqiter}{\ptaccum}}\\
            &\leq \opnorm{\Hess\objfun\paren*{\ptaccum} - \partxp[]{\ptaccum}{\ptotrsubseqiter} \circ \Hess\objfun\paren*{\ptotrsubseqiter} \circ \partxp[]{\ptotrsubseqiter}{\ptaccum}}\\
            &\quad + \sum_{\ineqidx\in\ineqset} \abs*{\ineqLagmultaccum[\ineqidx] - \sbra*{\ineqLagmultotrsubseqiter[]}_{\ineqidx}} \opnorm{\Hess\ineqfun[\ineqidx]\paren*{\ptaccum}}\\
            &\quad+ \sum_{\ineqidx\in\ineqset}\abs*{\sbra*{\ineqLagmultotrsubseqiter[]}_{\ineqidx}} \opnorm{\Hess\ineqfun[\ineqidx]\paren*{\ptaccum} - \partxp[]{\ptaccum}{\ptotrsubseqiter} \circ \Hess\ineqfun[\ineqidx]\paren*{\ptotrsubseqiter} \circ \partxp[]{\ptotrsubseqiter}{\ptaccum}}\\
            &\leq \constthr[\objfun]\Riemdist{\ptotrsubseqiter}{\ptaccum} + \sum_{\ineqidx\in\ineqset} \abs*{\ineqLagmultaccum[\ineqidx] - \sbra*{\ineqLagmultotrsubseqiter[]}_{\ineqidx}} \opnorm{\Hess\ineqfun[\ineqidx]\paren*{\ptaccum}} + \sum_{\ineqidx\in\ineqset} \abs*{\sbra*{\ineqLagmultotrsubseqiter[]}_{\ineqidx}} \constthr[{\ineqfun[\ineqidx]}]\Riemdist{\ptotrsubseqiter}{\ptaccum},
        \end{align}
        }{
        \begin{align}
            &\opnorm{\Hess[{\pt[]}]\Lagfun\paren*{\allvaraccum} - \partxp[]{\ptaccum}{\ptotrsubseqiter} \circ \Hess[{\pt[]}]\Lagfun\paren*{\allvarotrsubseqiter} \circ \partxp[]{\ptotrsubseqiter}{\ptaccum}} \leq \opnorm{\Hess\objfun\paren*{\ptaccum} - \\
            &\partxp[]{\ptaccum}{\ptotrsubseqiter} \circ \Hess\objfun\paren*{\ptotrsubseqiter} \circ \partxp[]{\ptotrsubseqiter}{\ptaccum}} + \sum_{\ineqidx\in\ineqset} \abs*{\ineqLagmultaccum[\ineqidx] - \sbra*{\ineqLagmultotrsubseqiter[]}_{\ineqidx}} \opnorm{\Hess\ineqfun[\ineqidx]\paren*{\ptaccum}}\\
            &+ \sum_{\ineqidx\in\ineqset}\abs*{\sbra*{\ineqLagmultotrsubseqiter[]}_{\ineqidx}} \opnorm{\Hess\ineqfun[\ineqidx]\paren*{\ptaccum} - \partxp[]{\ptaccum}{\ptotrsubseqiter} \circ \Hess\ineqfun[\ineqidx]\paren*{\ptotrsubseqiter} \circ \partxp[]{\ptotrsubseqiter}{\ptaccum}}\\[-5mm]
            &\leq \constthr[\objfun]\Riemdist{\ptotrsubseqiter}{\ptaccum} + \sum_{\ineqidx\in\ineqset} \abs*{\ineqLagmultaccum[\ineqidx] - \sbra*{\ineqLagmultotrsubseqiter[]}_{\ineqidx}} \opnorm{\Hess\ineqfun[\ineqidx]\paren*{\ptaccum}} + \sum_{\ineqidx\in\ineqset} \abs*{\sbra*{\ineqLagmultotrsubseqiter[]}_{\ineqidx}} \constthr[{\ineqfun[\ineqidx]}]\Riemdist{\ptotrsubseqiter}{\ptaccum},
        \end{align}}
        where the second inequality follows
        \isextendedversion{
        from \cref{ineq:LipschitzHessobjfun,ineq:LipschitzHessineqfun}.
        }{
        from \cref{assu:LipschitzHess}.
        }
        Since the right-hand side converges to zero as $\subseqidx\to\infty$ by the boundedness of $\brc*{\norm*{\ineqLagmultotrsubseqiter[]}}_{\subseqidx}$ and $\ineqLagmultotrsubseqiter[]\to\ineqLagmultaccum[]$, we have
        \begin{equation}\label{eq:diffHessLagpartxpHessLag}
            \lim_{\subseqidx\to\infty}\opnorm{\Hess[{\pt[]}]\Lagfun\paren*{\allvaraccum} - \partxp[]{\ptaccum}{\ptotrsubseqiter} \circ \Hess[{\pt[]}]\Lagfun\paren*{\allvarotrsubseqiter} \circ \partxp[]{\ptotrsubseqiter}{\ptaccum}} = 0.
        \end{equation}

        Using the results above, we now establish the weak second-order stationarity.
        For any $\tanvecone[\ptaccum]\in\weakcriticalcone\paren*{\ptaccum}$ and all $\subseqidx\in\setNz$ sufficiently large,
        we let
        \begin{align}
            \Akn\paren*{\allvarotrsubseqiter}\sbra*{\tanvecone[\ptaccum]}&\coloneqq \metr[\ptaccum]{\partxp[]{\ptaccum}{\ptotrsubseqiter} \circ \Hess[{\pt[]}]\Lagfun\paren*{\allvarotrsubseqiter} \circ \partxp[]{\ptotrsubseqiter}{\ptaccum}\sbra*{\tanvecone[\ptaccum]}}{\tanvecone[\ptaccum]}\\
            &= \metr[\ptotrsubseqiter]{\Hess[{\pt[]}]\Lagfun\paren*{\allvarotrsubseqiter}\sbra*{\partxp[]{\ptotrsubseqiter}{\ptaccum}\sbra*{\tanvecone[\ptaccum]}}}{\partxp[]{\ptotrsubseqiter}{\ptaccum}\sbra*{\tanvecone[\ptaccum]}}\\
            \Bkn\paren*{\allvarotrsubseqiter}\sbra*{\tanvecone[\ptaccum]}&\coloneqq \metr[\ptotrsubseqiter]{\Hess[{\pt[]}]\Lagfun\paren*{\allvarotrsubseqiter}\sbra*{\projwcconeaccum\paren*{\ptotrsubseqiter}\sbra*{\partxp[]{\ptotrsubseqiter}{\ptaccum}\sbra*{\tanvecone[\ptaccum]}}}}{\partxp[]{\ptotrsubseqiter}{\ptaccum}\sbra*{\tanvecone[\ptaccum]}}\\
            \Ckn\paren*{\allvarotrsubseqiter}\sbra*{\tanvecone[\ptaccum]}&\\
            \coloneqq \metr[\ptotrsubseqiter]{\Hess[{\pt[]}]&\Lagfun\paren*{\allvarotrsubseqiter}\sbra*{\projwcconeaccum\paren*{\ptotrsubseqiter}\sbra*{\partxp[]{\ptotrsubseqiter}{\ptaccum}\sbra*{\tanvecone[\ptaccum]}}}}{\projwcconeaccum\paren*{\ptotrsubseqiter}\sbra*{\partxp[]{\ptotrsubseqiter}{\ptaccum}\sbra*{\tanvecone[\ptaccum]}}}.
        \end{align}
        Then, it follows that 
        \begin{align}
            &\label{ineq:HessLagwSOSPbound}\\
            &\metr[\ptaccum]{\Hess[{\pt[]}]\Lagfun\paren*{\allvaraccum}\sbra*{\tanvecone[\ptaccum]}}{\tanvecone[\ptaccum]}= \metr[\ptaccum]{\Hess[{\pt[]}]\Lagfun\paren*{\allvaraccum}\sbra*{\tanvecone[\ptaccum]}}{\tanvecone[\ptaccum]} - \Akn\paren*{\allvarotrsubseqiter}\sbra*{\tanvecone[\ptaccum]}
            + \Akn\paren*{\allvarotrsubseqiter}\sbra*{\tanvecone[\ptaccum]}\\
            &\quad- \Bkn\paren*{\allvarotrsubseqiter}\sbra*{\tanvecone[\ptaccum]}
            + \Bkn\paren*{\allvarotrsubseqiter}\sbra*{\tanvecone[\ptaccum]}
            - \Ckn\paren*{\allvarotrsubseqiter}\sbra*{\tanvecone[\ptaccum]}
            + \Ckn\paren*{\allvarotrsubseqiter}\sbra*{\tanvecone[\ptaccum]}\\
            &\geq - \opnorm{\Hess[{\pt[]}]\Lagfun\paren*{\allvaraccum} - \partxp[]{\ptaccum}{\ptotrsubseqiter} \circ \Hess[{\pt[]}]\Lagfun\paren*{\allvarotrsubseqiter} \circ \partxp[]{\ptotrsubseqiter}{\ptaccum}}\Riemnorm[\ptaccum]{\tanvecone[\ptaccum]}^{2}\\
            &\quad-\opnorm{\Hess[{\pt[]}]\Lagfun\paren*{\allvarotrsubseqiter}}\Riemnorm[\ptotrsubseqiter]{\partxp[]{\ptotrsubseqiter}{\ptaccum}\sbra*{\tanvecone[\ptaccum]} - \projwcconeaccum\paren*{\ptotrsubseqiter}\sbra*{\partxp[]{\ptotrsubseqiter}{\ptaccum}\sbra*{\tanvecone[\ptaccum]}}}\\
            &\qquad \paren*{\Riemnorm[\ptotrsubseqiter]{\partxp[]{\ptotrsubseqiter}{\ptaccum}\sbra*{\tanvecone[\ptaccum]}} + \Riemnorm[\ptotrsubseqiter]{\projwcconeaccum\paren*{\ptotrsubseqiter}\sbra*{\partxp[]{\ptotrsubseqiter}{\ptaccum}\sbra*{\tanvecone[\ptaccum]}}}}\\
            &\quad+ \metr[\ptotrsubseqiter]{\Hess[{\pt[]}]\Lagfun\paren*{\allvarotrsubseqiter}\sbra*{\projwcconeaccum\paren*{\ptotrsubseqiter}\sbra*{\partxp[]{\ptotrsubseqiter}{\ptaccum}\sbra*{\tanvecone[\ptaccum]}}}}{\projwcconeaccum\paren*{\ptotrsubseqiter}\sbra*{\partxp[]{\ptotrsubseqiter}{\ptaccum}\sbra*{\tanvecone[\ptaccum]}}}.
        \end{align}
        By \cref{eq:limdiffpartxpprojpartxp,eq:diffHessLagpartxpHessLag} and the boundedness of $\brc*{\Riemnorm[\ptotrsubseqiter]{\projwcconeaccum\paren*{\ptotrsubseqiter}\sbra*{\partxp[]{\ptotrsubseqiter}{\ptaccum}\sbra*{\tanvecone[\ptaccum]}}}}_{\subseqidx}$,\\
        $\brc*{\Riemnorm[\ptotrsubseqiter]{\partxp[]{\ptotrsubseqiter}{\ptaccum}\sbra*{\tanvecone[\ptaccum]}}}_{\subseqidx}$, and $\brc*{\opnorm{\Hess\Lagfun\paren*{\allvarotrsubseqiter}}}_{\subseqidx}$ around $\allvaraccum$, the first two terms on the right-hand side of \cref{ineq:HessLagwSOSPbound} converge to zero as $\subseqidx\to\infty$.
        In addition, since the sequence $\brc*{\projwcconeaccum\paren*{\ptotrsubseqiter}\sbra*{\partxp[]{\ptotrsubseqiter}{\ptaccum}\sbra*{\tanvecone[\ptaccum]}}}_{\subseqidx}$ satisfies \cref{eq:asymptoticwctanvec} by definition, the last term on the right-hand side accumulates at a nonnegative value by taking a subsequence that realizes the limit inferior in \cref{eq:liminfHessLagfunotriter}.
        Therefore, \eqcref{ineq:HessLagwSOSPbound} implies that $\allvaraccum$ is a \wSOSP{}, which completes the proof.
    \end{proof}

    \section{Computation of dual variables}\label{subsec:computdualvar}
        In this appendix, we consider the clipping for the update of the dual variables in \cref{algo:RIPTRMInner} since the simple update for the dual variables $\ineqLagmultinriter[] + \dirineqLagmultinriter[]$ does not always preserve positivity, \cref{assu:ineqLagmultbounded}, and \cref{assu:ineqLagmultconverged}.
        The clipping is originally introduced in \cite[Section 4.3]{Connetal2000PrimalDualTRAlgoforNonconvexNLP} for the Euclidean \IPTRM{}.
        Here, we discuss its extension to the Riemannian setting.
        
        Letting $\constfou, \constfiv \in \setR$ be constants satisfying $0 < \constfou < 1 < \constfiv$, we define the interval for the $\ineqidx$-th element as $\clipintvlinriter[\ineqidx] \coloneqq \sbra*{\clipintvlinritermin[\ineqidx], \clipintvlinritermax[\ineqidx]}$, where
        \begin{align}\label{def:clipintvlminmax}
            \clipintvlinritermin[\ineqidx] \coloneqq \constfou\min\brc*{1, \ineqLagmultinriter[\ineqidx], \frac{\barrparam[]}{\pullback[\ptinriter]{\ineqfun[\ineqidx]}\paren*{\dirptinriter}}}, \quad \clipintvlinritermax[\ineqidx] \coloneqq \max\brc*{\constfiv, \ineqLagmultinriter[\ineqidx], \frac{\constfiv}{\barrparam[]}, \frac{\constfiv}{\pullback[\ptinriter]{\ineqfun[\ineqidx]}\paren*{\dirptinriter}}}
        \end{align}
        for each $\ineqidx\in\ineqset$. 
        We update the dual variable as follows:
        \begin{align}\label{def:clippingupdate}
            \ineqLagmultinriterp[] \leftarrow 
            \begin{cases}
                \clip[{\clipintvl[\inriteridx]}]\paren*{\ineqLagmultinriter[] + \dirineqLagmultinriter[]} & \text{if } \ptinriterp=\retr[\ptinriter]\paren*{\dirptinriter}, \\
               \ineqLagmultinriter[] & \text{if } \ptinriterp=\ptinriter,\\
            \end{cases}
        \end{align}
        where $\clip[{\clipintvl[\inriteridx]}]\colon\setR[\ineqdime]\to\setR[\ineqdime]$ is the operator whose $\ineqidx$-th component is defined by 
        \begin{align}
            \sbra*{\clip[{\clipintvl[\inriteridx]}]\paren*{\vecfou}}_{\ineqidx} = \max\brc*{\clipintvlinritermin[\ineqidx], \min\brc*{\clipintvlinritermax[\ineqidx], \vecfou_{\ineqidx}}}.
        \end{align}
        In \cref{sec:experiment}, we observe that the choices $\constfou=0.5$ and $\constfiv=10^{20}$ work satisfactorily in practice.

        Now, we provide theoretical analyses for \cref{def:clippingupdate}.
        Note that the update \cref{def:clippingupdate} satisfies the positivity of the dual variable by definition.
        We prove that the strategy satisfies \cref{assu:ineqLagmultbounded,assu:ineqLagmultconverged} in the following theorem: recall that $\brc*{\inrsuccsubseqiter}_{\succsubseqidx}\subseteq\setNz$ is the sequence of successful iterations.
        \begin{theorem}\label{lemm:dualvarsuffcond}
           Let $\brc*{\paren*{\ptinriter, \ineqLagmultinriter[]}}_{\inriteridx}$ be the sequence generated by \cref{algo:RIPTRMInner} with $\brc*{\ineqLagmultinriter[]}_{\inriteridx}$ updated by \cref{def:clippingupdate}.
           Suppose \cref{assu:Cauchydecreasing}, \assuenumicref{assu:ineqfunbounded}{assu:ineqfunboundedvalue}, and \assuenumicref{assu:objfunbounded}{assu:objfunboundedbelow}.
           Then, the following hold: 
           \begin{enumerate}
               \item \cref{assu:ineqLagmultbounded} holds. \label{lemm:dualvarsuffcondbounded}
               \item Under \assuenumicref{assu:ineqfunbounded}{assu:ineqfunboundedgradhess} and \cref{assu:pullbackineqfuncont}, if 
                   \begin{align}\label{eq:limdirptinriter}
                       \lim_{\succsubseqidx\to\infty}\Riemnorm[\ptinrsuccsubseqiter]{\dirinrsuccsubseqiter} = 0,
                   \end{align}
                   then \cref{assu:ineqLagmultconverged} holds.\label{lemm:dualvarsuffcondconverged}
           \end{enumerate}
        \end{theorem}
        \begin{proof}[Proof]
            First, we consider \cref{lemm:dualvarsuffcondbounded}.
            Since the dual variable is updated when the primal iterate is successful, we focus on the successful iterations.
            From \enumicref{lemm:ineqfunbound}{lemm:ienqfunptlowbd}, there exists $\epsfiv > 0$ such that $\frac{1}{\epsfiv} \geq \frac{1}{\ineqfun[\ineqidx]\paren*{\ptinriter}}$ for any $\ineqidx\in\ineqset$ and all $\inriteridx\in\setNz$, which, together with \cref{def:clipintvlminmax}, implies that 
            \begin{align}
                \ineqLagmultinriter[\ineqidx] \leq \max\brc*{\constfiv, \ineqLagmult[\ineqidx]^{0}, \frac{\constfiv}{\barrparam[]}, \frac{\constfiv}{\epsfiv}}
            \end{align}
            for each $\ineqidx\in\ineqset$.
            Thus, \cref{lemm:dualvarsuffcondbounded} holds.

            Next, we consider \cref{lemm:dualvarsuffcondconverged}.
            Recall that $\succset$ denotes the set of the successful iterations and $\brc*{\ptinrsuccsubseqiter}_{\succsubseqidx}$ is the ordered sequence of successful iterates.
            If $\card[\succset]$ is finite, then it follows from \cref{prop:existsuccidx} that $\barrKKTvecfld\paren*{\allvarinriter;\barrparam[]} = 0$ for all $\inriteridx\in\setNz$ sufficiently large.
            Thus, we have $\norm{\ineqLagmultinriter[]-\barrparam[]\inv{\Ineqfunmat[]\paren*{\ptinriter}}\onevec}=0$ for all $\inriteridx\in\setNz$ sufficiently large, meaning that \cref{assu:ineqLagmultconverged} holds.
            In the following, we consider the case where $\card[\succset]$ is infinite.
            It follows from $\pullback[\ptinrsuccsubseqiter]{\ineqfun[\ineqidx]}\paren*{\dirinrsuccsubseqiter} = \ineqfun[\ineqidx]\paren*{\pt^{\inrsuccsubseqiter + 1}}$,
            \cref{eq:pullbackineqfuncont},
            and \cref{eq:limdirptinriter} that
            \begin{align}\label{eq:limabsdiffineqfun}
                \lim_{\succsubseqidx\to\infty} \abs*{\ineqfun[\ineqidx]\paren*{\ptinrsuccsubseqiter} - \ineqfun[\ineqidx]\paren*{\pt^{\inrsuccsubseqiter + 1}}} = 0
            \end{align}
            for all $\ineqidx\in\ineqset$.
            Therefore, we have
            \begin{align}
                &\norm*{\ineqLagmultinrsuccsubseqiter[] + \dirineqLagmultinrsuccsubseqiter - \barrparam[]\inv{\Ineqfunmat[]\paren*{\pt^{\inrsuccsubseqiter + 1}}}\onevec}\\ &\leq \norm*{\ineqLagmultinrsuccsubseqiter[] + \dirineqLagmultinrsuccsubseqiter - \barrparam[]\inv{\Ineqfunmat[]\paren*{\pt^{\inrsuccsubseqiter}}}\onevec} + \barrparam[]\norm*{\inv{\Ineqfunmat[]\paren*{\ptinrsuccsubseqiter}}\onevec -  \inv{\Ineqfunmat[]\paren*{\pt^{\inrsuccsubseqiter+1}}}\onevec}\\
                &\leq \sum_{\ineqidx\in\ineqset}\abs*{\frac{\ineqLagmultinrsuccsubseqiter}{\ineqfun[\ineqidx]\paren*{\ptinrsuccsubseqiter}}\metr[\ptinrsuccsubseqiter]{\gradstr\ineqfun[\ineqidx]\paren*{\ptinrsuccsubseqiter}}{\dirinrsuccsubseqiter}} + \sum_{\ineqidx\in\ineqset}\frac{\barrparam[]\abs*{\ineqfun[\ineqidx]\paren*{\ptinrsuccsubseqiter} - \ineqfun[\ineqidx]\paren*{\pt^{\inrsuccsubseqiter+1}}}}{\ineqfun[\ineqidx]\paren*{\ptinrsuccsubseqiter}\ineqfun[\ineqidx]\paren*{\pt^{\inrsuccsubseqiter+1}}}\\
                &\leq \sum_{\ineqidx\in\ineqset}\frac{\abs*{\ineqLagmultinrsuccsubseqiter}}{\epsfiv}\Riemnorm[\ptinrsuccsubseqiter]{\gradstr\ineqfun[\ineqidx]\paren*{\ptinrsuccsubseqiter}}\Riemnorm[\ptinrsuccsubseqiter]{\dirinrsuccsubseqiter} + \sum_{\ineqidx\in\ineqset}\frac{\barrparam[]\abs*{\ineqfun[\ineqidx]\paren*{\ptinrsuccsubseqiter} - \ineqfun[\ineqidx]\paren*{\pt^{\inrsuccsubseqiter+1}}}}{{\epsfiv}^{2}},
            \end{align}
            where the second inequality follows from \cref{eq:dualNewtoneq} and the third one from \enumicref{lemm:ineqfunbound}{lemm:ienqfunptlowbd}.
            By \cref{lemm:dualvarsuffcondbounded}, \assuenumicref{assu:ineqfunbounded}{assu:ineqfunboundedgradhess}, \cref{eq:limabsdiffineqfun}, and \cref{eq:limdirptinriter}, the right-hand side converges to zero as $\subseqidx\to\infty$.
            Thus, under \assuenumicref{assu:ineqfunbounded}{assu:ineqfunboundedvalue}, letting $\tholdvalsix > 0$ be the scalar satisfying $\ineqfun[\ineqidx]\paren*{\ptinrsuccsubseqiter} \leq \tholdvalsix$ for any $\ineqidx\in\ineqset$, we obtain
            \begin{align}
                \abs*{\ineqLagmultinrsuccsubseqiter[\ineqidx] + \sbra*{\dirineqLagmultinrsuccsubseqiter}_{\ineqidx} - \frac{\barrparam[]}{\ineqfun[\ineqidx]\paren*{\pt^{\inrsuccsubseqiter + 1}}}\onevec} \leq \min\brc*{1-\constfou, \constfiv - 1} \frac{\barrparam[]}{\tholdvalsix}
            \end{align}
            for any $\ineqidx\in\ineqset$ and all $\succsubseqidx\in\setNz$ sufficiently large, which implies
            \begin{align}
                &\ineqLagmultinrsuccsubseqiter[\ineqidx] + \sbra*{\dirineqLagmultinrsuccsubseqiter}_{\ineqidx}\\
                &\quad \geq -\paren*{1-\constfou}\frac{\barrparam[]}{\tholdvalsix} + \frac{\barrparam[]}{\ineqfun[\ineqidx]\paren*{\pt^{\inrsuccsubseqiter + 1}}} \geq -\paren*{1-\constfou}\frac{\barrparam[]}{\ineqfun[\ineqidx]\paren*{\pt^{\inrsuccsubseqiter + 1}}} + \frac{\barrparam[]}{\ineqfun[\ineqidx]\paren*{\pt^{\inrsuccsubseqiter + 1}}} = \frac{\constfou\barrparam[]}{\ineqfun[\ineqidx]\paren*{\pt^{\inrsuccsubseqiter + 1}}},\\
                &\ineqLagmultinrsuccsubseqiter[\ineqidx] + \sbra*{\dirineqLagmultinrsuccsubseqiter}_{\ineqidx}  \leq \paren*{\constfiv - 1}\frac{\barrparam[]}{\tholdvalsix} + \frac{\barrparam[]}{\ineqfun[\ineqidx]\paren*{\pt^{\inrsuccsubseqiter + 1}}} \leq \paren*{\constfiv - 1}\frac{\barrparam[]}{\ineqfun[\ineqidx]\paren*{\pt^{\inrsuccsubseqiter + 1}}} + \frac{\barrparam[]}{\ineqfun[\ineqidx]\paren*{\pt^{\inrsuccsubseqiter + 1}}} = \frac{\constfiv\barrparam[]}{\ineqfun[\ineqidx]\paren*{\pt^{\inrsuccsubseqiter + 1}}}.
            \end{align}
            for any $\ineqidx\in\ineqset$.
            By the update rule \cref{def:clippingupdate}, $\ineqLagmult[]^{\inrsuccsubseqiter + 1} = \ineqLagmultinrsuccsubseqiter[] + \dirineqLagmultinrsuccsubseqiter$ holds for all $\subseqidx\in\setNz$ sufficiently large, which implies
            \begin{align}\label{eq:prodineqfunmatineqLagmult}
                \Ineqfunmat[]\paren*{\pt^{\inrsuccsubseqiter + 1}}\ineqLagmult[]^{\inrsuccsubseqiter + 1} = \Ineqfunmat[]\paren*{\pt^{\inrsuccsubseqiter + 1}}\inv{\Ineqfunmat[]\paren*{\ptinrsuccsubseqiter}}\paren*{\barrparam[]\onevec - \IneqLagmultmat[]^{\inrsuccsubseqiter}\coineqgradopr[\ptinrsuccsubseqiter]\sbra*{\dir[\ptinrsuccsubseqiter]}},
            \end{align}
            where we write $\IneqLagmultmat[]^{\inrsuccsubseqiter}$ for $\diag\paren*{\ineqLagmult[]^{\inrsuccsubseqiter}}$ and the equality follows from \cref{eq:dualNewtoneq}.
            Here, it follows from \enumicref{lemm:ineqfunbound}{lemm:ienqfunptlowbd} that 
            \begin{align}
                \norm*{\Ineqfunmat[]\paren*{\pt^{\inrsuccsubseqiter + 1}}\inv{\Ineqfunmat[]\paren*{\pt^{\inrsuccsubseqiter}}} - \eyemat[\ineqdime]} \leq \sum_{\ineqidx\in\ineqset} \abs*{\frac{\ineqfun[\ineqidx]\paren*{\pt^{\inrsuccsubseqiter + 1}}}{\ineqfun[\ineqidx]\paren*{\ptinrsuccsubseqiter}} - 1} \leq \sum_{\ineqidx\in\ineqset} \frac{\abs*{\ineqfun[\ineqidx]\paren*{\pt^{\inrsuccsubseqiter + 1}} - \ineqfun[\ineqidx]\paren*{\ptinrsuccsubseqiter}}}{\epsfiv},
            \end{align}
            where $\eyemat[\ineqdime]\in\setR[\ineqdime\times\ineqdime]$ is the identity matrix, and the right-hand side converges to zero as $\inrsuccsubseqiter\to\infty$ by \cref{eq:limabsdiffineqfun}.
            Therefore, \eqscref{eq:prodineqfunmatineqLagmult,eq:limdirptinriter} and the boundedness of $\brc*{\ineqLagmultinrsuccsubseqiter[]}_{\succsubseqidx}$ and $\brc*{\Riemnorm[\ptinrsuccsubseqiter]{\gradstr\ineqfun[\ineqidx]\paren*{\ptinrsuccsubseqiter}}}_{\succsubseqidx,\ineqidx}$ imply 
            \begin{align}\label{eq:limitprodIneqfunmatineqLagmult}
                \lim_{\succsubseqidx\to\infty}\Ineqfunmat[]\paren*{\pt^{\inrsuccsubseqiter + 1}}\ineqLagmult[]^{\inrsuccsubseqiter + 1} = \barrparam[]\onevec.
            \end{align}
            Since $\allvar[]^{\inrsuccsubseqiter + 1}=\allvar^{\inrsuccsubseqiterp}$ holds  and all iterates between $\inrsuccsubseqiter$ and $\inrsuccsubseqiterp$ are unsuccessful for all $\succsubseqidx\in\setNz$, \eqcref{eq:limitprodIneqfunmatineqLagmult} extends to all the iterates; that is, it holds that $\lim_{\inriteridx\to\infty}\Ineqfunmat[]\paren*{\ptinriter}\ineqLagmultinriter[] = \barrparam[]\onevec$.
            From $\ptinriter\in\strictfeasirgn$, for all $\inriteridx\in\setNz$, we have $\lim_{\inriteridx\to\infty}\norm{\ineqLagmultinriter[]-\barrparam[]\inv{\Ineqfunmat[]\paren*{\ptinriter}}\onevec}=0$.
            The proof is complete.
        \end{proof}

        We also provide a sufficient condition for \cref{eq:limdirptinriter}.
        To this end, we first introduce the following lemma:
        \begin{lemma}[{\cite[Lemma~2]{HuangAbsilGallivan2015RiemSymRankOneTRM}}]\label{lemm:retrtanvecequiv}
            Let $\mani$ be a Riemannian manifold endowed with a retraction $\retr[]$ and let $\pt[1]\in\mani$.
            There exist positive scalars $\constsix[1], \constsix[2], \tholdvalsev[{\constsix[1], \constsix[2]}] \in \setRpp[]$ such that, for all $\pt[2]$ in a sufficiently small neighborhood of $\pt[1]$ and all $\tanvecone[{\pt[2]}], \tanvectwo[{\pt[2]}] \in \tanspc[{\pt[2]}]\mani$ with $\Riemnorm[{\pt[2]}]{\tanvecone[{\pt[2]}]} \leq \tholdvalsev[{\constsix[1], \constsix[2]}]$ and $\Riemnorm[{\pt[2]}]{\tanvectwo[{\pt[2]}]} \leq \tholdvalsev[{\constsix[1], \constsix[2]}]$,
            \begin{align}
                \constsix[1] \Riemnorm[{\pt[2]}]{\tanvecone[{\pt[2]}] - \tanvectwo[{\pt[2]}]} \leq \Riemdist{\retr[{\pt[2]}]\paren*{\tanvecone[{\pt[2]}]}}{\retr[{\pt[2]}]\paren*{\tanvectwo[{\pt[2]}]}} \leq \constsix[2] \Riemnorm[{\pt[2]}]{\tanvecone[{\pt[2]}] - \tanvectwo[{\pt[2]}]}. 
            \end{align}
        \end{lemma}
        Then, the following sufficient condition holds:
        \begin{remark}
            If the maximal trust region radius $\maxtrradius > 0$ is sufficiently small and the sequence $\brc*{\ptinrsuccsubseqiter}_{\succsubseqidx}$ converges to some $\ptaccum\in\mani$, then \cref{eq:limdirptinriter} holds.
            Indeed, by \cref{lemm:retrtanvecequiv}, there exists $\constsix[] > 0$ such that, for all $\inrsuccsubseqiter\in\setNz$ sufficiently large, 
            \begin{align}\label{ineq:normdirptinrsuccsubseqiterdistpts}
                \constsix[] \Riemnorm[{\ptinrsuccsubseqiter}]{\dir[\ptinrsuccsubseqiter]} \leq \Riemdist{\retr[{\ptinrsuccsubseqiter}]\paren*{\dir[\ptinrsuccsubseqiter]}}{\ptinrsuccsubseqiter} = \Riemdist{\pt^{1 + \inrsuccsubseqiter}}{\ptinrsuccsubseqiter}
            \end{align}
            since the point $\ptinrsuccsubseqiter\in\strictfeasirgn$ is sufficiently close to $\ptaccum$ and the search direction $\dir[\ptinrsuccsubseqiter]\in\tanspc[\ptinrsuccsubseqiter]\mani$ satisfies $\Riemnorm[{\ptinrsuccsubseqiter}]{\dir[\ptinrsuccsubseqiter]} \leq \maxtrradius$.
            By the assumption, $\brc*{\ptinrsuccsubseqiter}_{\succsubseqidx}$ is a Cauchy sequence.
            Hence, \eqcref{ineq:normdirptinrsuccsubseqiterdistpts} implies \cref{eq:limdirptinriter} by taking the limit as $\inrsuccsubseqiter\to\infty$.
        \end{remark}


    \section{Further experiments}\label{sec:additionalexperiment}
        In this appendix, we present an additional experiment on the minimization of the Rosenbrock function on the Grassmann manifold.
        The experimental environment and setting are the same as those in \cref{sec:experiment}.

    \subsection{Problem settings: Rosenbrock function minimization}
        The Rosenbrock function is a well-known benchmark in the field of mathematical optimization.
        We formulate the following optimization problem on the \Grassmannian{}:
        \begin{mini!}[2]
            {\matone \in \Grassmann\paren*{\matrowdime, \matcoldime}}{\sum_{\vecidx=1}^{\matrowdime\matcoldime-1} \Rosenbrockcoeff\paren*{\sbra*{\vectorize\paren*{\matone}}_{\vecidx+1} - \sbra*{\vectorize\paren*{\matone}}_{\vecidx}}^{2} + \paren*{1 - \sbra*{\vectorize\paren*{\matone}}_{\vecidx}}^{2} }{\label{objfun:Rosenbrock}}
            {\label{Prob:illcondRosenbrock}}{}
            \addConstraint{\matone}{ \geq \Rosenbrockcoefftwo \cdot \onevec\trsp{\onevec},}
        \end{mini!}
        where $\vectorize\colon\setR[\matrowdime\times\matcoldime]\to\setR[\matrowdime\matcoldime]$ is the vectorization operator that rearranges the rows of a matrix into a vector,
        and $\Grassmann\paren*{\matrowdime, \matcoldime}\coloneqq\brc*{\spanstr\paren*{\matone}\relmiddle{\colon}\matone\in\setR[\matrowdime\times\matcoldime], \trsp{\matone}\matone=\eyemat[\matcoldime]}$ is the \Grassmannian{}. 
    
        \textbf{Input.} 
        We consider the case $\paren*{\matrowdime, \matcoldime} = \paren*{5, 3}$, $\Rosenbrockcoeff = 10^{7}$, and $\Rosenbrockcoefftwo = -0.01$.
        The initial point is set to $\trsp{\sbra*{\eyemat[\matcoldime] \relmiddle{|} \bf{0}}}\in\setR[\matrowdime\times\matcoldime]$.
        We initialize the dual variables with the vector of ones.

    \subsection{Results and discussion}

        We applied the algorithms to the artificial instance of \cref{Prob:illcondRosenbrock}.
        \begin{figure}[t]
            \centering
            \begin{minipage}{0.49\linewidth}
                \centering
                \includegraphics[width=\linewidth]{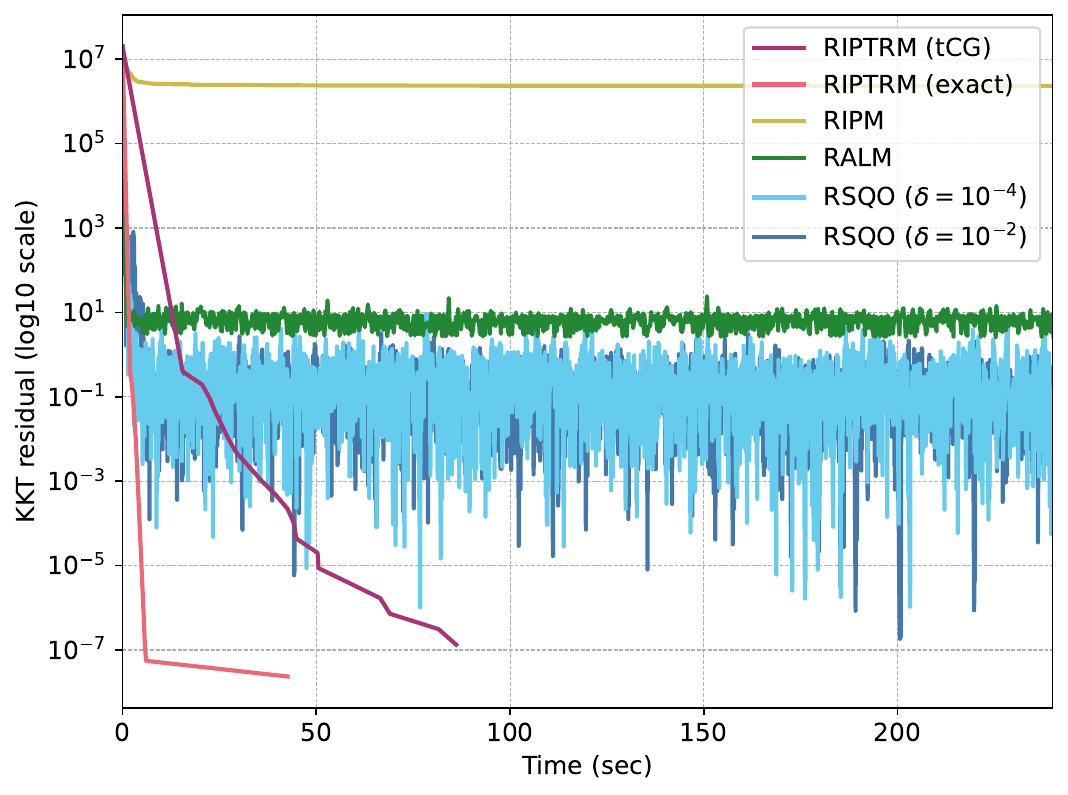}
                \caption{Residual over time for Rosenbrock function minimization}
                \label{fig:1stordresidRosenbrock}
            \end{minipage}
            \hfill
            \begin{minipage}{0.49\linewidth}
                \centering
                \includegraphics[width=\linewidth]{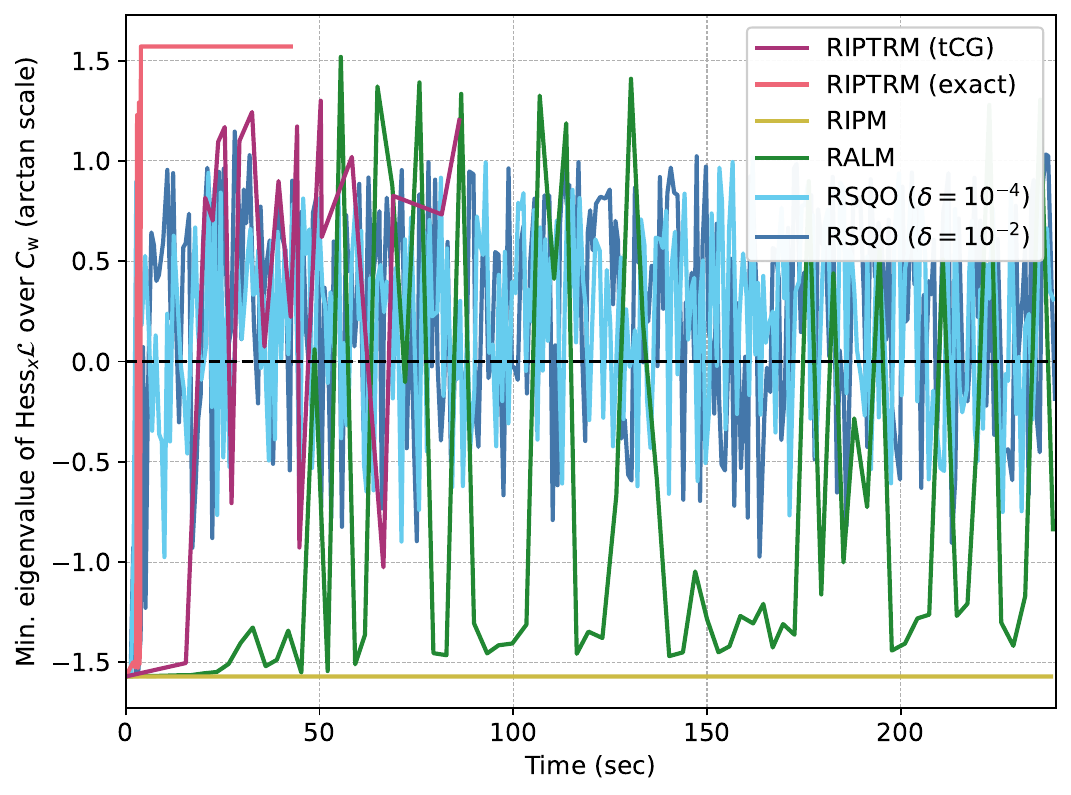}
                \caption{Second-order stationarity over time for Rosenbrock function minimization}
                \label{fig:2ndordresidRosenbrock}
            \end{minipage}
        \end{figure}
        \cref{fig:1stordresidRosenbrock} shows the residuals of the algorithms over time.
        We see that \RIPTRM{} (\Exact{}) successfully solved the instance with the fastest speed and the highest accuracy.
        In addition, we also measured the second-order stationarity \cref{eq:weaksecondordernecessarycond} in the following manner; the $\ineqidx$-th inequality constraint is regarded as active if it satisfies $\ineqfun[\ineqidx]\paren{\matone} < 10^{-6}$.
        We then identify $\weakcriticalcone$ \cref{{def:weakcriticalcone}} at $\matone$ and compute the minimum eigenvalue of the Hessian of the Lagrangian over $\weakcriticalcone\paren*{\matone}$ as the second-order stationarity.
        Note that, for this instance, the minimum eigenvalue at the initial point is $-2.000\times 10^{7}$, indicating that the Hessian of the Lagrangian has the large negative eigenvalue in this problem.
        \cref{fig:2ndordresidRosenbrock} illustrates the second-order stationarity using an arc-tangent scale.
        For clarity, we omit the inner iterations of \RIPTRM{}s and plot only their outer iterations.
        For \RIPM{}, \RSQO{}, and \RALM{}, values are plotted every 15 iterations to improve readability.
        The black dashed line in the figure represents zero, indicating that values above the line correspond to a positive minimum eigenvalue and signify that the second-order stationarity condition is satisfied.
        We observe that \RIPTRM{} (\Exact{}) rapidly achieved and maintained the minimum eigenvalues above the line, indicating successful computation of an \SOSP{}.
        This may suggest that the \exactsolution{}s are robust and efficient for problem instances where the Hessian of the Lagrangian has a large negative eigenvalue, as \exactsolution{}s can incorporate the information of the negative eigenvalues.
        In contrast, the second-order stationarity of \RIPTRM{} (\tCG{}), \RALM{}, and \RSQO{} oscillated between positive and negative ranges, and that of \RIPM{} stays negative.
        This behavior may be attributed to the lack of second-order convergence properties in these algorithms.
        

    }{}

\end{document}